\documentclass[a4paper,12pt,reqno,twoside,onecolumn,final]{amsart}

\usepackage{amsmath,amsfonts,amssymb,amscd}

\usepackage[foot]{amsaddr} 

\usepackage{graphicx}
\usepackage{epstopdf}

\usepackage{bm}  

\usepackage[utf8]{inputenc} 
\usepackage{lmodern}
\usepackage{subfig}
\usepackage{natbib}

\usepackage[]{algorithm}
\usepackage[]{algorithmic}
\usepackage{amssymb}
\usepackage{mathabx}

\usepackage{color}

\usepackage[skip=5pt,indent=15pt]{parskip}

\usepackage[width=.99\textwidth]{caption}

\usepackage{booktabs}            
\usepackage{array}               
\usepackage{enumitem}            
\usepackage{verbatim}            
\usepackage{subfig}
\usepackage{soul}
\usepackage{natbib}
\usepackage{hyperref}            
\usepackage{multirow}            
\usepackage[mathscr]{euscript}
\usepackage{color}
\usepackage{xcolor}
\usepackage{graphicx}
\usepackage{float}               
\usepackage{indentfirst}
\usepackage{placeins}            
\usepackage[normalem]{ulem}

\usepackage{setspace}

\makeatletter

\usepackage[pagewise]{lineno}

\newtheorem{theorem}{Theorem}[section]
\newtheorem{lemma}[theorem]{Lemma}
\newtheorem{corollary}[theorem]{Corollary}
\newtheorem{proposition}[theorem]{Proposition}
\newtheorem{remark}[theorem]{Remark}

\numberwithin{equation}{section}

\definecolor{ok}{rgb}{0.0, 0.0, 0.0}

\definecolor{giovanni}{rgb}{0.5, 0.1, 0.7}
\definecolor{maicon}{rgb}{0.8,0.7,0.}
\definecolor{maicon-red}{rgb}{1.00, 0.35, 0.21}

\def\0{\phantom{0}}
\def\ds{\displaystyle}

\def\nd{d}

\def\diamOm{\rho_{\Omega}}

\newcommand{\ahdp}{\mathfrak{a}}

\let\vector\relax
\newcommand{\vector}[1]{{\boldsymbol {#1}}}
\newcommand{\tensor}[1]{{\boldsymbol {#1}}}

\newcommand{\piola}[1]{\boldsymbol{P}_{F_K}{#1}}
\newcommand{\proj}[1]{{\pi_h^{#1}}}

\newcommand{\vertiii}[1]{{\vert\kern-0.25ex\vert\kern-0.25ex\vert #1 \vert\kern-0.25ex\vert\kern-0.25ex\vert}}

\newcommand\Tstrut{\rule{0pt}{2.5ex}}         
\newcommand\Bstrut{\rule[-1.2ex]{0pt}{0pt}}   

\newcommand{\TauH}{\mathcal{T}_h}
\let\div\relax
\DeclareMathOperator{\div}{\mathrm{div}}
\DeclareMathOperator{\des}{\mathrm{ds}}

\DeclareMathOperator{\dex}{\mathrm{dx}}

\DeclareMathOperator{\tr}{tr}

\DeclareMathOperator{\hdp}{\scriptscriptstyle \mathrm{HDP}}
\DeclareMathOperator{\ph}{\scriptscriptstyle \mathrm{PH}}
\DeclareMathOperator{\ap}{\scriptscriptstyle \mathrm{AP}}

\begin{document}


\title[Primal hybrid method for linear elasticity with stress recovery]{Numerical analysis of a locking-free primal hybrid method for linear elasticity with
$H(\mathrm{div})$-conforming stress recovery}



\author[G. Taraschi]{Giovanni Taraschi}
\address[G. Taraschi]{Universidade Estadual de Campinas (Unicamp), Center for Energy and Petroleum Studies, Brazil}
\curraddr{}
\email[G. Taraschi]{taraschi@unicamp.br}
\thanks{}

\author[M.R. Correa]{Maicon R. Correa}
\address[M.R. Correa]{Department of Applied Mathematics, Institute of Mathematics, Statistics and Computer Science, University of São Paulo, Brazil}
\email[M.R. Correa]{maicon@ime.usp.br}
\thanks{}

\subjclass[2020]{65N12, 65N22, 65N30, 35A35}

\date{}

\dedicatory{
{\rm This is a preprint of the original article published in} \\{\it  Journal of Scientific Computing, v.110, article 9 (2027) } \\
\rm 
https://doi.org/10.1007/s10915-026-03524-0
}

\begin{abstract}

In this work, we study a primal hybrid finite element method for the approximation of linear elasticity problems, posed in terms of displacement, an auxiliary pressure
field, and a Lagrange multiplier related to the traction.
We develop a general analysis for the existence and uniqueness of the solution for the discrete problem, which is applied to the construction of stable approximation spaces on triangular and quadrilateral meshes.
The use of these spaces lead to optimal convergence orders, resulting in a locking-free method capable of providing robust approximations
for nearly incompressible problems.
Finally, we propose a strategy for recovering the stress field from the hybrid solution by solving element-wise sub-problems.
The resulting stress approximation is $H(\mathrm{div})$-conforming, locally equilibrated, weakly symmetric, and robust to locking.
\end{abstract}


\maketitle


\paragraph{\em Keywords}
Primal Hybrid Finite Element Method, Linear elasticity, Stress recovery, Locking-free approximations

\section{Introduction}

Let $\Omega \subset \mathbb{R}^\nd$, $\nd=\{2,3 \}$, be a simply connected polyhedral domain occupied by an elastic body.
We assume that its boundary $\partial \Omega$ is partitioned into the surfaces/curves $\Gamma_D$ and $\Gamma_N$ that satisfy
\( \partial \Omega = \overline{\Gamma_D \cup \Gamma_N}, \quad \Gamma_D \cap \Gamma_N = \emptyset, \quad \textrm{ and } \quad \vert \Gamma_D \vert > 0,
\)
where $\vert \Gamma_D \vert$ denotes the Lebesgue measure of $\Gamma_D$.
Denoting by $\mathbb{M} = \mathbb{R}^{\nd \times \nd}$ the space of second-order real tensors and by $\mathbb{S} = \mathbb{R}^{\nd \times \nd}_\mathrm{sym}$ the subspace of symmetric tensors, the linear elasticity problem \citep{arnold1990mixed} consists in finding the displacement field $\vector{u}: \Omega \rightarrow \mathbb{R}^\nd$ and the stress tensor $\tensor{\sigma}: \Omega \rightarrow \mathbb{S}$ such that
\begin{subequations} \label{eq:model}
    \begin{align}
        \tensor{\sigma} & = \tensor{C} \tensor{\varepsilon}(\vector{u}) & \textrm{ in } \Omega, \label{eq:model_constitutive} \\
        {\color{ok}-\div \tensor{\sigma}} & = \vector{f} & \textrm{ in } \Omega, \label{eq:model_equilibrium} \\
        \vector{u} & = \vector{u}_D & \textrm{ on } \Gamma_D, \label{eq:model_dirichlet} \\
        \tensor{\sigma} \vector{n} & = \vector{t}_N & \textrm{ on } \Gamma_N. \label{eq:model_neumann}
    \end{align}
\end{subequations}

In the linear constitutive equation \eqref{eq:model_constitutive}, $\tensor{\varepsilon}(\vector{u})$ is the infinitesimal strain tensor, characterized by the symmetric part of the gradient of $\vector{u}$, i.e.,
\[
\tensor{\varepsilon}(\vector{u}) = \frac{\nabla \vector{u} + (\nabla \vector{u})^t}{2}, 
\]
and $\tensor{C}$ is the fourth-order, symmetric, bounded, and uniformly positive definite elasticity tensor, which describes the physical properties of the elastic body.
In the equilibrium equation \eqref{eq:model_equilibrium}, the divergence operator acts row-wise over $\tensor{\sigma}$, and $\vector{f}: \Omega \rightarrow \mathbb{R}^\nd$ is a given square-integrable function describing a distributed \textcolor{ok}{volume} load.
In the Dirichlet and Neumann boundary conditions, equations \eqref{eq:model_dirichlet} and \eqref{eq:model_neumann} respectively, $\vector{u}_D$ represents a prescribed displacement, $\vector{t}_N$ is a given traction field, and $\vector{n}$  denotes the outer unit normal vector.

For isotropic materials, the action of the elasticity tensor is fully characterized through the Lamé coefficients $\lambda \geq 0$ and $\mu > 0$ by the expression
\begin{equation} \label{eq:tensor_isotropic}
    \tensor{C} \tensor{\tau} = 2 \mu \tensor{\tau} +  \lambda \tr (\tensor{\tau}) \tensor{I}, \quad \forall \tensor{\tau} \in \mathbb{S},
\end{equation}
where $\tr(\cdot)$ denotes the trace of a matrix, and $\tensor{I}$ is the identity tensor.
Additionally, the Lamé coefficients can be written in terms of Young's modulus $E > 0$ and the Poisson ratio $0 \leq \nu < 0.5$ according to
\begin{equation} \label{eq:lame}
    \mu = \frac{E}{2(1+\nu)} \quad \textrm{ and } \quad \lambda = \frac{E \nu}{(1+\nu)(1-2\nu)}.
\end{equation}
We say that a material is nearly incompressible or near the incompressibility limit when $\nu$ tends to $0.5$.
For a fixed $E$, this results in the elasticity tensor $\tensor{C}$ becoming unbounded since $\lambda$ goes to infinity \citep{ainsworth2022unlocking}.

Nearly incompressible problems pose additional difficulties in the numerical approximation of \eqref{eq:model}.
For instance, it is well known that the standard low-order Galerkin approximations deteriorate for $\nu$ close to $0.5$ in a phenomenon called volumetric locking \citep{ainsworth2022unlocking, babuvska1992locking}.
In addition, the native approximations for the stress obtained by the standard Galerkin method have poor balance properties, are not $H(\div)$-conform, and tend to be more severely affected by the locking effect \cite{szabo1989stress, ainsworth2022unlocking, taraschi2024global}.
This is a major drawback since accurate approximations for the stress tensor are mandatory in many applications.

The development of Finite Element methods to accurately approximate both the displacement and stress fields on nearly incompressible problems has long been a goal of numerical analysts. 
The most natural way to achieve these requirements is using mixed displacement-stress formulations \citep{arnold1990mixed}, which approximate $\vector{u}$ and $\tensor{\sigma}$ simultaneously. 
In such formulations, it is expected that the stress approximation space 
\textcolor{ok}{is} $H(\div)$-conforming and symmetric. It should also satisfy an inf-sup condition to guarantee the well-posedness of the associated discrete problem.
Unfortunately, simultaneously meeting those restrictions proved to be quite challenging \cite{arnold2008finite, guzman2014symmetric, hu2016finite}.

When using polynomial shape functions, stable spaces that met both the symmetry and $H(\div)$-conformity requirements were first presented in \cite{arnold2002mixed} for triangular meshes. They were later developed for rectangular partitions in \cite{arnold2005rectangular}.
Regarding three-dimensional problems with tetrahedral discretizations, symmetric $H(\div)$-conforming spaces with polynomial shape functions were first devised in \cite{adams2005mixed, arnold2008finite}, and further developed in \cite{hu2015finite, hu2016finite}.
For stable and symmetric mixed spaces on cuboid meshes, we refer to \cite{hu2015new}.
Alternatively, one may consider composite elements, in which the stress approximation space is constructed using a subdivision of the mesh elements \citep{johnson1978some, arnold1984family, kvrivzek1982equilibrium, gopalakrishnan2025johnson} or elements augmented with rational bubble functions \citep{guzman2014symmetric}.
Another possible approach is to weaken the $H(\div)-$conformity of the stress approximation while keeping its symmetry, as done by \cite{JAY-GUZMAN-2011}.

A widely used approach involves altering the mixed formulation so that the symmetry of the stress tensor is enforced only in a weak sense \citep{fraeijs1965displacement, amara1979equilibrium, arnold1984peers}.
Those formulations introduce a third field related to the skew-symmetric part of an $H(\div)$-conforming tensor function (rotation).
Such a strategy 
\textcolor{ok}{has proven}
quite successful, as it significantly eases the construction of stable approximation spaces and has been widely explored \citep{boffi2009reduced, cockburn2010new, arnold2015mixed, lee2016towards}.
In particular, weakly symmetric mixed methods for general convex quadrilateral meshes were proposed, e.g., in \cite{arnold2015mixed,quinelato2019full,devloo2020enriched}. Such strategies are optimally convergent in the $H(\div)$ norm when adequate spaces are adopted.

Alongside mixed displacement-stress methods and their variations, other locking-free strategies for approximating the linear elasticity problem include:
Mixed displacement-pressure formulations \citep{franca1988new, bertrand2020equilibrated, bertrand2021weakly}, which add an auxiliary pressure-like variable;
Reduced integration methods \citep{malkus1978mixed, oden1982finite}, where the terms related to $\div \vector{u}$ are integrated with reduced quadrature schemes;
%
%
and Discontinuous Galerkin methods \citep{hansbo2002discontinuous, di2013locking, bramwell2012locking}.
{\color{ok}
The ideas behind the Discontinuous Galerkin methods can also be combined with hybridization techniques, as in the Hybridizable Discontinuous Galerkin (HDG) methods of \citep{cockburn2013superconvergent, cockburn2018devising, qiu2018hdg}.
Such methods start from a mixed formulation of the elasticity problem but relax the $H(\div)$-conformity requirement for the stress through the introduction of face unknowns (hybrid variables) associated with the displacement trace.
In that sense, HDG schemes are closely related to the hybrid version of the mixed stress-displacement formulation.
The main difference is that they introduce a numerical traction, computed from the element-wise stress and displacement unknowns, that acts as a stabilization term for the resulting discrete problem.
%
%
Among hybrid methods, Hybrid High-Order (HHO) methods \cite{di2015hybrid, carstensen2025locking} offer a particularly effective approach to nearly incompressible elasticity.
By introducing hybrid variables, local reconstruction operators, and stabilization terms, HHO methods achieve locking-free approximations based on a primal displacement formulation.
As in HDG methods, HHO hybrid variables approximate the displacement trace on the mesh skeleton.
This is not the only similarity between the two methods; deeper connections are discussed in \cite{cockburn2016bridging}.

In this work, we investigate an alternative hybrid approach to achieve locking-free approximations for the linear elasticity problem.
The proposed strategy employs a variant of the Primal Hybrid method \citep{raviart1977primal} to approximate the displacement field, followed by element-wise post-processing to recover the stress.
Primal Hybrid methods relax the $H^1(\Omega)$ conformity of the primal field (displacement) by introducing hybrid variables on element boundaries, which act as Lagrange multipliers.
 Unlike the HDG and HHO methods, the hybrid variables are no longer related to displacement traces, but to the normal component of the dual field (stress).
Furthermore, the well-posedness of the discrete problem is achieved by enriching the displacement approximation with bubble functions, rather than introducing stabilization terms in the variational form.
}
Some applications of such methodology can be found, for example, in \cite{taraschi2022convergence}, where the authors investigated the performance of the classical method of \cite{raviart1977primal} on general quadrilateral meshes,  in \cite{oliari2024posteriori}, where an \textit{a posteriori} error estimator is proposed, and in \cite{acharya2016primal, acharya2020primal}, where the authors present primal hybrid methods for parabolic problems. Multiscale counterparts of the primal hybrid method to Stokes and Brinkman equations can be found in \cite{araya2017multiscale, araya2025multiscale}, and some applications to non-linear problems are performed in \cite{milner1985primal,park1995primal}.

\textcolor{ok}{ 
For the linear elasticity problem, \cite{acharya2022primal} recently presented and analyzed a primal hybrid method that is proven to be locking-free for pure Dirichlet boundary conditions when the lowest-order approximation is employed on triangular meshes.
However, robustness to locking does not extend to higher-order spaces, as we will numerically demonstrate.
In this context, the Primal Hybrid formulation examined in the present work differs from that of \cite{acharya2022primal} by introducing an additional pressure-like field to more effectively address nearly incompressible materials when high-order spaces are used.
}
Although the introduction of an auxiliary pressure-like variable has long been used in the context of $H^1$-conforming \citep{taylor1968variational}, in the Primal Hybrid context, such an idea is less common, having been considered recently in \cite{pereira2017locking} to develop a multiscale method for the linear elasticity problem.
The resulting method was recently analyzed in \cite{gomes2024low}, where the authors adopt a stabilized strategy to solve the finer scales of the problem.
In the present work, we do not consider the multiscale aspects of the Primal Hybrid formulation and avoid the use of stabilization terms.
Hence, the development of a simple (yet general) framework for constructing stable approximation spaces for the displacement-multiplier-pressure formulation is one of our contributions.
We use this framework to construct compatible approximation spaces for two-dimensional problems discretized on triangular or convex quadrilateral meshes.

{\color{ok}
The connection between the hybrid variable in the Primal Hybrid method and the normal trace of the dual field can be leveraged to recover the dual field locally.
}
%
This was explored in \cite{chou2002flux} and \cite{correa2022optimal} to develop $H(\div)$-conforming approximations for the fluid velocity in the context of the Darcy problem.
The main contribution of this work is to adapt such ideas to the context of linear elasticity, where we devised a stress post-processing strategy for problems partitioned into triangular or convex quadrilateral elements.
Such a strategy is computed through element-wise local problems and uses standard tensor spaces generated by the Piola transformation.
As we shall demonstrate, the resulting approximations are globally $H(\div)$-conform, optimally convergent, and locking-free.
Moreover, despite not using point-wise symmetric approximation spaces, we show that the symmetry of the recovered stress holds in a weak (variational) sense.
{\color{ok}
As briefly discussed in Section \ref{sec:post-extensions}, an important aspect of the proposed post-processing is that it can be used in the context of HHO and HDG methods, extending its applicability within the broader landscape of locking-free hybrid methods.
}

In the {\color{ok}remainder} of this section, we introduce the notation and polynomial spaces used throughout the paper.
Section \ref{sec:ph} briefly discusses the standard Primal Hybrid method (PH) for linear elasticity
and then introduce the displacement-multiplier-pressure hybrid method (HDP). 
The analysis of the existence and uniqueness of solutions for the HDP method is developed in Section \ref{sec:hdp_analysis}, where some preliminary bounds for the approximation errors are shown.
Section \ref{sec:spaces} is dedicated to the construction of suitable approximation spaces for two-dimensional problems partitioned into triangular or quadrilateral meshes.
The stress recovery strategy for the HDP solution is developed and analyzed in Section \ref{sec:post}.
Finally, Section \ref{sec:exp} illustrates the analytical results from the previous sections with insightful numerical examples.

\subsection{Notations} \label{sec:notation}

For a bounded domain $D \subset \mathbb{R}^\nd$ with a Lipschitz continuous boundary $\partial D$ and $Y$ a finite-dimensional linear space, we denote by $L^2(D,Y)$ the space of square-integrable functions taking values in $Y$ and by $H^m(D,Y) \subset L^2(D,Y)$ the standard Sobolev spaces.
The $L^2$ norm is denoted by $\| \cdot \|_{0,D}$, while $\| \cdot \|_{m,D}$ and $\vert \cdot \vert_{m,D}$ are the standard norm and semi-norm over $H^m(D,Y)$, respectively.
We set $H^{1/2}(\partial D,Y)$ as the space made of the traces of the functions in $H^1(D,Y)$ and $H^{-1/2}(\partial D,Y)$ as its dual space.
Also, for $\Gamma \subset \partial D$ such that $\vert \Gamma\vert >0$, $H^{1}_{0,\Gamma}(D,Y)$ denotes the subset of $H^{1}(D,Y)$ of functions with vanishing trace on $\Gamma$.
In this work, $Y$ will be either $\mathbb{R}$, $\mathbb{R}^\nd$, $\mathbb{M}$, or $\mathbb{S}$ {\textcolor{ok}{(with $\mathbb{M}$ and $\mathbb{S}$ defined as in the first paragraph of the text)}}.

For $Y = \mathbb{R}^\nd$, $\mathbb{M}$, or $\mathbb{S}$, we also set $H(\div, D, Y)$ as the subspace of $L^2(D, Y)$ made of functions with square-integrable divergence (where the divergence operator is applied row-wise when $Y$ is a space of tensors), and define the norm
\[
\| \tensor{\tau} \|_{H(\div)}^2 = \| \tensor{\tau} \|^2_{0,D} + \| \div \tensor{\tau} \|^2_{0,D}. 
\]
Taking $\Gamma \subset \partial D$, we then define the subspace
\[ H_{0,\Gamma}(\div, D, \mathbb{S}) = \left\{ \tensor{\tau} \in H(\div, D, \mathbb{S}) \, : \, \int_{\partial D} (\tensor{\tau} \vector{n}) \cdot \vector{v} \des = 0, \; \forall \vector{v} \in H^1_{0,\partial D \setminus \Gamma}(D, \mathbb{R}^\nd) \right\}. 
\]
Here and thereafter, the boundary integrals are understood as the duality product between $H^{-1/2}(\partial D, \mathbb{R}^\nd)$ and $H^{1/2}(\partial D, \mathbb{R}^\nd)$.

Now let $\TauH$ be a partition of $\Omega$ into non-overlapping, open subdomains $K$, called elements.
For this work, we assume that the partitions $\TauH$, referred to as meshes, have no hanging nodes and that {\color{ok} the elements $K$ are convex polygons or polyhedra.} 
We shall make use of the broken Sobolev spaces $H^r(\TauH, Y)$, defined by
\begin{equation} \label{eq:broken_sobolev}
    H^r(\TauH, Y) = \{ v \in L^2(\Omega, Y) \, : \, v\vert_K \in H^r(K, Y), \; \forall K \in \TauH \}
\end{equation}
and equipped with the semi-norm
\begin{equation} \label{eq:broken_sobolev_semi-norm}
    \vert v \vert_{r,\TauH} = \left( \sum_{K \in \TauH} \vert v \vert_{r,K}^2 \right)^{1/2}.
\end{equation}
Furthermore, for given functions $\vector{v}, \vector{w} \in \ds \prod_{K \in \TauH} L^2(K, \mathbb{R}^\nd)$, we define the product
\[
(\vector{v}, \vector{w})_{\TauH} = \sum_{K \in \TauH} \int_{K} \vector{v} \cdot \vector{w} \dex. 
\]
For tensor functions in $\ds \prod_{K \in \TauH} L^2(K, \mathbb{M})$ or $\ds \prod_{K \in \TauH} L^2(K, \mathbb{S})$, the same notation is used to define the product
\[ (\vector{v}, \vector{w})_{\TauH} = \sum_{K \in \TauH} \int_{K} \vector{v} : \vector{w} \dex, 
\]
where `$:$' denotes the Frobenius inner product between two tensor fields.
\sloppy Analogously, for $\vector{l} \in \prod_{K \in \TauH} H^{-1/2}(\partial K, \mathbb{R}^\nd)$ and $\vector{v} \in \prod_{K \in \TauH} H^{1/2}(\partial K, \mathbb{R}^\nd)$, we set
\[
\langle \vector{l}, \vector{v} \rangle_{\partial \TauH} = \sum_{K \in \TauH} \int_{\partial K} \vector{l} \cdot \vector{v} \des. 
\]

We end this section by defining some two-dimensional polynomial spaces that will be used in Section \ref{sec:spaces}.
For an open polygon $D \subset \mathbb{R}^2$, we denote by $P_r(D,\mathbb{R})$ the space of polynomials on $D$ with total degree equal to or less than $r$, by $\tilde{P}_r(D, \mathbb{R})$ the set of homogeneous polynomials of exact degree $r$, and by $P_{r,s}(D, \mathbb{R})$ the space of polynomials with degree at most $r$ in the first coordinate and at most $s$ in the second one.
We also set $Q_r(D, \mathbb{R}) = P_{r,r}(D, \mathbb{R})$ and define $R_r(D, \mathbb{R})$ as the space spanned by the monomials of $Q_{r+1}(D, \mathbb{R})$ except for $x_1^{r+1} x_2^{r+1}$, where $\vector{x} = (x_1,x_2)$ 
is the vector of spatial coordinates.
Over the boundary $\partial D$, we consider the space $E_r(\partial D, \mathbb{R})$, defined as the subspace of $L^2(\partial D, \mathbb{R})$ such that its elements are polynomials of degree equal to or less than $r$ when restricted to each face (edge for two-dimensional domains) of $D$.
We set $P_r(D, \mathbb{R}^2) = P_r(D, \mathbb{R}) \times P_r(D, \mathbb{R})$ and denote by $P_r(D, \mathbb{M})$ the space of tensor-valued functions with entries in $P_r(D, \mathbb{R})$.
The same convection is used to extend the spaces $Q_r(D, \mathbb{R})$ and $E_r(\partial D, \mathbb{R})$ to vector or tensor-valued functions.
{\color{ok}
It is also useful to extend some of these spaces to the case where $D \subset \mathbb{R}^3$ is a tetrahedron, as commented in Remark \ref{rk:hdp3d}. In this setting, $P_r(D, \mathbb{R})$ 
is naturally extended from its two-dimensional counterpart, and $E_r(D, \mathbb{R})$ is the subspace of $L^2(\partial D, \mathbb{R})$ such that its elements are polynomials of total degree equal to or less than $r$ when restricted to each triangular face of $D$.
%
%
}

\section{Primal Hybrid methods for the linear elasticity} \label{sec:ph}

Before discussing the displacement-multiplier-pressure hybrid method, we briefly introduce the standard Primal Hybrid method for linear elasticity, written only in terms of the displacement and multiplier fields.

\subsection{The standard Primal Hybrid method}

Assuming that $\TauH$ is a partition of $\Omega$ as described in Section \ref{sec:notation}, we begin by introducing the {\color{ok}mesh-dependent} spaces
\begin{equation} \label{eq:def_X}
    \mathcal{X}(\TauH) = \{ \vector{v} \in L^2(\Omega, \mathbb{R}^\nd) \, : \, \vector{v}\vert_K \in H^1(K, \mathbb{R}^\nd), \; \forall K \in \TauH \}
\end{equation}
and
\begin{multline} \label{eq:def_M}
    \mathcal{M}(\TauH) = \left\{ \vector{l} \in \prod_{K \in \TauH} H^{-1/2}(\partial K, \mathbb{R}^\nd) \, : \, \exists \, \tensor{\tau} \in H_{0, \Gamma_N}(\div,\Omega, \mathbb{S}) \textrm{ such that } \right. \\
    \left. \tensor{\tau} \vector{n}^{\partial K} = \vector{l} \textrm{ over } \partial K, \; \forall K \in \TauH \right\}.
\end{multline}
{\color{ok}Notice that space $\mathcal{M}$ also depends on the Neumann boundary $\Gamma_N$ (see eq. \eqref{eq:model_neumann}).}
%
%
%
Denoting by $\diamOm$ the diameter of $\Omega$, and considering the following auxiliary set, for a given $\vector{l}\in \mathcal{M}$,
\begin{equation} \label{eq:Tmu}
  T(\vector{l})= \{ \tensor{\tau} \in H_{0, \Gamma_N}(\div, \Omega, \mathbb{S}) \, : \, \tensor{\tau} \vector{n}^{\partial K} = \vector{l}\vert_{\partial K} \textrm{ over } \partial K, \; \forall K \in \TauH \},
\end{equation}
we define the norms
\begin{equation} \label{eq:norm_X}
    \| \vector{v} \|_{\mathcal{X}} = \left[ \sum_{K \in \TauH}\left( \frac{1}{\diamOm^2} \| \vector{v} \|^2_{0,K} + \| \tensor{\varepsilon}(\vector{v}) \|^2_{0,K} \right)\right]^{1/2}
\end{equation}
and
\begin{equation} \label{eq:norm_M}
    \| \vector{l} \|_{\mathcal{M}} = \inf_{\tensor{\tau} \in T(\vector{l})} \left[ \sum_{K \in \TauH} \left(\| \tensor{\tau} \|_{0,K}^2 + \diamOm^2 \| \div \tensor{\tau} \|_{0,K}^2\right) \right]^{1/2}
\end{equation}
over $\mathcal{X}$ and $\mathcal{M}$, respectively.
Notice that $\mathcal{X}$ is the broken Sobolev space $H^1(\TauH, \mathbb{R}^\nd)$, but equipped with a special norm $\| \cdot \|_{\mathcal{X}}$.

The traditional Primal Hybrid formulation for the linear elasticity problem \eqref{eq:model} reads as: 
{\it find the  displacement field $\vector{u} \in \mathcal{X}$ and the Lagrange multiplier $\vector{m} \in \mathcal{M}$ such that
\begin{subequations} \label{eq:formulation_PH}
    \begin{align}
        a_{\ph}(\vector{v}, \vector{u}) + b(\vector{v}, \vector{m}) & = (\vector{f}, \vector{v})_{\Omega} + \langle \vector{t}_N, \vector{v} \rangle_{\Gamma_N} & \forall \vector{v} \in \mathcal{X}, \label{eq:formulation_PH_a} \\
        b(\vector{u}, \vector{l}) & =-  \langle \vector{l} ,\vector{u}_D  \rangle_{\Gamma_D} & \forall \vector{l} \in \mathcal{M}, \label{eq:formulation_PH_b}
    \end{align}
\end{subequations}
}
where the bilinear forms $a_{\ph}: \mathcal{X} \times \mathcal{X} \rightarrow \mathbb{R}$ and $b: \mathcal{X} \times \mathcal{M} \rightarrow \mathbb{R}$ are given by
\begin{equation} \label{eq:a_PH}
    a_{\ph}(\vector{u}, \vector{v}) = (\tensor{C} \tensor{\varepsilon}(\vector{u}), \tensor{\varepsilon}(\vector{v}))_{\TauH},
\end{equation}
\begin{equation} \label{eq:b}
    b(\vector{v}, \vector{l}) = - \langle \vector{l}, \vector{v} \rangle_{\partial \TauH}.
\end{equation}
{\color{ok}
Notice that the Dirichlet condition over $\Gamma_D$ is imposed weakly on the right-hand side of \eqref{eq:formulation_PH_b}.
}

Following the ideas of \cite{raviart1977primal}, the existence and uniqueness of the solution for \eqref{eq:formulation_PH} is discussed in \cite{harder2016hybrid}, where it is shown the following characterization for the Lagrange multiplier
\begin{equation} \label{eq:charac_m_PH}
    \vector{m}\vert_{\partial K} = \tensor{C} \tensor{\varepsilon}(\vector{u}) \vector{n}^{\partial K} = \tensor{\sigma} \vector{n}^{\partial K} \quad \textrm{ over } \partial K \setminus \Gamma_N, \; \forall K \in \TauH,
\end{equation}
connecting $\vector{m}$ with the traction field.
Also, equation \eqref{eq:formulation_PH_b} guarantees that the solution $\vector{u}$ belongs to $H^1(\Omega, \mathbb{R}^2)$ (see Lemma 1 of \cite{raviart1977primal}).

The Primal Hybrid method (PH) is obtained by replacing the spaces $\mathcal{X}$ and $\mathcal{M}$ with finite-dimensional subspaces $\mathcal{X}_h \subset \mathcal{X}$ and $\mathcal{M}_h \subset \mathcal{M}$ in \eqref{eq:formulation_PH}, resulting in a saddle-point type discrete problem.
Therefore, to guarantee the existence and uniqueness of the solution, the approximation spaces $\mathcal{X}_h \times \mathcal{M}_h$ need to be inf-sup stable, and the bilinear form $a_{\ph}(\cdot, \cdot)$ must be coercive on a discrete kernel space \citep{boffi2013mixed}.
Those two conditions are summarized below
\begin{enumerate}[label=\textit{(A\arabic*)}]
    \item Defining the kernel space associated with $b(\cdot, \cdot)$ as \label{it:PH_coercivity}
    \begin{equation} \label{eq:kernel_ph}
        \mathcal{N}_h^{b}
        = \{ \vector{v} \in \mathcal{X}_h \, : \, {\color{ok}b(\vector{v}, \vector{l})} = 0, \; \forall \vector{l} \in \mathcal{M}_h \},
    \end{equation}
    the bilinear form $a_{\ph}$ is said to be coercive on
    $\mathcal{N}_h^{b}$
    if there is a constant $\alpha_{\ph} > 0$ such that
    \[
    a_{\ph}(\vector{v}, \vector{v}) \geq \alpha_{\ph} \| \vector{v} \|_{\mathcal{X}}^2 \quad \forall \vector{v} \in \mathcal{N}_h^{b}. 
    \]
    \item A pair of spaces $\mathcal{X}_h \times \mathcal{M}_h$ is inf-sup stable for formulation \eqref{eq:formulation_PH} if there is a constant $\beta_{\ph} > 0$ such that \label{it:PH_inf-sup}
    \[
    \sup_{\vector{v} \in \mathcal{X}_h} \frac{b(\vector{v}, \vector{l})}{\| \vector{v} \|_{\mathcal{X}}} \geq \beta_{\ph} \| \vector{l} \|_{\mathcal{M}} \quad \forall \vector{l} \in \mathcal{M}_h.
    \]
\end{enumerate}
In the condition \ref{it:PH_inf-sup}, and in the remainder of the text, we assume that the supremum is always taken excluding the zero function.
The construction of compatible approximation spaces satisfying both \ref{it:PH_coercivity} and \ref{it:PH_inf-sup} is postponed to Section \ref{sec:inf-sup_PH}.

\

\paragraph{\bf The AP method}
For pure Dirichlet problems on isotropic materials, the authors of \cite{acharya2022primal}
replace the bilinear form $a_{\ph}(\cdot, \cdot)$ in \eqref{eq:formulation_PH_a} by
\begin{equation} \label{eq:a_aph}
    a_{\ap}(\vector{u}, \vector{v}) = (2 \mu \nabla \vector{u} + (\lambda + \mu)\left(\div \vector{u}\right)\bm{I}, \nabla \vector{v})_{\TauH},
\end{equation}
introducing a primal hybrid method in which the Lagrange multiplier $\vector{m}$ is given by
\begin{equation} \label{eq:charac_m_AP}
    \vector{m}\vert_{\partial K} = 2 \mu (\nabla \vector{u}) \vector{n}^{\partial K} + (\mu + \lambda) (\div \vector{u}) \, \vector{n}^{\partial K} \quad \textrm{ over } \partial K, \; \forall K \in \TauH.
\end{equation}
In the present work, we refer to this hybrid method as AP (Acharya-Porwal).

One advantage of using the AP method is that it allows for using lower-order approximation spaces.
This choice is not available for the PH method, as discussed in Remark \ref{rk:coercivity}.
On the other hand, the Lagrange multiplier $\vector{m}$ can no longer be directly associated with the traction field.
{\color{ok}
It is also important to notice that the proof of convergence and locking robustness presented in \cite{acharya2022primal} is strongly tailored to
the piecewise-constant nature of the Lagrange multipliers (see, for instance, Lemma 4.2 of \cite{acharya2022primal}). Therefore, their results cannot be directly extended to higher-order spaces, even though the formulation allows their use.
}

\subsection{The hybrid displacement-multiplier-pressure  method (HDP)} \label{sec:hdp}

 The 
{\color{ok} PH method provides accurate approximations of both $\vector{u}$ and $\vector{m}$ for compressible problems.
}
On nearly incompressible problems, however, approximation errors may increase as $\nu$ tends to $0.5$, especially for the Lagrange multiplier.
The AP method can accurately approximate nearly incompressible problems if the lowest-order possible approximation space is employed on triangular partitions \citep{acharya2022primal}.
Unfortunately, this robustness does not seem to hold for higher-order spaces, as indicated by the numerical experiments performed in Section \ref{sec:exp_locking}.

To better deal with the incompressibility limit, we adopt a hybrid method for the problem \eqref{eq:model} involving three fields. 
The first two fields are the displacement vector $\vector{u}$ and the Lagrange multiplier $\vector{m}$, as in \eqref{eq:formulation_PH}.
The third field is an auxiliary pressure-like variable denoted by $p$. This variable belongs to $\mathcal{P} = L^2(\Omega, \mathbb{R}) \backslash \mathbb{R}$, when $\Gamma_D = \partial \Omega$, and $\mathcal{P} = L^2(\Omega, \mathbb{R})$, otherwise.
With the introduction of $p$, the hybrid displacement-multiplier-pressure formulation (HDP) for isotropic elasticity problems reads: {\it  find the {\color{ok}triple} $(\vector{u}, \vector{m}, p) \in \mathcal{X} \times \mathcal{M} \times \mathcal{P}$ such that
\begin{subequations} \label{eq:formulation_HDP}
    \begin{align}
        a(\vector{v}, \vector{u}) + b(\vector{v}, \vector{m}) + c(\vector{v}, p) & = (\vector{f}, \vector{v})_{\Omega} + \langle \vector{t}_N, \vector{v} \rangle_{\Gamma_N} & \forall \vector{v} \in \mathcal{X}, \label{eq:formulation_HDP_a} \\
        b(\vector{u}, \vector{l}) & = - \langle \vector{l} ,\vector{u}_D  \rangle_{\Gamma_D}& \forall \vector{l} \in \mathcal{M}, \label{eq:formulation_HDP_b}\\
        c(\vector{u}, q) - \frac{1}{\lambda}(p,q)_{\Omega} & = 0 & \forall q \in \mathcal{P}, \label{eq:formulation_HDP_c}
    \end{align}
\end{subequations}
}
where $b(\cdot, \cdot)$ is given as in \eqref{eq:b}, and the bilinear forms $a: \mathcal{X} \times \mathcal{X} \rightarrow \mathbb{R}$ and $c:\mathcal{X} \times \mathcal{P} \rightarrow \mathbb{R}$ are defined according to
\begin{equation} \label{eq:a_and_c}
    a(\vector{v}, \vector{u}) = (2 \mu \, \tensor{\varepsilon}(\vector{u}), \tensor{\varepsilon}(\vector{v}))_{\TauH} \quad \textrm{ and } \quad c(\vector{v}, q) = (\div \vector{v}, q)_{\TauH}.
\end{equation}

Formulation \eqref{eq:formulation_HDP} is employed in \cite{pereira2017locking} and can be seen as a penalized version of the Primal Hybrid formulation for the Stokes problem \citep{araya2017multiscale, araya2025multiscale}.
{\color{ok}
In fact, the similarity with Stokes formulations motivates referring to $p$ as ‘pressure’, even though this is not its physical meaning in the elasticity equations.
}
Furthermore, following closely the proof of Theorem 6 from \cite{araya2017multiscale}, it is possible to show that the variational problem \eqref{eq:formulation_HDP} admits a unique solution, and the following characterizations hold
\begin{equation} \label{eq:charac_p}
    p = \lambda \div \vector{u} \quad \textrm{ in } \Omega,
\end{equation}
\begin{equation} \label{eq:charac_mult_HDP}
    \vector{m} = (2 \mu \tensor{\varepsilon}(\vector{u}) + p \tensor{I}) \vector{n}^{\partial K} = \tensor{\sigma} \vector{n}^{\partial K} \quad \textrm{ over } \partial K \setminus \Gamma_N, \; \forall K \in \TauH.
\end{equation}

The hybrid displacement-multiplier-pressure method (HDP) is obtained by replacing the spaces $\mathcal{X} \times \mathcal{M} \times \mathcal{P}$ with finite-dimensional subspaces $\mathcal{X}_h \subset \mathcal{X}$, $\mathcal{M}_h \subset \mathcal{M}$, and $\mathcal{P}_h \subset \mathcal{P}$.
The associated discrete problem then reads: {\it find $(\vector{u}_h, \vector{m}_h, p_h) \in \mathcal{X}_h \times \mathcal{M}_h \times \mathcal{P}_h$ satisfying
\begin{subequations} \label{eq:method_HDP}
    \begin{align}
        a(\vector{u}_h, \vector{v}) + b(\vector{v}, \vector{m}_h) + c(\vector{v}, p_h) & = (\vector{f}, \vector{v})_{\Omega} + \langle \vector{t}_N, \vector{v} \rangle_{\Gamma_N} & \forall \vector{v} \in \mathcal{X}_h, \label{eq:method_HDP_a} \\
        b(\vector{u}_h, \vector{l}) & = - \langle \vector{l}, \vector{u}_D \rangle_{\Gamma_D} & \forall \vector{l} \in \mathcal{M}_h, \label{eq:method_HDP_b}\\
        c(\vector{u}_h, q) - \frac{1}{\lambda}(p_h,q)_{\Omega} & = 0 & \forall q \in \mathcal{P}_h. \label{eq:method_HDP_c}
    \end{align}
\end{subequations}
}
Similarly to the standard PH method, the spaces $\mathcal{X}_h$, $\mathcal{M}_h$, and $\mathcal{P}_h$ need to be carefully chosen to guarantee the existence and uniqueness of solution for \eqref{eq:method_HDP}.

{\color{ok}
\begin{remark}[Comments on nonlinear problems]
 It is worth mentioning that the proposed primal hybrid displacement–pressure 
 framework may offer a basis for solving nonlinear problems such as hyperelasticity, which
is (almost) incompressible by nature. In general, such problems admit a volumetric-isochoric decomposition of the strain-energy density, which leads to a mixed displacement–pressure formulation \cite{kadapa2022linearized}.
As an alternative to standard mixed formulations, hybridization principles can be employed to weakly enforce the $H^1$-conformity of the displacement field through traction hybrid variables. 
Typically, the resulting nonlinear system is solved by using an iterative scheme.
Each iteration relies on solving a problem similar to the linear elasticity system, meaning the HDP method can be applied to determine the increments, assuring robustness in the incompressibility limit.
\end{remark}
}

\section{Analysis of the HDP method} \label{sec:hdp_analysis}

In this section, we analyze the well-posedness of the HDP method, showing that stable approximations can be obtained 
by combining stable spaces for the PH method with compatible spaces for the standard Galerkin method for Stokes.
A possible approach for analyzing the HDP method is to restrict the multiscale-based analysis of \cite{araya2025multiscale} to the case where the sub-meshes are composed of only one element.
Here, we adopt a different strategy, however, showing that when no multiscale aspects are considered, stability can be proven in a simpler way.

To analyze the well-posedness of the HDP method, we rephrase formulation \eqref{eq:method_HDP} as: {\it find $((\vector{u}_h, p_h), \vector{m}_h) \in \mathcal{H}_h \times \mathcal{M}_h$ such that
\begin{subequations} \label{eq:method_HDP_aux}
    \begin{align}
    \ahdp((\vector{u}_h,p_h), (\vector{v},q)) + \mathfrak{b}((\vector{v},q), \vector{m}_h) & = (\vector{f}, \vector{v})_{\Omega} + \langle \vector{t}_N, \vector{v} \rangle_{\Gamma_N} & \forall (\vector{v}, q) \in \mathcal{H}_h, 
    \label{eq:method_HDP_aux_a} \\
    \mathfrak{b}((\vector{u}_h, p_h), \vector{l} ) & = - \langle \vector{l}, \vector{u}_D \rangle_{\Gamma_D}, & \forall \vector{l} \in\mathcal{M}_h, \label{eq:method_HDP_aux_b}
    \end{align}
\end{subequations}
}
where $\mathcal{H}_h = \mathcal{X}_h \times \mathcal{P}_h$ is a finite dimensional subspace of $\mathcal{H} = \mathcal{X} \times \mathcal{P}$, 
equipped with the norm
\begin{equation} \label{eq:norm_H}
	\| (\vector{v}, q) \|_{\mathcal{H}} = \| \vector{v} \|_{\mathcal{X}} + \| q \|_{0,\Omega}.
\end{equation}
The bilinear forms $\ahdp(\cdot, \cdot): \mathcal{H} \times \mathcal{H} \rightarrow \mathbb{R}$ and $\mathfrak{b}(\cdot, \cdot): \mathcal{H} \times \mathcal{M} \rightarrow \mathbb{R}$ are given by
\begin{equation} \label{eq:form_a_hdp_aux}
	\ahdp((\vector{u},p), (\vector{v},q)) = a(\vector{u}, \vector{v}) + c(\vector{v},p) - c(\vector{u}, q) + \frac{1}{\lambda}(p, q)_{\Omega}
\end{equation}
and
\begin{equation} \label{eq:form_b_hdp_aux}
	\mathfrak{b}((\vector{v},q), \vector{l}) = b(\vector{v}, \vector{l}).
\end{equation}

\begin{theorem}[Well-posedness of \eqref{eq:method_HDP_aux}] \label{theo:wp_discrete}
The discrete problem \eqref{eq:method_HDP_aux} (and consequently its equivalent form \eqref{eq:method_HDP}) has a unique solution $((\vector{u}_h, p_h), \vector{m}_h) \in \mathcal{H}_h \times \mathcal{M}_h$ if the following conditions hold:
\begin{enumerate}
[label=\textit{(B\arabic*)}]
    %
    \item There is a constant $\alpha_{\hdp} > 0$ such that the inequality \label{it:hdp_coercivity}
    \begin{equation}
    	\sup_{(\vector{v}, q) \in \mathcal{N}_h^{\hdp}} \frac{\ahdp((\vector{w},z) , (\vector{v},q))}{\| (\vector{v}, q) \|_{\mathcal{H}}} \geq \alpha_{\hdp} \| (\vector{w},z) \|_{\mathcal{H}} 
    \end{equation}
	is satisfied for every $(\vector{w}, z)$ in the kernel of $\mathfrak{b}(\cdot,\cdot)$
	\begin{equation} \label{eq:kernel_hdp}
		\mathcal{N}^{\hdp}_h = \{ (\vector{w}, z) \in \mathcal{H}_h \, : \, \mathfrak{b}((\vector{w}, z), \vector{l}) = 0, \; \forall \vector{l} \in \mathcal{M}_h \}= \mathcal{N}^{b}_h\times \mathcal{P}_h.
	\end{equation}
    \item There is a constant $\beta_{\hdp} > 0$ such that the inequality
    \begin{equation}
		\sup_{(\vector{v}, q) \in \mathcal{H}_h} \frac{\mathfrak{b}((\vector{v}, q), \vector{l})}{\| (\vector{v},q) \|_{\mathcal{H}}} \geq \beta_{\hdp} \| \vector{l} \|_{\mathcal{M}}
	\end{equation}
	holds for every $\vector{l} \in \mathcal{M}_h$. \label{it:hdp_inf-sup} 
\end{enumerate}
In this case, the following error bounds can be shown
\begin{equation}
\label{eq:firstbound_u_w_hdp}
\| (\vector{u},p) - (\vector{u}_h, p_h) \|_{\mathcal{H}} \leq C_1 \inf_{(\vector{v}, q) \in \mathcal{H}_h} 
\| (\vector{u},p) - (\vector{v},q) \|_{\mathcal{H}} + C_2 \inf_{\vector{l} \in \mathcal{M}_h} 
\| \vector{m} - \vector{l} \|_{\mathcal{M}},
\end{equation}
\begin{equation} \label{eq:firstbound_m_hdp}
\| \vector{m} - \vector{m}_h \|_{\mathcal{M}} \leq C_3 \inf_{(\vector{v},q) \in \mathcal{H}_h} 
\| (\vector{u},p) - (\vector{v},q) \|_{\mathcal{H}} + C_4 \inf_{\vector{l} \in \mathcal{M}_h} 
\| \vector{m} - \vector{l} \|_{\mathcal{M}},
\end{equation}
where  $(\vector{u}, \vector{m}, p) \in \mathcal{X} \times \mathcal{M} \times \mathcal{P}$ are the exact solutions of \eqref{eq:formulation_HDP}, and the constants $C_1$ to $C_4$ depend only on $\alpha_{\hdp}$, $\beta_{\hdp}$ and on the continuity constants $C_{\ahdp}$ and $C_{\mathfrak{b}}$ of the bilinear forms $\ahdp(\cdot, \cdot)$ and $\mathfrak{b}(\cdot, \cdot)$, respectively.
{\color{ok} In particular, $C_1-C_4$ remain bounded when $\lambda \rightarrow \infty$ provided that $\alpha_{\hdp}$ and $\beta_{\hdp}$ are bounded away from zero in that limit.}
\end{theorem}
\begin{proof}
    Such results are obtained by applying the standard theory of saddle-point problems approximation. See for instance Theorem 2.34 from \cite{ern2004theory}.
\end{proof}

In the {\color{ok}remainder} of this section, we explore sufficient conditions over $\mathcal{X}_h \times \mathcal{M}_h \times \mathcal{P}_h$ to satisfy \ref{it:hdp_coercivity} and \ref{it:hdp_inf-sup}.
As a result, we obtain a straightforward procedure to construct stable approximation spaces, which is latter applied for two-dimensional problems in Section \ref{sec:spaces}.

\subsection{Continuity and Inf-Sup} \label{sec:hdp_analysis_1}

Before addressing the conditions of Theorem \ref{theo:wp_discrete}, we first establish some bounds for the continuity constants.

\begin{lemma}[Continuity constants]
The bilinear forms $\ahdp(\cdot, \cdot)$ and $\mathfrak{b}(\cdot, \cdot)$, defined in \eqref{eq:form_a_hdp_aux} and \eqref{eq:form_b_hdp_aux}, are continuous with continuity constants given by
\[
C_{\ahdp} = \max \left\{2 \mu, \sqrt{\nd}, \frac{1}{\lambda} \right\} \quad \textrm{ and } \quad C_{\mathfrak{b}} = 1,
\]
respectively.
\end{lemma}
\begin{proof}
    We first show the continuity for the base bilinear forms $a(\cdot, \cdot)$, $b(\cdot, \cdot)$, and $c(\cdot, \cdot)$.
    Considering $\vector{v}, \vector{w} \in \mathcal{X}$ and $q \in \mathcal{P}$ arbitrary functions, it follows from the Cauchy-Schwarz inequality that
    \begin{equation*}
	   a(\vector{v}, \vector{w}) \leq 2 \mu \left(\sum_{K \in \TauH} \| \tensor{\varepsilon}(\vector{v}) \|^2_{0,K}\right)^{1/2} \left(\sum_{K \in \TauH} \| \tensor{\varepsilon}(\vector{w}) \|^2_{0,K} \right)^{1/2} \leq 2 \mu \|\vector{v}\|_{\mathcal{X}} \| \vector{w} \|_{\mathcal{X}}
	\end{equation*}
    and
	\begin{equation*}
		c(\vector{v}, q) \leq \left( \sum_{K \in \TauH} \| \div \vector{v} \|^2_{0,K} \right)^{1/2} \| q \|_{0,\Omega} \leq \sqrt{\nd} \| \vector{v} \|_{\mathcal{X}} \| q \|_{0,\Omega},		
	\end{equation*}
 where we used the relation
$\|\div \vector{v}\|_{0,K}\leq \sqrt{\nd}\|\tensor{\varepsilon}(\vector{v})\|_{0,K},\,\forall \vector{v}\in H^1(K,\mathbb{R}^{\nd}).
$
Next, for an arbitrary function $\vector{l} \in \mathcal{M}$, let $\tensor{\tau} \in T(\vector{l})$, as defined in \eqref{eq:Tmu}.
Applying the Green's formula element-wise, we get that
\begin{equation*}
b(\vector{v}, \vector{l}) = \langle \tensor{\tau}\vector{n}^{\partial K},  \vector{v}\rangle_{\partial \TauH}
= (\div \tensor{\tau} ,\vector{v})_{\TauH} + ( \tensor{\tau} , \tensor{\varepsilon}(\vector{v}))_{\TauH}.
\end{equation*}
Then, it follows from the Cauchy-Schwarz and the triangular inequalities that
\begin{equation*}
\begin{split}
b(\vector{v}, \vector{l}) & \leq \sum_{K \in \TauH} \left( \| \div \tensor{\tau} \|_{0,K} \| \vector{v} \|_{0,K} + \|\tensor{\tau} \|_{0,K} \| \tensor{\varepsilon}(\vector{v}) \|_{0,K}  \right) \\
& \leq \left( \sum_{K \in \TauH} \| \tensor{\tau} \|^2_{0, K} + \diamOm^2 \| \div \tensor{\tau} \|^2_{0,K} \right)^{1/2} \left( \sum_{K \in \TauH} \frac{1}{\diamOm^2} \| \vector{v} \|^2_{0, K} + \| \tensor{\varepsilon}(\vector{v}) \|^2_{0,K} \right)^{1/2} \\
& = \left( \sum_{K \in \TauH} \| \tensor{\tau} \|^2_{0, K} + \diamOm^2 \| \div \tensor{\tau} \|^2_{0,K} \right)^{1/2}\| \vector{v} \|_{\mathcal{X}},
\end{split}
\end{equation*}
and the continuity of $b(\cdot, \cdot)$ is obtained taking the infimum over $T(\vector{l})$.
Using the definitions \eqref{eq:norm_H}--\eqref{eq:form_b_hdp_aux}, the continuity for $\ahdp(\cdot, \cdot)$ and $\mathfrak{b}(\cdot, \cdot)$ follow directly from the continuity of the base bilinear forms.
\end{proof}

Now, we show that Condition \ref{it:hdp_inf-sup} (inf-sup) is satisfied under a suitable choice of the discrete spaces  $\mathcal{X}_h$ and $\mathcal{M}_h$.
\begin{lemma}[Condition \ref{it:hdp_inf-sup}]
If $\mathcal{X}_h \times \mathcal{M}_h$ is inf-sup stable for the standard PH method (condition \ref{it:PH_inf-sup}), then \ref{it:hdp_inf-sup} holds for any choice of $\mathcal{P}_h$ with
\[
\beta_{\hdp} = \beta_{\ph}. 
\]
{\color{ok}Moreover, $\beta_{\hdp}$ does not depend on $\lambda$.}
\end{lemma}
\begin{proof}
In fact, if $\mathcal{X}_h \times \mathcal{M}_h$ satisfies \ref{it:PH_inf-sup}, then it is direct to verify that
\[
\beta_{\ph} \| \vector{l} \|_{\mathcal{M}} \leq \sup_{\vector{v} \in \mathcal{X}_h} \frac{b(\vector{v}, \vector{l})}{\| \vector{v} \|_{\mathcal{X}}} = \sup_{(\vector{v},0) \in \mathcal{X}_h\times\mathcal{P}_h} \frac{b(\vector{v}, \vector{l})}{\| \vector{v} \|_{\mathcal{X}} + {\| 0 \|_{0, \Omega}}} 
=
\sup_{(\vector{v}, q) \in \mathcal{H}_h} \frac{\mathfrak{b}((\vector{v}, q), \vector{l})}{\| (\vector{v}, q) \|_{\mathcal{H}}}
\]
for every $\vector{l} \in \mathcal{M}_h$, from which we conclude \ref{it:hdp_inf-sup}. 
\end{proof}

\subsection{Coercivity on the kernel} \label{sec:hdp_analysis_2}

To address condition \ref{it:hdp_coercivity}, we first establish a coercivity result for the bilinear form $a(\cdot, \cdot)$.
For that, consider the local rigid-body motion spaces
\[
\mathcal{X}_{\mathrm{rm}}^{K} = \{ \vector{v} \in H^1(K, \mathbb{R}^\nd) \, : \, \tensor{\varepsilon}(\vector{v}) = \vector{0} \textrm{ in } K \}
\]
and 
\[
\mathcal{X}_{\mathrm{rm}}^{\partial K} = \{ \vector{\mu} \in H^{1/2}(\partial K, \mathbb{R}^\nd) \, : \, \vector{\mu}\vert_e = \vector{v}_\mathrm{rm}\vert_e \textrm{ for every  edge/face } e \textrm{ of } \partial K, \, \textrm{ with } \vector{v}_\mathrm{rm} \in \mathcal{X}_{\mathrm{rm}}^{K} \}, 
\]
from which we construct the global spaces
\begin{equation} \label{eq:def_Xrm}
    \mathcal{X}_\mathrm{rm}^{\TauH} = \{ \vector{v} \in \mathcal{X} \, : \, \vector{v}\vert_K \in \mathcal{X}_{\mathrm{rm}}^{K} , \; \forall K \in \TauH \}
\end{equation}
and
\begin{equation} \label{eq:def_Mrm}
    \mathcal{X}_{\mathrm{rm}}^{\partial \TauH} = \{ \vector{\mu} \in \mathcal{M} \, : \, \vector{\mu}\vert_{\partial K} \in \mathcal{X}_{\mathrm{rm}}^{\partial K} , \; \forall K \in \TauH \}.
\end{equation}
Setting $\tilde{\mathcal{X}}=\left(\mathcal{X}_\mathrm{rm}^{\TauH}\right)^\perp$, the $L^2$-orthogonal complement of $\mathcal{X}_\mathrm{rm}^{\TauH}$ in $\mathcal{X}$, the following broken 
{\color{ok} Poincaré and Korn} inequalities hold \citep{gomes2022mhm,brenner2004korn}
\begin{equation} \label{eq:broken_korn_local}
    \| \nabla \tilde{\vector{v}} \|_{0,K} \geq \frac{1}{C_{\mathrm{Pc}} h_K} \| \tilde{\vector{v}} \|_{0,K} \quad \textrm{ and } \quad \| \tensor{\varepsilon}( \tilde{\vector{v}}) \|_{0,K} \geq \frac{1}{C_{\mathrm{Kr}}} \| \nabla \tilde{\vector{v}} \|_{0,K}, \quad \forall \tilde{\vector{v}} \in \tilde{\mathcal{X}}.
\end{equation}
Also, using the assumption that $\Omega$ is a simply connected polyhedron, it follows from Theorem 2.2 
%
%
{\color{ok} of}
\cite{gomes2022mhm} (see also Lemma 4.1 of \cite{harder2016hybrid}) that there is a constant $C_\mathrm{rm}$ such that
\begin{equation} \label{eq:result_vrm}
    \| \vector{v}_\mathrm{rm} \|_{\mathcal{X}} \leq C_\mathrm{rm} \left(\sup_{\vector{\mu}_\mathrm{rm} \in \mathcal{X}_\mathrm{rm}^{\partial \TauH}} \frac{b(\vector{v}_\mathrm{rm}, \vector{\mu}_\mathrm{rm})}{\| \vector{\mu}_\mathrm{rm} \|_{\mathcal{M}}}\right), \quad \forall \vector{v}_\mathrm{rm} \in \mathcal{X}_\mathrm{rm}^{\TauH}.
\end{equation}

\begin{lemma}[Coercivity of $a(\cdot, \cdot)$] \label{lemma:coercivity_a}
    Assuming that $\Omega$ is a simply connected polyhedron, if  $\mathcal{X}_\mathrm{rm}^{\TauH} \subset \mathcal{X}_h$ and $\mathcal{X}_\mathrm{rm}^{\partial \TauH} \subset \mathcal{M}_h$, then the bilinear form $a(\cdot, \cdot)$ is coercive on $\mathcal{N}_h^b$.
\end{lemma}
\begin{proof}
    Let $\vector{v} \in \mathcal{N}_h^{b}$ be decomposed into
    \[
    \vector{v} = \vector{v}_\mathrm{rm} + \tilde{\vector{v}}, 
    \]
    with $\vector{v}_\mathrm{rm} \in \mathcal{X}_\mathrm{rm}^{\TauH}$ and $\tilde{\vector{v}} \in \tilde{\mathcal{X}}$.
    From the definition \eqref{eq:kernel_ph}, we have, in particular, that
    \[ 
    b(\vector{v}_\mathrm{rm}, \vector{l}) = -b(\tilde{\vector{v}}, \vector{l}), \quad \forall \vector{l} \in \mathcal{M}_h.
    \]
    Using this identity and the assumption $\mathcal{X}_\mathrm{rm}^{\partial \TauH} \subset \mathcal{M}_h$, it follows from \eqref{eq:result_vrm} and the continuity of $b(\cdot, \cdot)$ that
    \begin{equation} \label{eq:lemma_coercivity_aux1}
            \| \vector{v}_\mathrm{rm} \|_{\mathcal{X}} \leq C_\mathrm{rm} \sup_{\vector{l}_\mathrm{rm} \in \mathcal{X}_\mathrm{rm}^{\partial \TauH}} \frac{b(\vector{v}_\mathrm{rm},\vector{l}_\mathrm{rm})}{\| \vector{l}_\mathrm{rm} \|_{\mathcal{M}}}
            = C_\mathrm{rm} \sup_{\vector{l}_\mathrm{rm} \in \mathcal{X}_\mathrm{rm}^{\partial \TauH}} \frac{-b( \tilde{\vector{v}},\vector{l}_\mathrm{rm})}{\| \vector{l}_\mathrm{rm} \|_{\mathcal{M}}}
            \leq C_\mathrm{rm} \| \tilde{\vector{v}} \|_{\mathcal{X}}.
    \end{equation}
    Finally, combining the broken Korn and Poincaré inequalities \eqref{eq:broken_korn_local} with \eqref{eq:lemma_coercivity_aux1}, we conclude that for every $\vector{v} \in \mathcal{N}_h^{b}$
    \begin{equation} \label{eq:coercivity_a}
        a(\vector{v}, \vector{v}) = (2\mu \tensor{\varepsilon}(\tilde{\vector{v}}), \tensor{\varepsilon}(\tilde{\vector{v}}))_{\TauH} 
        \geq \frac{\mu}{( C_{\mathrm{Pc}} C_{\mathrm{Kr}} \diamOm)^2} \| \tilde{\vector{v}} \|^2_{\mathcal{X}} 
        \geq \frac{\mu}{(C_{\mathrm{Pc}} C_{\mathrm{Kr}})^2 (1+C_\mathrm{rm}^2)} \| \vector{v} \|^2_{\mathcal{X}},
    \end{equation}
    from which we conclude the coercivity of $a(\cdot,\cdot)$ with the constant
    \begin{equation}
        \tilde{\alpha} = \frac{\mu}{(C_{\mathrm{Pc}} C_{\mathrm{Kr}})^2 (1+C_\mathrm{rm}^2)}.
    \end{equation}
\end{proof}

Under the same conditions of Lemma \ref{lemma:coercivity_a}, the coercivity of $a(\cdot, \cdot)$ directly implies that
\begin{equation} \label{eq:coercivity_hdp_compressible}
	\ahdp((\vector{v},q), (\vector{v},q)) = a(\vector{v}, \vector{v}) + \frac{1}{\lambda}(q,q)_{\Omega} \geq \tilde{\alpha} \| \vector{v} \|_{\mathcal{X}}^2 + \frac{1}{\lambda} \| q \|_{0,\Omega}^2, \quad \forall (\vector{v}, q) \in 
    \mathcal{N}_h^{b} \times \mathcal{P}_h,
\end{equation}
%
from which we conclude that \ref{it:hdp_coercivity} holds with
\[
\alpha_{\hdp} = \frac{\sqrt{2}}{2} \min \left\{ \tilde{\alpha}, \frac{1}{\lambda} \right\}. 
\]

This approach, however, leads to a coercivity constant $\alpha_{\hdp}$ that goes to zero as $\lambda$ tends to infinity.
Since our goal is to obtain stable and convergent approximations in nearly-incompressible problems, it is imperative to prove \ref{it:hdp_coercivity} with a constant that is bounded away from zero when $\lambda$ goes to infinity.
This 
{\color{ok} requires} additional conditions 
{\color{ok} on} the spaces $\mathcal{X}_h$ and $\mathcal{P}_h$, which are addressed in the following lemma.

\begin{lemma} \label{lemma:coercivity_hdp}
    Under the same conditions of Lemma \ref{lemma:coercivity_a}, assume, additionally, that there is a positive constant $\tilde{\beta}$, \textcolor{ok}{independent of $\lambda$,} such that
	\begin{equation} \label{eq:coercivity_cond_hdp}
		\sup_{\vector{v} \in \mathcal{N}_h^{b}} \frac{c(\vector{v}, z)}{\| \vector{v} \|_{\mathcal{X}}} \geq \tilde{\beta} \| z \|_{0, \Omega}, \quad \quad \forall z \in \mathcal{P}_h,
	\end{equation}
	where $\mathcal{N}_h^{b}$ is the space defined in \eqref{eq:kernel_ph}.
    Then condition \ref{it:hdp_coercivity} holds with a constant $\alpha_{\hdp}$ bounded away from zero when $\lambda$ goes to infinity.
\end{lemma}

\begin{proof} 
    Let $(\vector{w}, z) \in \mathcal{N}_h^{\hdp}$ be an arbitrary element. 
	Using Lemma A.42 from \cite{ern2004theory}, condition \eqref{eq:coercivity_cond_hdp} implies that there is $\vector{y} \in \mathcal{N}_h^{b}$ such that
\begin{equation*}
c(\vector{w}, q) = - c(\vector{y}, q) \quad \forall q \in \mathcal{P}_h \quad \textrm{ and } \quad \| \vector{y} \|_{\mathcal{X}} \leq \frac{1}{\tilde{\beta}} \sup_{q \in \mathcal{P}_h} \frac{c(\vector{w}, q)}{\| q \|_{0, \Omega}}.
\end{equation*}
In particular, we have
\begin{equation} \label{eq:lemma_coercivity_hdp_aux2}
\begin{split}
\| \vector{y} \|_{\mathcal{X}} & \leq \frac{1}{\tilde{\beta}} \sup_{q \in \mathcal{P}_h} \frac{a(\vector{w}, \vector{0}) + c(\vector{0}, z) - c(\vector{w}, q) + \dfrac{1}{\lambda}(q,z)_{\Omega} - \dfrac{1}{\lambda}(q,z)_{\Omega}}{\|(\vector{0},q) \|_\mathcal{H}} \\
& \leq \frac{1}{\tilde{\beta}} \left( \sup_{(\vector{v}, q) \in \mathcal{N}_h^{\hdp}} \frac{\ahdp((\vector{w}, z),(\vector{v}, q))}{\|(\vector{v}, q) \|_{\mathcal{H}}} + \dfrac{1}{\lambda} \| z \|_{0, \Omega} \right).		
\end{split} 
\end{equation}
Consider now the auxiliary space $\mathcal{N}_h^{b,c} \subset \mathcal{N}_h^{b}$ defined by
\[
\mathcal{N}_h^{b,c} = \{ \vector{v} \in \mathcal{N}_h^{b} \, : \, c(\vector{v}, q) = 0, \; \forall q \in \mathcal{P}_h \}. 
\]
Using the continuity of $a(\cdot, \cdot)$, the coercivity result \eqref{eq:coercivity_a}, and noticing that $\vector{w} + \vector{y}$ belongs to $\mathcal{N}_h^{b,c}$, we verify that
\begin{equation} \label{eq:lemma_coercivity_hdp_aux3}
\begin{split}
		\| \vector{w} + \vector{y} \|_{\mathcal{X}} & \leq \frac{1}{\tilde{\alpha}} \sup_{\vector{v} \in \mathcal{N}_h^{b,c}} \frac{a(\vector{w} + \vector{y}, \vector{v})}{\| \vector{v} \|_{\mathcal{X}}} \\ 
		& = \frac{1}{\tilde{\alpha}} \sup_{\vector{v} \in \mathcal{N}_h^{b,c}} \frac{a(\vector{w} + \vector{y}, \vector{v}) + c(\vector{v},z) - c(\vector{w}, 0) + \frac{1}{\lambda}(0, z)_{\Omega}}{\| (\vector{v}, 0) \|_{\mathcal{H}}} \\
		& \leq \frac{1}{\tilde{\alpha}} \left( \sup_{(\vector{v}, q) \in \mathcal{N}_h^{\hdp}} \frac{\ahdp((\vector{w}, z),(\vector{v}, q))}{\|(\vector{v}, q) \|_{\mathcal{H}}} + 2\mu \|  \vector{y} \|_{\mathcal{X}} \right).
	\end{split}
	\end{equation}
	Then, by combining \eqref{eq:lemma_coercivity_hdp_aux2} with \eqref{eq:lemma_coercivity_hdp_aux3}, it follows from the triangular inequality that
	\begin{equation} \label{eq:lemma_coercivity_hdp_aux4}
		\| \vector{w} \|_{\mathcal{X}} \leq \left( \frac{\tilde{\alpha} + \tilde{\beta} + 2\mu}{ \tilde{\alpha} \tilde{\beta}} \right) \sup_{(\vector{v}, q) \in \mathcal{N}_h^{\hdp}} \frac{\ahdp((\vector{w}, z),(\vector{v}, q))}{\|(\vector{v}, q) \|_{\mathcal{H}}} + \left( \frac{2\mu + \tilde{\alpha}}{\tilde{\alpha} \tilde{\beta} \lambda} \right) \| z \|_{0, \Omega}.
	\end{equation}
	Next, using once more \eqref{eq:coercivity_cond_hdp} and the continuity for  $a(\cdot, \cdot)$, we get
	\begin{equation*} 
	\begin{split}
		\| z \|_{0, \Omega} & \leq \frac{1}{\tilde{\beta}} \sup_{\vector{v} \in \mathcal{N}_h^{b}} \frac{c(\vector{v}, z)}{\| \vector{v} \|_{\mathcal{X}}} \\
        & =
        \frac{1}{\tilde{\beta}} \sup_{\vector{v} \in \mathcal{N}_h^{b}} \frac{a(\vector{w}, \vector{v}) - a(\vector{w}, \vector{v}) + c(\vector{v}, z) - c(\vector{w}, 0) + \frac{1}{\lambda}(0, z)_{\Omega}}{\| (\vector{v}, 0) \|_{\mathcal{H}}}
        \\
		& \leq \frac{1}{\tilde{\beta}} \left( \sup_{(\vector{v}, q) \in \mathcal{N}_h^{\hdp}} \frac{\ahdp((\vector{w}, z),(\vector{v}, q))}{\|(\vector{v}, q) \|_{\mathcal{H}}} + 2\mu \| \vector{w} \|_{\mathcal{X}} \right).
	\end{split} 
	\end{equation*} 
	By multiplying the previous inequality by $\tilde{\beta}/(4 \mu)$ and summing it to \eqref{eq:lemma_coercivity_hdp_aux4}, we then get
	\[\frac{\| \vector{w} \|_{\mathcal{X}}}{2} + \frac{\tilde{\beta} \| z \|_{0}}{4 \mu} \leq \left( \frac{\tilde{\alpha}+ \tilde{\beta} + 2\mu}{\tilde{\alpha} \tilde{\beta}} + \frac{1}{4 \mu} \right) \sup_{(\vector{v}, q) \in \mathcal{N}_h^{\hdp}} \frac{\ahdp((\vector{w}, z),(\vector{v}, q))}{\|(\vector{v}, q) \|_{\mathcal{H}}} + \left( \frac{\tilde{\alpha} + 2\mu}{\tilde{\alpha} \tilde{\beta} \lambda} \right) \| z \|_{0}. \]
	For $\lambda$ sufficiently large, say 
	\[ \lambda \geq \frac{8\mu(\tilde{\alpha} + 2\mu)}{\tilde{\alpha} \tilde{\beta}^2} = \frac{1}{\epsilon}, \]
	the previous relation implies that
	\[ 
    \sup_{(\vector{v}, q) \in \mathcal{N}_h^{\hdp}} \frac{\ahdp((\vector{w}, z),(\vector{v}, q))}{\|(\vector{v}, q) \|_{\mathcal{H}}} \geq \alpha_{\hdp} \| (\vector{w}, z) \|_{\mathcal{H}}, 
    \]
	with $\alpha_{\hdp}$ given by
	\[
    \alpha_{\hdp} = \min \left\{\frac{1}{2}, \frac{1}{8\mu}  \right\} \left( \frac{\tilde{\alpha}+ \tilde{\beta} + 2\mu}{\tilde{\alpha} \tilde{\beta}} + \frac{1}{4 \mu} \right)^{-1}.
    \]
	If $\lambda$ is less than $1/\epsilon$, the same procedure used in \eqref{eq:coercivity_hdp_compressible} furnishes
	\[
    \alpha_{\hdp} = \frac{\sqrt{2}}{2} \min \left\{ \tilde{\alpha}, \frac{\tilde{\alpha}\tilde{\beta}^2}{8\mu(\tilde{\alpha} + 2\mu)} \right\}.
    \]
    {\color{ok} From those two cases, we conclude that \ref{it:hdp_coercivity} is verified with a constant $\alpha_{\hdp}$ that, despite depending on $\lambda$, does not deteriorate as $\lambda$ goes to infinity.}
\end{proof}

As a final result of this section, we show that condition \eqref{eq:coercivity_cond_hdp} is easily satisfied if we construct $\mathcal{H}_h = \mathcal{X}_h \times \mathcal{P}_h$ based on 
{\color{ok} stable} Stokes spaces.

\begin{lemma} \label{lemma:coercivity_hdp_stokes}
	Consider $\mathcal{U}_h \subset H^1_{0, \Gamma_{D}}(\Omega, \mathbb{R}^{\nd})$, $\mathcal{P}_h \subset \mathcal{P}$ and $\mathcal{M}_h \subset \mathcal{M}$ finite dimensional subspaces.
	Assume that $\, \mathcal{U}_h \times \mathcal{P}_h$ form an inf-sup stable combination for the Stokes problem, meaning that there is a constant $\beta_{S}$ such that the inequality
    \begin{equation} \label{eq:inf-sup_stokes}
        \sup_{\vector{v} \in \mathcal{U}_h} \frac{c(\vector{v}, q)}{\| \vector{v} \|_{1,\Omega}} \geq \beta_S \| q \|_{0,\Omega}
    \end{equation}
    is satisfied for every $q \in \mathcal{P}_h$.
	Then, {\color{ok} the} combination $\mathcal{X}_h \times \mathcal{M}_h \times \mathcal{P}_h$ satisfies \eqref{eq:coercivity_cond_hdp} for every $\mathcal{X}_h$ such that
	\begin{equation} \label{eq:coercivity_cond_hdp2}
		\mathcal{U}_h \subset \mathcal{X}_h.
	\end{equation}
\end{lemma}
\begin{proof}
	First notice that there is a constant $C$, depending only on $\Omega$ and $\Gamma_D$, such that
	\[
    \| \vector{v} \|_{\mathcal{X}} \leq C \| \vector{v} \|_{1,\Omega}, \quad \quad \forall \vector{v} \in H^1_{0,\Gamma_{D}}(\Omega, \mathbb{R}^d).
    \]
	Furthermore, from definitions \eqref{eq:def_M} and \eqref{eq:formulation_PH_b}, it follows that
    \begin{equation} \label{eq:property_b_H1}
    b(\vector{v}, \vector{l}) = 0 \quad \forall \vector{v} \in H^1_{0,\Gamma_D}(\Omega, \mathbb{R}^\nd), \quad \forall \vector{l} \in \mathcal{M}.
    \end{equation}
    Therefore, since $\mathcal{U}_h \subset \mathcal{X}_h \cap H^1_{0,\Gamma_D}(\Omega, \mathbb{R}^d)$, we conclude that $\mathcal{U}_h \subset \mathcal{N}_h^{b}$.
	Finally, from the hypothesis that $\, \mathcal{U}_h \times \mathcal{P}_h$ is Stokes stable, it follows that
	\[
    \beta_S \| w \|_{0, \Omega} \leq \sup_{\vector{v} \in \mathcal{U}_h} \frac{c(\vector{v}, w)}{\| \vector{v} \|_{1,\Omega}} \leq C \sup_{\vector{v} \in \mathcal{U}_h} \frac{c(\vector{v}, w)}{\| \vector{v} \|_{\mathcal{X}}} \leq C \sup_{\vector{v} \in \mathcal{N}^{b}_h} \frac{c(\vector{v}, w)}{\| \vector{v} \|_{\mathcal{X}}}, 
    \]
	from which we conclude \eqref{eq:coercivity_cond_hdp} with $\tilde{\beta} = \beta_S/C$. 
\end{proof}

\begin{remark}[Coercivity for the PH and AP methods] \label{rk:coercivity}
    It is possible to show that the same condition $\mathcal{X}_\mathrm{rm}^{\partial \TauH} \subset \mathcal{M}_h$, established in Lemma \ref{lemma:coercivity_a}, is also needed to prove the Coercivity of $a_{\ph}(\cdot, \cdot)$ in the PH method.
    This condition prevents the use of the lowest-order approximation for the Lagrange multiplier, in which $\mathcal{M}_h$ is made of piece-wise constant functions on the faces (edges for two-dimensional domains) of $\partial K$.
    On the other hand, as proven in \cite{acharya2022primal}, such piece-wise constant space is enough to ensure coercivity for $a_{\ap}(\cdot, \cdot)$, allowing the lowest-order approximation in the AP method.
\end{remark}

\section{Compatible spaces on triangular and quadrilateral partitions} \label{sec:spaces}

{\color{ok}
The results of Section \ref{sec:hdp_analysis} provide a straightforward recipe for constructing HDP-compatible approximation spaces.
Although these results hold for both two- and three-dimensional domains with different mesh topologies, we illustrate their application only in the two-dimensional case.
Therefore, in this section, we shall assume a two-dimensional domain partitioned into triangular or convex quadrilateral elements.
In this context, we construct stable HDP approximations by combining stable spaces for the PH and Stokes formulations.
Extension to three dimensions follows the same ideas, but selecting suitable stable spaces for the PH and Stokes formulations can become more challenging.
In Remark \ref{rk:hdp3d}, we comment on some options assuming a tetrahedral mesh.
}

We start by defining a reference element $\hat{K}$ as
\begin{subequations}
    \begin{align*}
        & \hat{K} = \{ (\hat{x}_1,\hat{x}_2) \in \mathbb{R}^2 \, : \, \hat{x}_1, \hat{x}_2 > 0 \textrm{ and } \hat{x}_1+\hat{x}_2 < 1 \} & \textrm{ when } \TauH \textrm{ is a triangular partition}, \\
       & \hat{K} = (-1,1) \times (-1,1) & \textrm{ when } \TauH \textrm{ is a quadrilateral partition},
    \end{align*}
\end{subequations}
which can be mapped to each geometrical element $K \in \TauH$ through the bijective mappings $F_K$, i.e.,
\[
K = F_K(\hat{K}). 
\]
As a general convention, we denote by $\hat{\vector{x}} = (\hat{x}_1, \hat{x}_2)$ the spatial coordinates on the reference element $\hat{K}$.

For triangular meshes, the mappings $F_K$ are always affine transformations.
On general convex quadrilateral partitions, the transformation $F_K$ is affine if and only if $K$ is a parallelogram element \citep{boffi2002quadrilateral}.
The mappings $F_K$ need to be bilinear to cover general convex quadrilateral elements.
We say that a quadrilateral mesh is affine if every transformation is affine (i.e., every element is a parallelogram) and bilinear if at least one transformation is bilinear and non-affine (i.e., there is at least one non-parallelogram quadrilateral element).

Denote by $h_K$ the diameter of an element $K \in \TauH$ and define the mesh parameter $h$ as
\[
h = \max_{K \in \TauH} \{ h_K \}. 
\]
In order to prove stability and convergence for the HDP method, we need to set some regularity conditions in addition to the assumption that $\TauH$ has no hanging nodes.
For each element $K \in \TauH$, define
\begin{subequations}
    \begin{align*}
        \rho_K & = \textrm{diameter of the circle inscribed in } K, & \textrm{ for } K \textrm{ triangular,} \\
        \rho_K & = 2 \min_{1 \leq i \leq 4} \{ \textrm{diameter of the circle inscribed in } T_i \} & \textrm{ for } K \textrm{ quadrilateral,}
    \end{align*}
\end{subequations}
with $T_i$ denoting the four possible triangles obtained by choosing three of the vertices of a convex quadrilateral element $K$.
For both triangular and quadrilateral meshes, it follows  that  \citep[Appendix A]{girault2012finite} there are constants $C$, independent of the element geometry, satisfying
\begin{equation} \label{eq:limJF}
    \sup_{\hat{\vector{x}} \in \hat{K}} JF_K(\hat{\vector{x}}) \leq C h_K^2 \quad \textrm{ and } \quad \sup_{\vector{x} \in K} JF^{-1}_K(\vector{x}) \leq C \frac{1}{\rho_K^2},
\end{equation}
where $JF_K = \textrm{det} (DF_K)$, with $DF_K$ denoting the Jacobian matrix of $F_K$.
We define the shape-constant of a mesh by
\[
\theta = \max_{K \in \TauH} \dfrac{h_K}{\rho_K}, 
\]
and say that a family of meshes $\{ \TauH \}_h$ (indexed by their mesh parameter) is regular if their shape-constants can be uniformly bounded, say by a constant $\theta_{\max}$.

Under the assumption that the triangular and quadrilateral meshes are regular, Sections \ref{sec:inf-sup_PH} and \ref{sec:inf-sup_stokes} present inf-sup stable spaces for the standard PH method and the mixed Stokes method, respectively.
Such spaces are combined in Section \ref{sec:spaces_HDP} to generate stable and optimally convergent spaces for the HDP method.

\subsection{Stable spaces for the standard PH method} \label{sec:inf-sup_PH}

{ \color{ok}
By choosing finite-dimensional subspaces $\hat{\mathcal{X}} \subset H^1(\hat{K}, \mathbb{R}^2)$ and $\hat{\mathcal{M}} \subset L^2(\partial \hat{K}, \mathbb{R}^2)$ over the reference element $\hat{K}$, we  obtain
 local spaces
 $\mathcal{X}_K$ and $\mathcal{M}_K$  by mapping $\hat{\mathcal{X}}$ and $\hat{\mathcal{M}}$ to the geometrical elements $K$, according to
\begin{equation*}
    \mathcal{X}_K = \{ \vector{v} \in H^1(K, \mathbb{R}^2) \, : \, \vector{v} = \hat{\vector{v}} \circ F_K^{-1}, \; \hat{\vector{v}} \in \hat{\mathcal{X}} \},
\end{equation*}
\begin{equation*}
    \mathcal{M}_K = \{ \vector{l} \in L^2(\partial K, \mathbb{R}^2) \, : \, \vector{l} = \hat{\vector{l}} \circ F_K^{-1}, \; \hat{\vector{l}} \in \hat{\mathcal{M}} \}.
\end{equation*}
Then,
the global displacement and multiplier approximation spaces are constructed as follows
\begin{equation} \label{eq:constr_Xh}
    \mathcal{X}_h = \{ \vector{v} \in L^2(\Omega, \mathbb{R}^2) \, : \, \vector{v}\vert_K \in \mathcal{X}_K, \; \forall K \in \TauH \},
\end{equation}
\begin{multline} \label{eq:constr_Mh}
    \mathcal{M}_h = \left\{  \vector{l} \in \prod_{K \in \TauH} \mathcal{M}_K \, : \, \vector{l} \vert_{\Gamma_N} = 0 \, \textrm{ and } \, \vector{l} \vert_{\partial K_1} + \vector{l} \vert_{\partial K_2} = 0 \textrm{ on } \partial K_1 \cap \partial K_2 \right. \\
    \left. \textrm{for every pair of neighboring elements } K_1,K_2 \in \TauH \right\}.
\end{multline}

On two-dimensional domains, we take $\hat{\mathcal{M}} = E_r(\partial \hat{K}, \mathbb{R}^2)$ for some integer $r \geq 1$.
Here, $E_r(\partial \hat{K}, \mathbb{R}^2)$ is the space of vector-valued functions whose components are polynomials of degree at most $r$ on each edge of $\partial \hat{K}$, as defined in Section \ref{sec:notation}.
The combination of this space with $\hat{\mathcal{X}} = P_{r+1}(\hat{K}, \mathbb{R}^2)$ for $r$ even on triangular meshes or with $\hat{\mathcal{X}} = P_{r+1}(\hat{K}, \mathbb{R}^2)$ for any $r\geq 1$ on quadrilateral meshes is not stable \cite{raviart1977primal}.
To achieve stability, it is enough to enrich the displacement polynomial spaces with additional bubbles, as follows.

For triangular meshes set $P_r^{\ph}(\hat{K}, \mathbb{R})$ as $P_r(\hat{K}, \mathbb{R})$ if $r$ is odd and as the space spanned by the polynomials of $P_r(\hat{K}, \mathbb{R})$ plus the additional function
}
\begin{equation*}
    v_0^T(\hat{x}_1, \hat{x}_2) = (\phi_1 - \phi_2)(\phi_2 - \phi_3)(\phi_3 - \phi_1)[ (\phi_1 \phi_2)^{\frac{r-2}{2}} + (\phi_2 \phi_3)^{\frac{r-2}{2}} + (\phi_3 \phi_1)^{\frac{r-2}{2}}]
\end{equation*}
if $r$ is even.
In the definition of $v_0^T$, the functions $\phi_i$ are the linear polynomials 
\begin{equation} \label{eq:barycentric_coordinates}
    \phi_1(\hat{x}_1, \hat{x}_2) = 1 - \hat{x}_1 - \hat{x}_2, \quad \phi_2(\hat{x}_1, \hat{x}_2) = \hat{x}_1, \quad \phi_3(\hat{x}_1, \hat{x}_2) = \hat{x}_2.
\end{equation}
On quadrilateral meshes, denote by $Q_r^{\ph}(\hat{K}, \mathbb{R})$ the space spanned by the polynomials from $Q_r(\hat{K}, \mathbb{R})$ and the additional function
\begin{subequations} \label{eq:v0_quadrilateral}
    \begin{align*}
        v_0^Q(\hat{x}_1, \hat{x}_2) & = [(1 + \hat{x}_1)(1-\hat{x}_1) - (1+ \hat{x}_2)(1-\hat{x}_2)] \\
        & \hspace{3cm} [ ((1+\hat{x}_1)(1-\hat{x}_1))^{\frac{r-1}{2}} + ((1+\hat{x}_2)(1-\hat{x}_2))^{\frac{r-1}{2}} ], & \mbox{ for } r \mbox{ odd}, \\[0.5cm]
        v_0^Q(\hat{x}_1, \hat{x}_2) & = \hat{x}_1 \hat{x}_2 [(1 + \hat{x}_1)(1-\hat{x}_1) - (1+ \hat{x}_2)(1-\hat{x}_2)] \\ 
        & \hspace{3cm} [ ((1+\hat{x}_1)(1-\hat{x}_1))^{\frac{r-2}{2}} + ((1+\hat{x}_2)(1-\hat{x}_2))^{\frac{r-2}{2}} ], & \mbox{ for } r \mbox{ even.}
    \end{align*}
\end{subequations}
The vector-valued spaces $P_r^{\ph}(\hat{K}, \mathbb{R}^2)$ and $Q_r^{\ph}(\hat{K}, \mathbb{R}^2)$ are defined using the same convention as in Section \ref{sec:notation}.

\begin{proposition} \label{eq:prop_stability_PH}
    For $r \geq 0$ an integer, consider $\mathcal{M}_h$ the space constructed through \eqref{eq:constr_Mh} with $\hat{\mathcal{M}} = E_r(\partial \hat{K}, \mathbb{R}^2)$, and $\mathcal{X}_h$ a space constructed according to  \eqref{eq:constr_Xh}.
    The resulting approximation space $\mathcal{X}_h \times \mathcal{M}_h$ satisfies \ref{it:PH_inf-sup} if
    \begin{subequations}
        \begin{align*}
            P^{\ph}_{r+1}(\hat{K}, \mathbb{R}^2) & \subset \hat{\mathcal{X}} & \textrm{ when } \TauH \textrm{ is a triangular partition}, \\
            Q^{\ph}_{r+1}(\hat{K}, \mathbb{R}^2) & \subset \hat{\mathcal{X}} & \textrm{ when } \TauH \textrm{ is a quadrilateral partition}.
        \end{align*}
    \end{subequations}
\end{proposition}
\begin{proof} 
    Analogous results for scalar-valued functions derive from Lemmas 6 and 8 of \cite{raviart1977primal}. 
    The extension for vector-valued functions follows naturally by applying the scalar results on each entry of the vector functions.
\end{proof}

\subsection{Stable spaces for Stokes} \label{sec:inf-sup_stokes}

In addition to the polynomial spaces previously defined, for $r \geq 2$ an integer, we also consider 
\[
P_r^S(\hat{K}, \mathbb{R}) = P_r(\hat{K}, \mathbb{R}) \oplus (\phi_1 \phi_2 \phi_3) \tilde{P}_{r-2}(\hat{K}, \mathbb{R}), 
\]
where the linear polynomials $\phi_i$ are defined as in \eqref{eq:barycentric_coordinates}.

The global velocity and pressure approximation spaces for the Stokes formulation are constructed according to
\begin{equation} \label{eq:constr_Uh}
    \mathcal{U}_h = \{ \vector{v} \in \mathcal{C}^0(\Omega, \mathbb{R}^2) \, : \, \vector{v}\vert_K \in \mathcal{U}_K, \; \forall K \in \TauH \},
\end{equation}
\begin{equation} \label{eq:constr_Ph}
    \mathcal{P}_h = \{ q \in \mathcal{P} \, : \, q\vert_K \in \mathcal{P}_K, \, \forall K \in \TauH \}.
\end{equation}
In \eqref{eq:constr_Uh}, the local spaces $\mathcal{U}_K$ are obtained by mapping functions defined over the reference element $\hat{K}$ to the geometrical elements $K$, i.e.,
\[
\mathcal{U}_K = \{ \vector{v} \in H^1(K, \mathbb{R}^2) \, : \, \vector{v} = \hat{\vector{v}} \circ F_K^{-1}, \;  \hat{\vector{v}}  \in \hat{\mathcal{U}} \},
\]
with $\hat{\mathcal{U}} \subset H^1(\hat{K}, \mathbb{R}^2)$ being a finite-dimensional subspace.
For the construction of $\mathcal{P}_h$, on the other hand, the space $\mathcal{P}_K$ is defined directly over the geometrical elements $K$ rather than being mapped from a reference element.

\begin{proposition}
    Let $r \geq 1$ be an integer, $\mathcal{P}_h$ the space constructed through \eqref{eq:constr_Ph} with $\mathcal{P}_K = P_r(K, \mathbb{R})$, and $\mathcal{U}_h$ a space constructed according to \eqref{eq:constr_Uh}.
    If we have
    \begin{subequations}
        \begin{align*}
            P_{r+1}^S(\hat{K}, \mathbb{R}^2) & \subset \hat{\mathcal{U}} & \textrm{ when } \TauH \textrm{ is a triangular partition}, \\
            Q_{r+1}(\hat{K}, \mathbb{R}^2) & \subset \hat{\mathcal{U}} & \textrm{ when } \TauH \textrm{ is a quadrilateral partition},
        \end{align*}
    \end{subequations}
    then the resulting approximation space $\, \mathcal{U}_h \times \mathcal{P}_h$ satisfies \eqref{eq:inf-sup_stokes}.
\end{proposition}
\begin{proof}
    A proof for this result can be found in \citep[Chapter II, Sections 2 and 3]{girault2012finite} or \citep[Sections 8.5.5 and 8.6.3]{boffi2013mixed}, for instance.
\end{proof}

\begin{remark}[Mapped vs. unmapped $\mathcal{P}_h$] \label{rk:mapped_vs_geo}
    It is also possible to construct the pressure approximation space by mapping $P_r(\hat{K}, \mathbb{R})$ from the reference element to the geometrical elements.
    Such construction is equivalent to directly set $\mathcal{P}_K = P_r(K, \mathbb{R})$ on the geometrical elements when the transformations $F_K$ are affine (triangular or parallelogram meshes).
    On general quadrilateral meshes, however, where the transformations $F_K$ can be non-affine, the mapped construction can lead to sub-optimal convergence orders {\color{ok}\citep{boffi2002quadrilateral,ABDALA-MADY-CORREA-2025}.}
    Section \ref{sec:exp_convergence_quad} illustrates this loss of optimality in the context of the HDP method.
\end{remark}

\subsection{Approximation spaces for the HDP method} \label{sec:spaces_HDP}

Considering $r \geq 1$ an integer, to simplify our presentation, we denote from now on:
\begin{itemize}
    \item $\mathcal{U}_h^{r+1}$ the space constructed through \eqref{eq:constr_Uh} by setting $\hat{\mathcal{U}} = P^S_{r+1}(\hat{K}, \mathbb{R}^2)$, for triangular meshes, and $\hat{\mathcal{U}} = Q_{r+1}(\hat{K}, \mathbb{R}^2)$, for the quadrilateral ones;
    \item $\mathcal{X}_h^{r+1}$ the space constructed through \eqref{eq:constr_Xh} by setting $\hat{\mathcal{X}} = P^{\ph}_{r+1}(\hat{K}, \mathbb{R}^2) \oplus (\phi_1 \phi_2 \phi_3) \tilde{P}_{r-1}(\hat{K}, \mathbb{R})$, for triangular meshes, and $\hat{\mathcal{X}} = Q^{\ph}_{r+1}(\hat{K}, \mathbb{R}^2)$, for the quadrilateral ones;
    \item $\mathcal{M}_h^r$ the space constructed through \eqref{eq:constr_Mh} setting $\hat{\mathcal{M}} = E_r(\partial K, \mathbb{R}^2)$;
    \item $\mathcal{P}_h^r$ the space defined in \eqref{eq:constr_Ph} with $\mathcal{P}_K = P_r(K, \mathbb{R})$.
\end{itemize}
\begin{corollary} \label{cor:convergence_rates}
    Let $\Omega \subset \mathbb{R}^2$ be a simply connected polygonal domain partitioned into regular triangular or quadrilateral meshes.
    The spaces $\mathcal{X}^{r+1}_h \times \mathcal{M}_h^r \times \mathcal{P}_h^r$, {\color{ok} for an integer }$r \geq 1$, as defined above, guarantee a unique solution $(\vector{u}_h, \vector{m}_h, p_h)$ for the discrete problem \eqref{eq:method_HDP}.
    Also, if the exact solution $(\vector{u}, \vector{m}, p)$ for \eqref{eq:formulation_HDP} is regular enough, i.e., $\vector{u} \in H^1(\Omega, \mathbb{R}^2) \cap H^{r+2}(\TauH, \mathbb{R}^2)$, and $2\mu \tensor{\varepsilon}(\vector{u}) +  p \tensor{I} \in H(\div, \Omega, \mathbb{S}) \cap H^r(\TauH, \mathbb{S})$, the following convergence  estimate holds
\begin{equation} \label{eq:estimate_all}
        \| \vector{u} - \vector{u}_h \|_{\mathcal{X}} + \| \vector{m} - \vector{m}_h \|_{\mathcal{M}} + \| p - p_h \|_{0, \Omega} \leq C h^{r+1} (\vert \vector{u} \vert_{r+2, \TauH} + \vert p \vert_{r+1, \TauH}), 
    \end{equation}
    where $C$ is a constant independent of the mesh parameter $h$, which remains bounded as the Lamé coefficient $\lambda$ goes to infinity.
\end{corollary}
\begin{proof}
    First, notice that $\mathcal{U}_h^{r+1} \subset \mathcal{X}_h^{r+1}$, $\mathcal{X}_\mathrm{rm}^{\TauH} \subset \mathcal{X}_h^{1}$, and $\mathcal{X}_\mathrm{rm}^{\partial \TauH} \subset \mathcal{M}_h^1$.
    Therefore, since $\mathcal{X}_h^{r+1} \times \mathcal{M}_h^{r}$ and $\mathcal{U}_h^{r+1} \times \mathcal{P}_h^{r}$ are inf-sup stable spaces for the PH and mixed Stokes methods, the analysis of Section \ref{sec:hdp_analysis} guarantees the existence and uniqueness of solution.

    A direct application of the Bramble-Hilbert lemma \citep{dupont1980polynomial} furnishes the following approximation result for the space $\mathcal{P}^r_h$
    \begin{equation*} 
        \inf_{q \in \mathcal{P}^r_h} \| p - q \|_{0, \Omega} \leq C h^{r+1} \vert p \vert_{r+1, \TauH}.
    \end{equation*}
    From Theorem 3 of \cite{arnold2002approximation} and applying the Bramble-Hilbert lemma again, the following $\mathcal{X}^{r+1}_h$ approximation result holds
    \begin{equation*} 
        \inf_{\vector{v} \in \mathcal{X}^{r+1}_h} \| \vector{u} - \vector{v} \|_{\mathcal{X}} \leq C h^{r+1} \vert \vector{u} \vert_{r+2, \TauH}.
    \end{equation*}
    Such a result is also presented in Lemma 3.2 of \cite{taraschi2022convergence}.
    Finally, the following approximation result on $\mathcal{M}^r_h$
    \begin{equation*} 
        \inf_{\vector{l} \in \mathcal{M}^r_h} \| \vector{m} - \vector{l} \|_{\mathcal{M}} \leq C h^{r+1} (\vert \vector{u} \vert_{r+2, \TauH} + \vert p \vert_{r+1, \TauH}),
    \end{equation*}
    derives from Lemma 4.1 of \cite{araya2025multiscale}, where the constant $C$ may depend on $\mu$ but not on $\lambda$.

    The convergence estimate \eqref{eq:estimate_all} is then obtained by combining \eqref{eq:firstbound_u_w_hdp} and \eqref{eq:firstbound_m_hdp} with the approximation results above.
    Since $\alpha_{\hdp}$ is bounded away from zero (Lemmas \ref{lemma:coercivity_hdp} and \ref{lemma:coercivity_hdp_stokes}), the constants appearing in \eqref{eq:firstbound_u_w_hdp} and \eqref{eq:firstbound_m_hdp} remain bounded for $\lambda$ going to infinity, and so does the constant $C$ appearing in \eqref{eq:estimate_all}.
\end{proof}

\begin{remark}[Alternative norm for $\mathcal{M}^r_h$] \label{rk:norm_M}
    Let $\vertiii{ \cdot }_{\partial \TauH}$ denote the following norm over $\ds \prod_{K \in \TauH} L^2(\partial K, \mathbb{R}^2)$
    \begin{equation} \label{eq:norm_alt_Mh}
        \vertiii{\vector{l}}_{\partial \TauH} = \left( \sum_{K \in \TauH} h_K \| \vector{l} \|^2_{0, \partial K}\right)^{1/2}.
    \end{equation}
    According to Remark 6 from \cite{raviart1977primal}, there is a constant $C$ (independent of $h$) such that
    \[
    \vertiii{\vector{l}}_{\partial \TauH} \leq C \| \vector{l} \|_{\mathcal{M}}, \quad \forall \vector{l} \in \mathcal{M}^r_h, \; r \geq 0. 
    \]
    Whenever the exact solution $\vector{m}$ of \eqref{eq:formulation_HDP} is regular enough, we will use \eqref{eq:norm_alt_Mh} instead of \eqref{eq:norm_M} to measure the approximation errors for the Lagrange multiplier.
\end{remark}

\begin{remark}[$L^2$ convergence for the displacement] \label{rk:l2_rate}
    Estimate \eqref{eq:estimate_all} immediately gives the following convergence result for the displacement
    \[
    \| \vector{u} - \vector{u}_h \|_{0,\Omega} \leq C h^{r+1} (\vert \vector{u} \vert_{r+2, \TauH} + \vert p \vert_{r+1, \TauH}),  
    \]
    leading to an $\mathcal{O}(h^{r+1})$ convergence order in the $L^2$ norm.
    Under the additional assumptions that $\vector{u} \in H^2(\Omega, \mathbb{R}^2)$ and $\Omega$ is convex, however, it is possible to improve this result using a duality argument \citep{raviart1977primal, harder2016hybrid}.
    This leads to an $\mathcal{O}(h^{r+2})$ convergence order for the displacement in $L^2$, as we shall verify in the numerical experiments of Section \ref{sec:exp_convergence}.
\end{remark}

\begin{remark}[Static condensation for the HDP method] \label{rk:condensation}
    Notice that, for the approximation spaces constructed in this Sections, $\mathcal{X}_h$ and $\mathcal{P}_h$ are made of functions with no inter-element continuity.
    This allows to statically condense most of the degrees of freedom related to $\vector{u}_h$ and $p_h$, as discussed in {\color{ok}Appendix \ref{ap:condensation}}.
\end{remark}

{\color{ok}
\begin{remark}[Comments on the three-dimensional case] \label{rk:hdp3d}
    Stable spaces for three-dimensional domains can be constructed following the ideas of the two-dimensional case.
    According to Section \ref{sec:hdp_analysis}, it suffices to pick a PH stable pair $\mathcal{X}_h \times \mathcal{M}_h$ and combine it to a stable Stokes pair $\mathcal{U}_h \times \mathcal{P}_h$ such that $\mathcal{U}_h \subset \mathcal{X}_h$.
    Assuming the simpler case of a tetrahedral mesh, the multiplier space can be constructed following \eqref{eq:constr_Mh} with $\hat{\mathcal{M}} = E_r(\partial \hat{K}, \mathbb{R}^3)$. 
    To obtain a PH stable pair, it is sufficient to follow \eqref{eq:constr_Xh} setting $\hat{\mathcal{X}} = P_{r+3}(\hat{K}, \mathbb{R}^3)$ \cite{gomes2022mhm}.
    For the Stokes space, we can proceed as in \eqref{eq:constr_Uh} and \eqref{eq:constr_Ph} taking $\hat{\mathcal{U}} = P_{r+3}(\hat{K}, \mathbb{R}^3)$ and $\mathcal{P}_K = P_r(K,\mathbb{R})$, which are proven stable in \cite[Chapter 2, Section 2.3]{girault2012finite}.
    Notice that in this construction, the displacement space is quite large.
    We conjecture that it is possible to reduce it by only adding some higher-dimensional bubble functions $\vector{v}_i \in P_{r+3}(\hat{K}, \mathbb{R}^3)$ to the space $P_{r+1}(\hat{K}, \mathbb{R}^3)$.
    This, however, demands further investigation.
    Also, the somewhat inflated displacement space is amended because most of its degrees of freedom may be statically condensed. 
\end{remark}
}

\section{Local 
{ H}({\rm{div}})-conforming 
stress recovery} \label{sec:post}

In this section, we propose a post-processing strategy to recover an approximated stress field from the HDP solution.
Such a strategy extends the  flux recovery method for the Darcy equation, discussed in \cite{chou2002flux, correa2022optimal}, to the linear elasticity problem.
Following Section \ref{sec:spaces}, we shall focus on two-dimensional problems partitioned into triangular or quadrilateral meshes.
{\color{ok}
In Section \ref{sec:post-extensions}, we comment on the extension to tetrahedral meshes, along side the possibility to apply the post-processing on other Finite Element methods.
}

To introduce the stress recovery, we first need to construct appropriate $H(\div)$-conforming tensorial approximation spaces.
The spaces employed in this work are based on the Raviart-Thomas family \citep{raviart1977mixed, boffi2013mixed} of vector-valued spaces, for triangular meshes, and on the Arnold-Boffi-Falk family \citep{arnold2005quadrilateral}, for quadrilateral partitions.
In this sense, each row of the tensors in our approximation spaces can be seen as an $H(\div)$-conforming vector function belonging to either the Raviart-Thomas or the Arnold-Boffi-Falk spaces.

Given $\hat{\tensor{\tau}} \in H(\div, \hat{K}, \mathbb{M})$ a tensor function on $\hat{K}$, we transform it to the geometrical elements $K$ using the Piola transformation, defining $\tensor{\tau}=\piola{\hat{\tensor{\tau}}}:K\to\mathbb{M}$ by
\begin{equation}
\label{eq:piolatau}
\tensor{\tau}(\vector{x})=\piola{\hat{\tensor{\tau}}}(\hat{\vector{x}}) = (JF_K(\hat{\vector{x}}))^{-1} \hat{\tensor{\tau}}(\hat{\vector{x}}) DF^T_K(\hat{\vector{x}}),
\end{equation}
where, once more, $\hat{\vector{x}} = F_K^{-1}(\vector{x})$ denotes the spatial coordinates in the reference element $\hat{K}$.
Setting $\vector{\tau}$ as in \eqref{eq:piolatau} and $\vector{v}(\vector{x}) = \hat{\vector{v}} \circ F_K^{-1}(\vector{x})$, for some function $\hat{\vector{v}} \in H^1(\hat{K}, \mathbb{R}^\nd)$, two key properties of the Piola transform for the construction of $H(\div)$-conforming subspaces are
\begin{equation} \label{eq:piola_property_div}
    \div \tensor{\tau}(\vector{x}) = \frac{1}{JF_K(\hat{\vector{x}})}  \hat{\div} \hat{\tensor{\tau}}(\hat{\vector{x}}),
\end{equation}
\begin{equation} \label{eq:piola_property_normal}
    \int_{\partial K} (\tensor{\tau} \vector{n}^{\partial K}) \cdot \vector{v} \des = \int_{\partial \hat{K}} (\hat{\tensor{\tau}} \vector{n}^{\partial \hat{K}}) \cdot \hat{\vector{v}} \, \hat{\des},
\end{equation}
as demonstrated in \cite[Section 2.1.3]{boffi2013mixed}, for instance.

The Raviart-Thomas based tensor space of index $r \geq 0$ on a triangular element $K$ is defined as $RT_r(K,\mathbb{M}) = \piola{(RT_r(\hat{K}, \mathbb{M}))}$, with
\[
RT_r(\hat{K}, \mathbb{M}) = \begin{bmatrix}
P_r(\hat{K}, \mathbb{R}^2) + \hat{\vector{x}} \tilde{P}_r(\hat{K}, \mathbb{R}) \\
P_r(\hat{K}, \mathbb{R}^2) + \hat{\vector{x}} \tilde{P}_r(\hat{K}, \mathbb{R})
\end{bmatrix}. 
\]
For a quadrilateral element $K$, the Arnold-Boffi-Falk based tensor space of index $r \geq 0$ is defined as $ABF_r(K,\mathbb{M}) = \piola{(ABF_r(\hat{K}, \mathbb{M}))}$, with
\[
ABF_r(\hat{K}, \mathbb{M}) =  \begin{bmatrix}
P_{r+2,r}(\hat{K}, \mathbb{R}) \times P_{r,r+2}(\hat{K}, \mathbb{R}) \\
P_{r+2,r}(\hat{K}, \mathbb{R}) \times P_{r,r+2}(\hat{K}, \mathbb{R})
\end{bmatrix}. 
\]
The global approximation spaces are constructed according to
\begin{equation} \label{eq:hdiv_space}
    \mathcal{S}^r_h = \{ \tensor{\tau} \in H(\div, \Omega, \mathbb{M}) \, : \, \tensor{\tau}\vert_{K} \in \mathcal{S}^r_K, \; \forall K \in \TauH \},
\end{equation}
where $\mathcal{S}^r_K$ is either $RT_r(K,\mathbb{M})$ or $ABF_r(K,\mathbb{M})$. 
Notice that we use the same notation for the RT-based spaces on triangles and the ABF-based spaces on quadrilaterals since the context will suffice to distinguish between them.

\subsection{The local post-processing} \label{sec:post-init}

For two-dimensional problems, assume that the HDP solution $(\vector{u}_h, \vector{m}_h, p_h)$ was obtained using either a triangular or quadrilateral mesh with the approximation space $\mathcal{X}_h^{r+1} \times \mathcal{M}_h^r \times \mathcal{P}_h^r$, $r \geq 1$, as defined in Section \ref{sec:spaces_HDP}. 
We denote by $\vector{t}_h$ the following traction approximation
\begin{subequations} \label{eq:th}
    \begin{align}
        \vector{t}_h & = \vector{m}_h && \textrm{over } \partial K \setminus \Gamma_N, \\
        \vector{t}_h\vert_{\partial K \cap \Gamma_N} & = \vector{t}_N && \textrm{whenever } \partial K \textrm{ intersects } \Gamma_N. 
    \end{align}
\end{subequations}
In this context, the post-processing strategy is defined through the element-wise local problems: {\it for each $K \in \TauH$, find $\tensor{\sigma}_{h,K} \in \mathcal{S}_K^r$ such that
\begin{subequations} \label{eq:stress_post-processing}
    \begin{align}
        \int_{\partial K} (\tensor{\sigma}_{h,K} \vector{n}^{\partial K}) \cdot \vector{l} \des & = \int_{\partial K} \vector{t}_h \cdot \vector{l} \des, & \forall \vector{l} \in E_r(\partial K, \mathbb{R}^2), \label{eq:stress_post-processing_a} \\
        \int_K \tensor{\sigma}_{h,K} : \tensor{\tau} \dex & = \int_K (2 \mu \tensor{\varepsilon}(\vector{u}_h) + p_h \tensor{I}) : \tensor{\tau}\dex , & \forall \tensor{\tau} \in \Psi_r(K, \mathbb{M}), \label{eq:stress_post-processing_b} \\
        \nonumber
        \mbox{with the additional equation}&  &
        \\
        \int_K \div \tensor{\sigma}_{h,K} \cdot \vector{\phi} \dex & = -\int_K \vector{f} \cdot \vector{\phi} \dex, & \forall \vector{\phi} \in \Phi_r(K, \mathbb{R}^2), 
        \label{eq:stress_post-processing_c}
        \\
        \nonumber \mbox{ if } K \mbox{ is a quadrilateral.} \qquad  &&
    \end{align}
\end{subequations}
}
In \eqref{eq:stress_post-processing_b}, the auxiliary space $\Psi_r(K, \mathbb{M})$ is given by
\begin{equation} \label{eq:def_psi_trig}
    \Psi_r(K,\mathbb{M}) = P_{r-1}(K, \mathbb{M})
\end{equation}
when $K$ is a triangle and by
\begin{equation} \label{eq:def_psi_quad}
    \Psi_r(K,\mathbb{M}) = \{ \tensor{\tau}(\vector{x}) = (\hat{\tensor{\tau}} \circ F_K^{-1}(\vector{x})) \, DF_K^{-1}(\vector{x}) \, : \, \hat{\tensor{\tau}}_i \in P_{r-1,r}(\hat{K}) \times P_{r,r-1}(\hat{K}), \; i=1,2 \}
\end{equation}
if $K$ is a convex quadrilateral.
For quadrilateral elements, Eq. \eqref{eq:stress_post-processing_c} makes use of the additional space
\begin{equation}
    \Phi_r(K, \mathbb{R}^2) = \{ \vector{\phi}(\vector{x}) = \hat{\vector{\phi}} \circ F_K^{-1}(\vector{x}) \, : \, \hat{\vector{\phi}}_i \in \textrm{span} \{ \hat{x}_1^{r+1} \hat{x}_2^{s}, \; \hat{x}_1^{s} \hat{x}_2^{r+1}, \; s=0,1, \dots , r \}, \; i=1,2 \}.
\end{equation}
In such definitions, the sub-index $i$ denotes either the $i$-th row of a tensor function or the $i$-th entry of a vector function.

The well-posedness of the local problems follows directly from the definition of the degrees of freedom for the Raviart-Thomas and Arnold-Boffi-Falk spaces (see \cite{boffi2013mixed} and \cite{arnold2005quadrilateral} for more details).
From the definition of space $\mathcal{M}_h^r$ and the property \eqref{eq:piola_property_normal}, Eq. \eqref{eq:stress_post-processing_a} guarantees that the global stress approximation
\[
\tensor{\sigma}_h\vert_K = \tensor{\sigma}_{h,K}, \quad \forall K \in \TauH 
\]
has {\color{ok}continuous normal component} across inter-element boundaries \citep{boffi2013mixed, correa2022optimal}. 
Therefore, $\tensor{\sigma}_h$ belongs to $H(\div, \Omega, \mathbb{M})$ despite being computed through independent local problems.
The following Lemma shows that $\tensor{\sigma}_h$ is also locally equilibrated.
\begin{lemma}[Local equilibrium] \label{lemma:conservation}
    For a given element $K \in \TauH$, consider the space $\mathcal{R}_K^{r}$ defined as
    \begin{equation*}
        \mathcal{R}_K^{r} = \{ \vector{v} = \hat{\vector{v}} \circ F_K^{-1} \, : \, \hat{\vector{v}} \in \hat{\mathcal{R}} \},
    \end{equation*}
    with $\hat{\mathcal{R}} = P_{r}(\hat{K}, \mathbb{R}^2)$, when $K$ is a triangle, and $\hat{\mathcal{R}} = R_{r}(\hat{K}, \mathbb{R}^2)$, for $K$ quadrilateral.
    We recall that, as defined in Section \ref{sec:notation}, $R_{r}(\hat{K}, \mathbb{R}^2)$ stands for the space of vector-valued functions spanned by the monomials of $Q_{r+1}(D, \mathbb{R}^2)$ except for $[x_1^{r+1} x_2^{r+1}, 0]^t$ and $[0, x_1^{r+1} x_2^{r+1}]^t$
    It follows that the approximated stress $\tensor{\sigma}_h$ computed through \eqref{eq:stress_post-processing} satisfies
    \begin{equation} \label{eq:local_conservative}
        -\int_K \div \tensor{\sigma}_h \cdot \vector{v} \dex = \int_{K} \vector{f} \cdot \vector{v} \dex, \quad \forall \vector{v} \in \mathcal{R}_K^{r}, \; \forall K \in \TauH.
    \end{equation}
    In particular, taking $\vector{v} = (1,0)$ and $\vector{v} = (0,1)$, we get by the divergence theorem that
    \begin{equation} \label{eq:local_conservative_2}
        -\int_{\partial K}  \tensor{\sigma}_{h,i} \cdot \vector{n}^{\partial K} \des = \int_{K} \vector{f}_i \dex, \; \forall K \in \TauH,
    \end{equation}
    where the sub-index $i$ denotes the $i$-th row of $\tensor{\sigma}_h$ or the $i$-th entry of $\vector{f}$, meaning that the traction field is at equilibrium on each element.
\end{lemma}
\begin{proof}
    Given an element $K \in \TauH$, consider ${\color{ok}\mathcal{Z}_{K}^{r}} = \mathcal{R}_{K}^{r}$, when $K$ is a triangle and
    \begin{equation*}
        \mathcal{Z}_K^{r} = \{ \vector{v} = \hat{\vector{v}} \circ F_K^{-1} \, : \, \hat{\vector{v}} \in Q_{r}(\hat{K}, \mathbb{R}^2) \}
    \end{equation*}
    for $K$ a quadrilateral element.
    Take $\vector{v} \in \mathcal{X}_h^{r+1}$ such that $\vector{v}\vert_K$ belongs to $\mathcal{Z}_K^{r}$ and vanishes on the remaining elements.
    For every $\vector{v}$ in $\mathcal{Z}^{r}_K$, it follows from \eqref{eq:method_HDP_a} that
    \begin{equation*}
        \int_K 2 \mu \tensor{\varepsilon}(\vector{u}_h) : \tensor{\varepsilon}(\vector{v}) \dex + \int_{K} p_h \div \vector{v} \dex - \int_{\partial K} \vector{m}_h \cdot \vector{v} \des = \int_{K} \vector{f} \cdot \vector{v} \dex + \int_{\partial K \cap \Gamma_N} \vector{t}_N \cdot \vector{v} \des,
    \end{equation*}
    where the last term in the previous identity is included only if $\partial K$ intersects the Neumann boundary.
    Since $\vector{v}\vert_{\partial K} \in E_r(\partial K, \mathbb{R}^2)$ for $\vector{v}$ in $\mathcal{Z}_K^r$, the Green's formula and equations \eqref{eq:th} and \eqref{eq:stress_post-processing_a} lead to
    \begin{equation*}
        \int_K (2 \mu \tensor{\varepsilon}(\vector{u}_h) + p_h \tensor{I}): \tensor{\varepsilon}(\vector{v}) \dex - \int_{K} \tensor{\sigma}_h : \nabla \vector{v} \dex - \int_{K} \div \tensor{\sigma}_h \cdot \vector{v} \dex = \int_{K} \vector{f} \cdot \vector{v} \dex, \quad \forall \vector{v} \in \mathcal{Z}^{r}_K.
    \end{equation*}
    Because $(2 \mu \tensor{\varepsilon}(\vector{u}_h) + p_h \tensor{I})$ is symmetric, we can replace $\tensor{\varepsilon}(\vector{v})$ with $\nabla \vector{v}$ in the first integral of the previous equation.
    Now, from the fact that $ \nabla \vector{v} \in \Psi_{r}(K, \mathbb{M})$ for every $\vector{v} \in \mathcal{Z}^{r}_K$, eq. \eqref{eq:stress_post-processing_b} implies that 
    \begin{equation} \label{eq:lemma_conservation_aux1}
        -\int_{K} \div \tensor{\sigma}_h \cdot \vector{v} \dex = \int_{K} \vector{f} \cdot \vector{v} \dex, \quad \forall \vector{v} \in \mathcal{Z}^{r}_K.
    \end{equation}
    For triangular elements, equation \eqref{eq:lemma_conservation_aux1} gives the desired result. 
    For quadrilateral elements, the final result is obtained by combining \eqref{eq:lemma_conservation_aux1} with \eqref{eq:stress_post-processing_c}.
\end{proof}

One may notice that in the construction of $\mathcal{S}_h^r$, the symmetry of the tensor elements is not enforced.
Therefore, the stress approximation $\tensor{\sigma}_h$ obtained through \eqref{eq:stress_post-processing} is not point-wise symmetric.
However, the next Lemma shows that a weaker result holds, where the symmetry of the stress tensor is satisfied in a variational sense.

\begin{lemma}[A symmetry result] \label{lemma:ht_symmetry}
	Defining the scalar quantity $\mathrm{asym}(\tensor{\tau}) = \tensor{\tau}_{1,2} - \tensor{\tau}_{2,1}$ as a measure of the skew-symmetric part of $\vector{\tau} \in \mathbb{R}^{2 \times 2}$, it follows that
	\begin{equation} \label{eq:ht_symmetry}
	\int_{K} \mathrm{asym}(\tensor{\sigma}_h) \, q \dex = 0 \quad \forall q \in P_{r-1}(K, \mathbb{R}),
	\end{equation}
	In particular, taking $q=1$, we conclude that
	\[
    \int_{K} \mathrm{asym}(\tensor{\sigma}_h) \dex = 0, 
    \]
    implying that the stress approximation is element-wise symmetric on average.
\end{lemma}
\begin{proof}
    First, we verify that the space
    \begin{equation} \label{eq:aux_symmetry1}
        P_{r-1}(K, \mathbb{R}) \begin{bmatrix}
	    0 & 1 \\
	    -1 & 0
	\end{bmatrix}   
    \end{equation}
	belongs to $\Psi_r(K, \mathbb{M})$ for every $r \geq 1$.
    From definition \eqref{eq:def_psi_trig}, this result is direct for triangular meshes.
    Lets now prove it for quadrilateral meshes.
    Let $\tensor{\tau}$ be a tensor-valued function in \eqref{eq:aux_symmetry1}.
    Since $F_K$ is bilinear, its Jacobian matrix has the form
    \[
    DF_K(\hat{x}) = \begin{bmatrix}
        d_1 + d_2 \hat{x}_2 & d_3 + d_2 \hat{x}_1 \\
        d_4 + d_5 \hat{x}_2 & d_6 + d_5 \hat{x}_1
    \end{bmatrix},
    \]
    where $d_i$ are constants, and
    \[
    \tensor{\tau} \circ F_K(\hat{\vector{x}}) \in Q_{r-1}(\hat{K}) \begin{bmatrix}
        0 & 1 \\
        -1 & 0
    \end{bmatrix}.
    \]
    Therefore, the function $\hat{\tensor{\tau}} := (\tensor{\tau} \circ F_K(\hat{\vector{x}})) DF_K(\hat{\vector{x}})$ belongs to 
    \[
    \begin{bmatrix}
        P_{r-1,r}(\hat{K}) & P_{r,r-1}(\hat{K}) \\
        P_{r-1,r}(\hat{K}) & P_{r,r-1}(\hat{K})
    \end{bmatrix},
    \]
    and from definition \eqref{eq:def_psi_quad} we conclude that $\tensor{\tau} \in \Psi(K, \mathbb{M})$.

    Once this inclusion is verified, the desired result follows directly from eq. \eqref{eq:stress_post-processing_b} taking $\tensor{\tau}$ in \eqref{eq:aux_symmetry1} and noticing that $(2 \mu \, \tensor{\varepsilon}(\vector{u}_h) + p_h \tensor{I})$ is symmetric.
\end{proof}

\subsection{Convergence of the post-processed stress} \label{sec:post-conv}

We now establish some convergence results for the stress recovery strategy \eqref{eq:stress_post-processing}.
To do so, we briefly introduce the standard projection operators over the spaces $S_h^r$ and establish their approximation properties in \ref{prop:projection}.
Next, Proposition \ref{prop:stress_pre} presents an auxiliary result, and Theorem \ref{theo:stress} proves optimal convergence orders for the stress approximation in the $H(\div, \Omega, \mathbb{M})$ norm.
\begin{proposition}[Projection operators] \label{prop:projection}
    For a tensor $\tensor{\omega} \in H(\div, \Omega, \mathbb{M})$ and an integer $r \geq 1$, the standard projection operators $\proj{\mathcal{S}_h^r} : H(\div, \Omega, \mathbb{M}) \rightarrow \mathcal{S}_h^r$ are defined by the local systems
    \begin{subequations} \label{eq:proj_H(div)}
    \begin{align}
    \int_{\partial K} (\proj{\mathcal{S}_h^r} \tensor{\omega} \, \vector{n}^{\partial K}) \cdot \vector{l} \des & = \int_{\partial K} (\tensor{\omega} \vector{n}^{\partial K}) \cdot \vector{l} \des & \forall \vector{l} \in E_r(\partial K, \mathbb{R}^2), \label{eq:proj_H(div)_a} \\
    \int_K \proj{\mathcal{S}_h^r} \tensor{\omega} : \tensor{\tau} \dex & = \int_K \tensor{\omega} : \tensor{\tau} \dex & \forall \tensor{\tau} \in \Psi_r(K, \mathbb{M}), \label{eq:proj_H(div)_b} \\
    \int_K \div \left(\proj{\mathcal{S}_h^r} \tensor{\omega}\right) \cdot \vector{\phi} \dex & = \int_K \div \tensor{\omega} \cdot \vector{\phi} \dex & \forall \vector{\phi} \in \Phi_r(K, \mathbb{R}^2), && \textrm{(only for } K \textrm{ quadrilateral)}. \label{eq:proj_H(div)_c}
    \end{align}
    \end{subequations}
    Also, if $\tensor{\omega} \in H^{r+1}(\Omega, \mathbb{M})$, $\div \tensor{\omega} \in H^{r+1}(\Omega, \mathbb{R}^2)$, and the mesh $\TauH$ is regular (as established in Section \ref{sec:spaces}), the following estimates hold
    \begin{equation} \label{eq:approx_proj_L2}
        \| \tensor{\omega} - \proj{\mathcal{S}^r_h} \tensor{\omega}  \|_{0,\Omega} \leq C h^{r+1}\vert\tensor{\omega} \vert_{r+1,\Omega},
    \end{equation}
    \begin{equation} \label{eq:approx_proj_div}
        \| \div (\tensor{\omega} - \proj{\mathcal{S}^r_h} \tensor{\omega})  \|_{0,\Omega} \leq C h^{r+1}\vert\div \tensor{\omega} \vert_{r+1,\Omega},
    \end{equation}
    where the constants $C$ are independent of the mesh parameter $h$.
\end{proposition}
\begin{proof}
    Analogous results for the vector Raviart-Thomas and Arnold-Boffi-Falk spaces are proven in \cite{boffi2013mixed, arnold2005quadrilateral}.
    The results for the tensor-valued spaces follow by applying the vector projections row-wise.
\end{proof}
\begin{proposition} \label{prop:stress_pre}
    Let $\vector{t} \in L^2(\partial K, \mathbb{R}^2)$, $\tensor{\xi} \in L^2(K, \mathbb{M})$, and $\vector{z} \in L^2(K, \mathbb{R}^2)$ be given functions.
    If $\tensor{\omega} \in \mathcal{S}_K^r$ ($r \geq 1$) satisfies 
    \begin{subequations} \label{eq:stress_aux}
    \begin{align}
        \int_{\partial K} (\tensor{\omega} \vector{n}^{\partial K}) \cdot \vector{l} \des & = \int_{\partial K} \tensor{t} \cdot \vector{l} \des & \forall \vector{l} \in E_r(\partial K, \mathbb{R}^2), \label{eq:stress_aux_a} \\
        \int_K \tensor{\omega} : \tensor{\tau} \dex & = \int_K \tensor{\xi} : \tensor{\tau} \dex & \forall \tensor{\tau} \in \Psi_r(K, \mathbb{M}), \label{eq:stress_aux_b} \\
        \int_K \div \tensor{\omega} \cdot \vector{\phi} \dex & = \int_K \vector{z} \cdot \vector{\phi} \dex & \forall \vector{\phi} \in \Phi_r(K, \mathbb{R}^2), && \textrm{(only for } K \textrm{ quadrilateral)}, \label{eq:stress_aux_c}
    \end{align}
    \end{subequations}
    then there is a constant $C$, independent of $h_K$, such that
    \begin{subequations}
        \begin{align}
            \| \tensor{\omega} \|_{0,K} & \leq C (h_k^{1/2} \| \vector{t} \|_{0, \partial K} + \| \tensor{\xi} \|_{0,K}) & \textrm{ for } K \textrm{ a triangle}, \\
            \| \tensor{\omega} \|_{0,K} & \leq C (h_k^{1/2} \| \vector{t} \|_{0, \partial K} + \| \tensor{\xi} \|_{0,K} + h_K \| \vector{z} \|_{0,K}) & \textrm{ for } K \textrm{ a convex quadrilateral}.
        \end{align}
    \end{subequations}
\end{proposition}
\begin{proof}
    For $K$ a triangle and $\mathcal{S}_K^r = RT_r(K, \mathbb{M})$, this result follows from Lemma 3.3 of \cite{chou2002flux} by applying it on each row of the tensor $\tensor{\omega}$.
    Similarly, Lemma 3.1 from \cite{correa2022optimal} proves the result for $K$ a convex quadrilateral and $\mathcal{S}_K^r = ABF_r(K, \mathbb{M})$.
\end{proof}

\begin{theorem}[$H(\div)$ convergence] \label{theo:stress}
    Let $\TauH$ be a regular partition into triangular or quadrilateral elements, and consider $(\vector{u}_h, \vector{m}_h, p_h)$ the solution of \eqref{eq:method_HDP} using the approximation spaces $\mathcal{X}_h^{r+1} \times \mathcal{M}_h^r \times \mathcal{P}_h^r$, $r \geq 1$, as defined in Section \ref{sec:spaces_HDP}.
    Assuming that the exact solutions are regular enough, the post-processed stress $\tensor{\sigma}_h$ obtained through \eqref{eq:stress_post-processing} satisfies
    \begin{equation} \label{eq:strees_rate}
        \| \tensor{\sigma} - \tensor{\sigma}_h \|_{H(\div)} \leq C h^{r+1} (\vert \tensor{\sigma} \vert_{r+1, \Omega} + \vert \div \tensor{\sigma} \vert_{r+1, \Omega} + \vert \vector{u} \vert_{r+2, \Omega} + \vert p \vert_{r+1, \Omega}),
    \end{equation}
    with $C$ being a constant independent of the mesh parameter $h$ that remains bounded as the Lamé parameter $\lambda$ goes to infinity.
\end{theorem}
\begin{proof}
    We begin by proving the $L^2$ estimate
    \begin{equation} \label{eq:theo_stress_aux1}
        \| \tensor{\sigma} - \tensor{\sigma}_h \|_{0, \Omega} \leq C h^{r+1} (\vert \tensor{\sigma} \vert_{r+1, \Omega} + \vert \vector{u} \vert_{r+1, \Omega} + \vert p \vert_{r+1,\Omega}).
    \end{equation}
    Subtracting \eqref{eq:stress_post-processing} from \eqref{eq:proj_H(div)} and using the constitutive equation
    \[
    \tensor{\sigma} = 2 \mu \tensor{\varepsilon}(\vector{u}) + p \tensor{I},
    \]
    we get that 
    \begin{subequations}
    \begin{align*}
        \int_{\partial K} (\proj{\mathcal{S}_h^r} \tensor{\sigma} - \tensor{\sigma}_{h,K})\vector{n}^{\partial K} \cdot \vector{l} \des & = \int_{\partial K} (\tensor{\sigma} \vector{n}^{\partial K} - \vector{t}_h) \cdot \vector{l} \des & \forall \vector{l} \in E_r(\partial K, \mathbb{R}^2), \\
        \int_K (\proj{\mathcal{S}_h^r} \tensor{\sigma} - \tensor{\sigma}_{h,K}) : \tensor{\tau} \dex & = \int_K (2 \mu \tensor{\varepsilon}(\vector{u} - \vector{u}_h) + (p - p_h)\tensor{I}) : \tensor{\tau} \dex & \forall \tensor{\tau} \in \Psi_r(K, \mathbb{M}), \\
        \int_K \div (\proj{\mathcal{S}_h^r} \tensor{\sigma} - \tensor{\sigma}_{h,K}) \cdot \vector{\phi} \dex & = 0 \qquad \forall \vector{\phi} \in \Phi_r(K, \mathbb{R}^2), & \textrm{(only for } K \textrm{ quad.)}.
    \end{align*}
    \end{subequations}
    Setting $\tensor{\omega} = \proj{\mathcal{S}_h^r} \tensor{\sigma} - \tensor{\sigma}_{h,K}$ and applying Lemma \ref{prop:stress_pre} in the previous system, it follows from the definition \eqref{eq:th} and characterization \eqref{eq:charac_mult_HDP} that
    \begin{equation*}
        \| \proj{\mathcal{S}_h^r} \tensor{\sigma} - \tensor{\sigma}_{h,K} \|_{0,K} \leq C (h_K^{1/2} \| \vector{m} - \vector{m}_h \|_{0, \partial K} + \| \tensor{\varepsilon}(\vector{u} - \vector{u}_h) \|_{0,K} + \| p - p_h \|_{0,K}).
    \end{equation*}
    Summing the previous inequality over all elements, using the HDP estimate \eqref{eq:estimate_all}, and considering Remark \ref{rk:norm_M}, we then obtain
    \begin{equation*}
        \| \proj{\mathcal{S}_h^r} \tensor{\sigma} - \tensor{\sigma}_{h} \|_{0,\Omega} \leq C h^{r+1} (\vert \vector{u} \vert_{r+1, \Omega} + \vert p \vert_{r+1,\Omega}).
    \end{equation*}
    Combining the previous estimate with the projection result \eqref{eq:approx_proj_L2}, the $L^2$ estimate \eqref{eq:theo_stress_aux1} follows from the triangular inequality.

    Next, we prove a similar estimate for the divergence of $\tensor{\sigma}_h$.
    Given $\tensor{\tau} \in \mathcal{S}_h^r$, define for each $K \in \TauH$ the vector function
    \[
    \vector{z}(\vector{x}) = \begin{cases}
        \lvert JF_K(\hat{\vector{x}}) \rvert \div \tensor{\tau}(\vector{x}), \; \textrm{ if } \vector{x} \in K, \\
        0, \; \textrm{ otherwise.}
    \end{cases}
    \]
    Since $\hat{\div} (RT_r(\hat{K}, \mathbb{M})) = P_r(\hat{K}, \mathbb{R}^2)$ for $\hat{K}$ the reference triangle, and $\hat{\div} (ABF_r(\hat{K}, \mathbb{M})) = R_r(\hat{K}, \mathbb{R}^2)$ for $\hat{K}$ the reference square, property \eqref{eq:piola_property_div} ensures that $\vector{z}(\vector{x})\rvert_K \in \mathcal{R}^r_K$.
    Therefore, it follows from Lemma \ref{lemma:conservation} that
    \begin{equation*}
        \int_{K} (-\vector{f} - \div \tensor{\sigma}_h) \cdot \vector{z} \dex = \int_{K} \div( \tensor{\sigma} - \tensor{\sigma}_h) \cdot \vector{z} \dex =  0.
    \end{equation*}
    Setting $\tensor{\tau} = \tensor{\sigma}_h - \proj{\mathcal{S}_h^r} \tensor{\sigma}$, we get from the previous equation that
    \begin{equation} \label{eq:aux_H(div)_estimate}
        \begin{split}
            \| \lvert JF_K \rvert^{1/2} \div(\tensor{\sigma} - \tensor{\sigma}_h) \|^2_{0,K} & = \int_K (\lvert JF_K \rvert^{1/2} \div(\tensor{\sigma} - \tensor{\sigma}_h) ) ( \lvert JF_K \rvert^{1/2} \div(\tensor{\sigma} - \tensor{\sigma}_h) ) \dex \\ 
            & + \int_{K} \div( \tensor{\sigma} - \tensor{\sigma}_h) \cdot \vector{z} \dex \\
            & = \int_K (\lvert JF_K \rvert^{1/2} \div(\tensor{\sigma} - \tensor{\sigma}_h) ) ( \lvert JF_K \rvert^{1/2} \div(\tensor{\sigma} - \proj{\mathcal{S}_h^r} \tensor{\sigma}) ) \dex.
        \end{split}
    \end{equation}
    Using the upper bounds \eqref{eq:limJF}, the hypothesis that $\TauH$ is regular, and the Cauchy-Schwarz inequality, eq. \eqref{eq:aux_H(div)_estimate} implies in the existence of a constant $C$, depending only on the mesh shape constant, such that
    \[
    \| \div(\tensor{\sigma} - \tensor{\sigma}_h) \|_{0,K} \leq C \| \div(\tensor{\sigma} - \proj{\mathcal{S}_h^r} \tensor{\sigma}) \|_{0,K}.
    \]
    Summing over all elements and using the projection estimate \eqref{eq:approx_proj_div}, we conclude that
    \begin{equation} \label{eq:theo_stress_aux2}
        \| \div (\tensor{\sigma} - \tensor{\sigma}_h) \|_{0,\Omega} \leq C h^{r+1} \vert \div \tensor{\sigma} \vert_{r+1, \Omega}.
    \end{equation}
    Combining estimates \eqref{eq:theo_stress_aux1} and \eqref{eq:theo_stress_aux2}, we get the desired $H(\div)$ norm.
    Since the constant $C$ in estimate \eqref{eq:estimate_all} is bounded for $\lambda$ going to infinity, the one appearing in \eqref{eq:strees_rate} will be as well. 
\end{proof}

\begin{remark}[RT-based spaces on quadrilaterals] \label{rk:RT_quadrilateral}
    The Raviart-Thomas tensor spaces are also defined on quadrilateral meshes \citep{boffi2013mixed}, and we could use them in the stress post-processing instead of the Arnold-Boffi-Falk based spaces.
    In that case, the definition of the local problems is similar to \eqref{eq:stress_post-processing}, but equation \eqref{eq:stress_post-processing_c} must be removed.
    We remark, however, that the use of Raviart-Thomas based approximations leads to sub-optimal convergence orders in the $H(\div)$ norm on general quadrilateral meshes \citep{arnold2005quadrilateral, correa2022optimal}.
    We illustrate that fact in the numerical experiments of Section \ref{sec:exp_convergence_quad}.
\end{remark}

{\color{ok}

\subsection{Comments on possible extensions} \label{sec:post-extensions}

Although we have presented the stress post-processing in the context of the HDP method, it can be applied to other methods as well.
The core requirement is that the discretization provides an approximation $\vector{t}_h$ for the traction field along the boundaries $\partial K$, and a local approximation $\tilde{\tensor{\sigma}}_h$ for the stress in the element interior. 
In the HDP case, $\vector{t}_h$ comes directly from the hybrid variable $\vector{m}_h$, and we have set $\tilde{\tensor{\sigma}}_h = 2\mu \tensor{\varepsilon}(\vector{u}_h) + p_h \tensor{I}$.

}

{\color{ok}
The weak symmetry result of Lemma \ref{lemma:ht_symmetry} holds whenever $\tilde{\tensor{\sigma}}_h$ is symmetric. Convergence results follow from the convergence properties of $\vector{t}_h$ and $\tilde{\tensor{\sigma}}_h$, applying the same steps of the proof of Theorem \ref{theo:stress}.
In that note, robustness in the incompressibility limit is guaranteed if $\vector{t}_h$ and $\tilde{\tensor{\sigma}}_h$ are locking-free.
Finally, linear momentum conservation can be proven following the same steps of Lemma \ref{lemma:conservation} if
\begin{equation} \label{eq:conservation_condition}
    \int_K \tilde{\tensor{\sigma}}_h : \tensor{\varepsilon}(v) \dex - \int_{\partial K} \vector{t}_h \cdot \vector{v} \des = \int_{K} \vector{f} \cdot \vector{v} \dex, \quad \forall \vector{v} \in \mathcal{Z}_K^r.
\end{equation}
}

{\color{ok}
In this direction, an obvious extension is to the PH method, in which $\vector{t}_h$ is still given by the hybrid variable, and we take $\tilde{\tensor{\sigma}}_h = \tensor{C} \tensor{\varepsilon}(\vector{u}_h)$.
All the properties proven in Sections \ref{sec:post-init} and \ref{sec:post-conv} remain valid, with the exception of the robustness in the incompressibility limit.
It is also straightforward to extend the post-processing to HHO and HDG methods.
Although in these methods the hybrid variables are usually associated with displacement traces, both furnish locally conservative traction fields.
In the HDG formulations \cite{cockburn2013superconvergent, cockburn2018devising, qiu2018hdg}, a numerical traction is computed using the element-wise stress approximation plus a stabilization term, and is locally conservative by design. 
Therefore, condition \eqref{eq:conservation_condition} can be met by taking $\tilde{\tensor{\sigma}}_h$ as the element-wise stress approximation and $\vector{t}_h$ as the numerical traction.
For the HHO method, \cite{di2015equilibrated, cicuttin2021hybrid} show how to recover a conservative traction field from local procedures that satisfies \eqref{eq:conservation_condition} for a suitable choice of $\tilde{\tensor{\sigma}}$.

Finally, we remark that the core ideas of the post-processing strategies generalize to three-dimensional domains as well.
The precise form of the local systems to be solved depends on the element topology and on the $H(\div)$ space employed \cite{correa2022optimal}.
For tetrahedral meshes, for instance, the classical Raviart-Thomas-Nédeléc spaces yield local systems analogous to \eqref{eq:stress_post-processing_a} -- \eqref{eq:stress_post-processing_b}, with the obvious modifications in the test spaces.
Although not detailed here, all the results proven for triangles, regarding convergence, weak symmetry, and linear momentum conservation, naturally extend to tetrahedral elements.
}

\section{Numerical tests} \label{sec:exp}

This section is dedicated to validating the results of Sections \ref{sec:spaces} and \ref{sec:post} through some simple numerical experiments.
Our experiments are divided into {\color{ok}three} blocks.
In the first one, we verify the convergence orders of the HDP method with the post-processing strategy \eqref{eq:stress_post-processing} proven in Corollary \ref{cor:convergence_rates} and Theorem \ref{theo:stress}.
We also illustrate the importance of using unmapped pressure approximation spaces and Arnold-Boffi-Falk-based tensor approximation spaces to achieve optimal convergence orders on non-affine quadrilateral meshes.

The second block of experiments aims to illustrate the robustness of the HDP method and the stress post-processing strategy with respect to locking.
In this block, we also compare the performance of the HDP, PH, and AP methods in the incompressibility limit, showing that the HDP is the only one robust on nearly incompressible problems for all meshes and spaces considered.
{\color{ok}
Finally, the third block, also related to near-incompressibility behavior, showcases the locking-free property of the HDP method in Cook's membrane benchmark.
}

\subsection{Convergence test} \label{sec:exp_convergence}

In this first block of experiments, we solve the model problem \eqref{eq:model} on the square domain {\color{ok}$\Omega = (-1,1) \times (-1,1)$}, setting Dirichlet boundary conditions on the top and bottom, and Neumann boundary conditions on the remaining of $\partial \Omega$.
We assume an isotropic material with Lamé coefficients given by $\mu = 1$ and $\lambda = 0.3$. 
The values for $\vector{f}$, $\vector{u}_D$, and $\vector{t}_N$ were chosen so the exact solution for the displacement is
\[
\vector{u}(\vector{x}) = \begin{bmatrix}
    \sin (\pi x_1) \sin (\pi x_2) \\
    \sin (\pi x_1) \sin (\pi x_2)
\end{bmatrix}.
\]

\begin{figure}[hbt]
\centering
\subfloat{
\includegraphics[width=.235\linewidth]{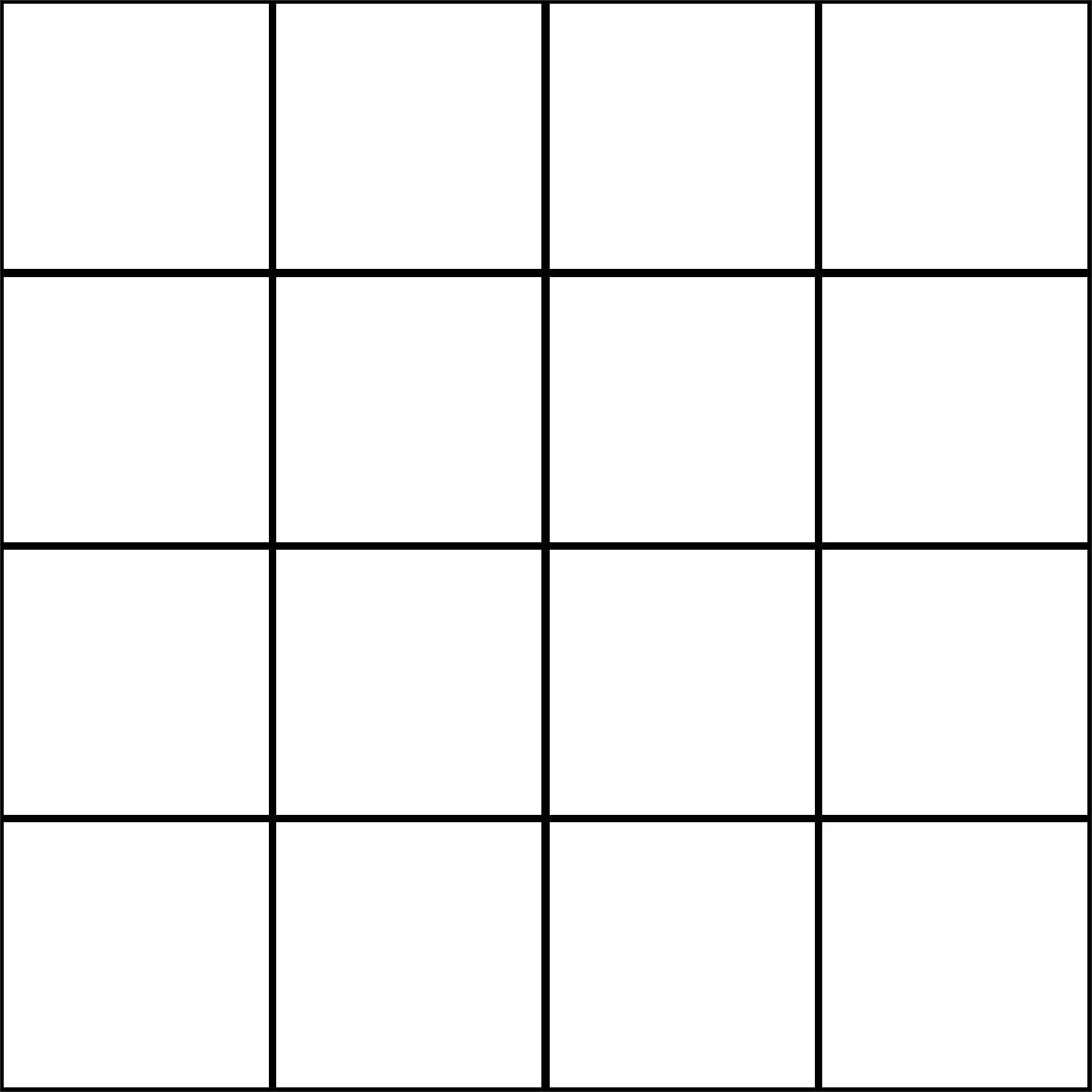}
}
\hspace{2.5pt}
%
\subfloat{
\includegraphics[width=.235\linewidth]{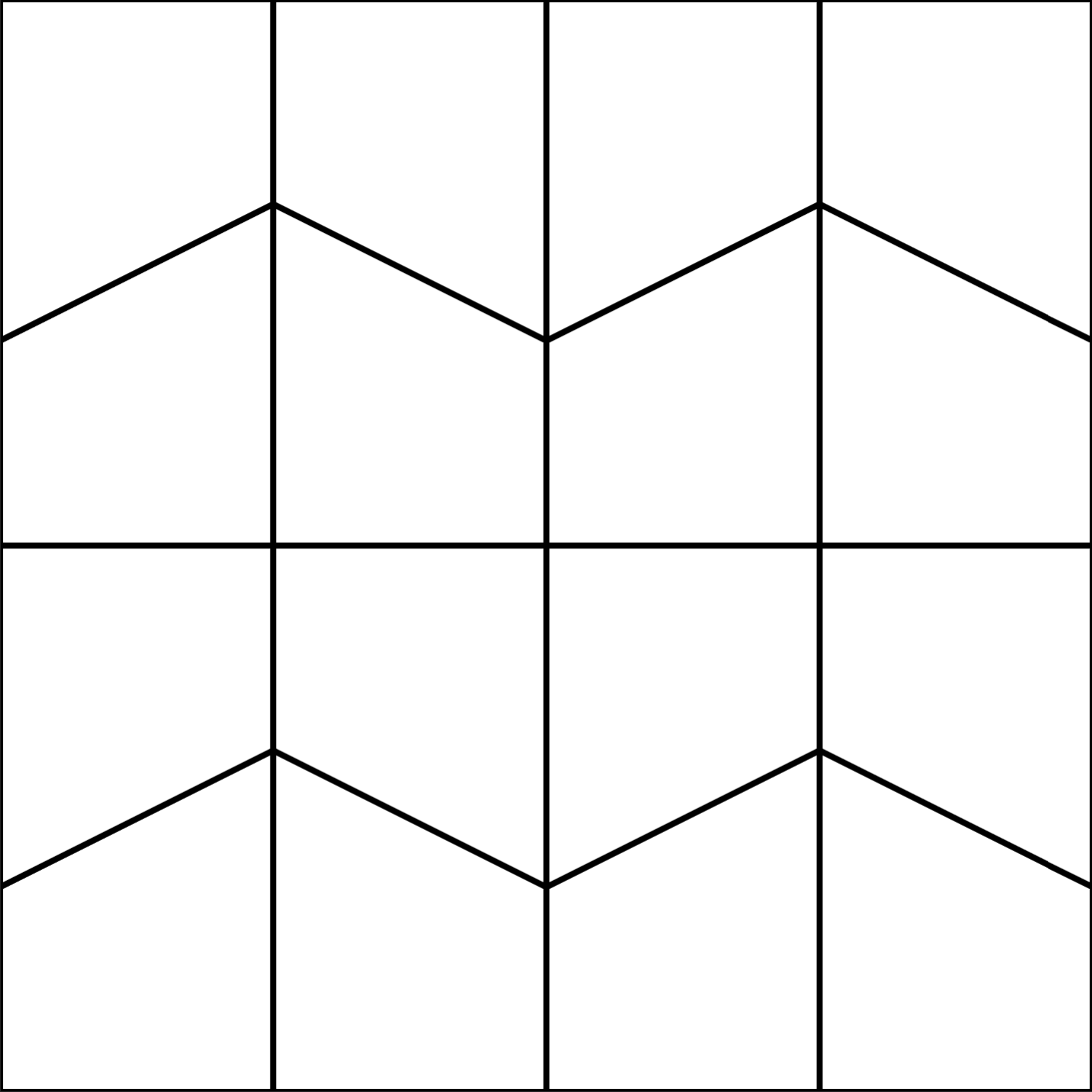}
}
\hspace{2.5pt}
%
\subfloat{
\includegraphics[width=.235\linewidth]{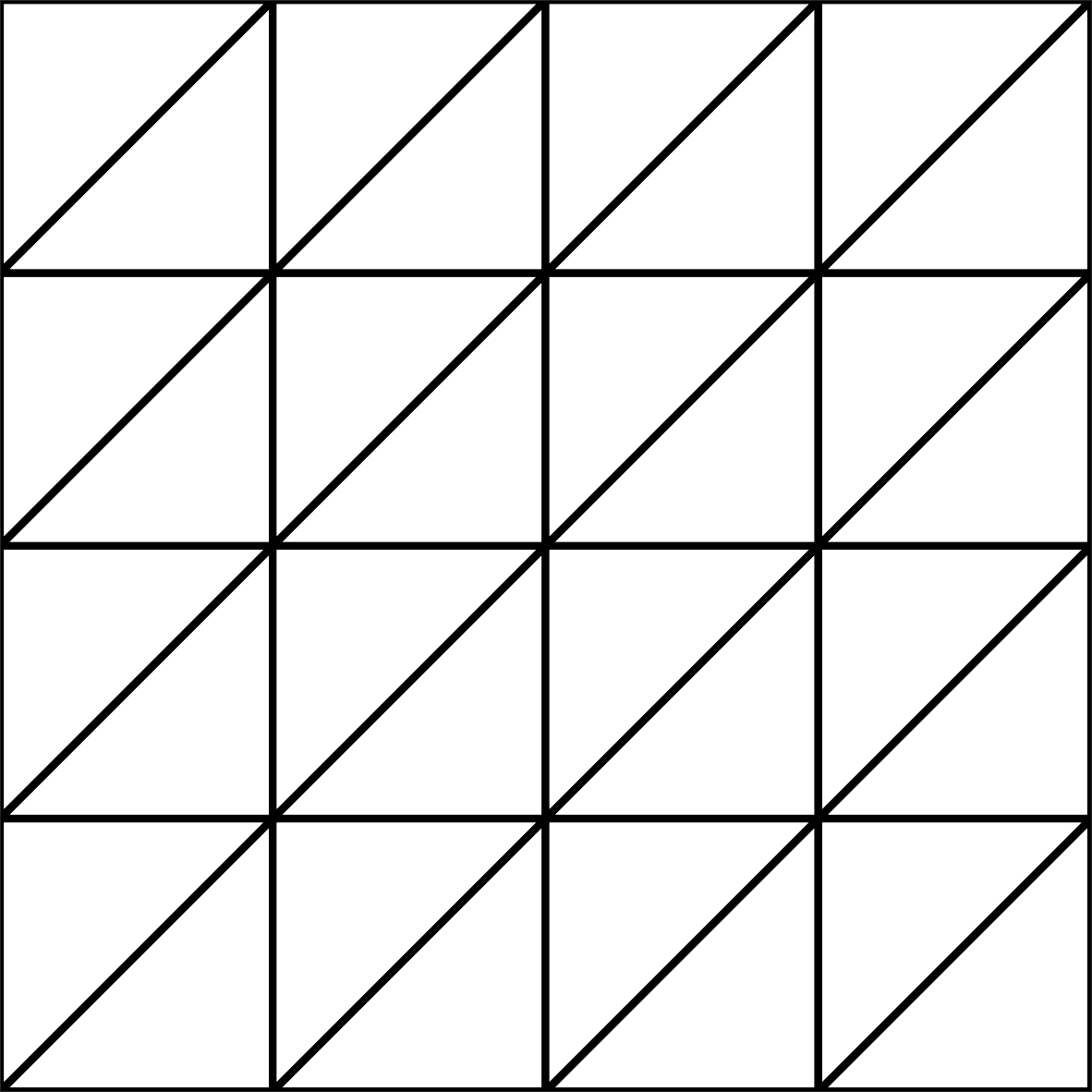}
}
\\
\vspace{0.5em}
\subfloat{
\includegraphics[width=.235\linewidth]{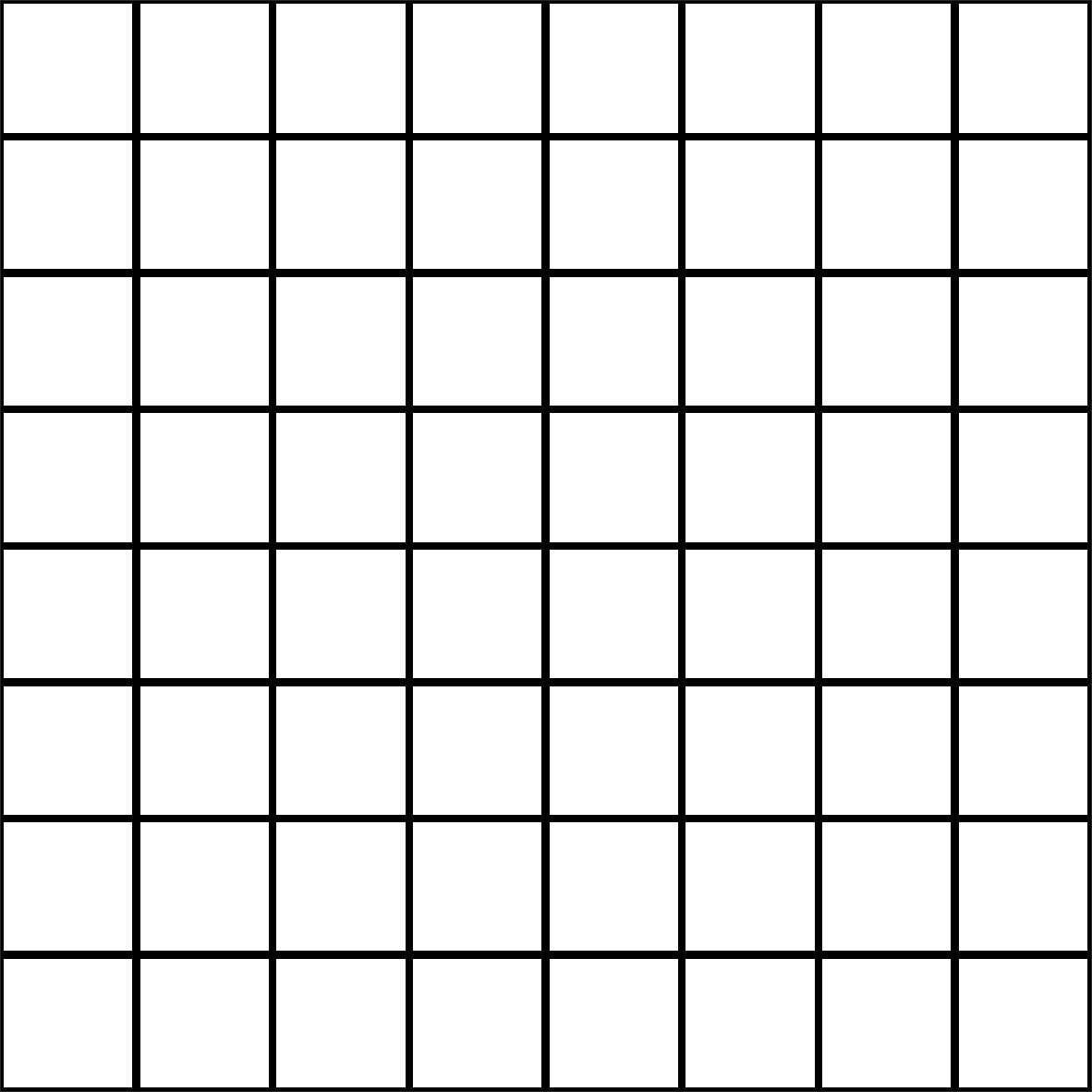}
}
\hspace{2.5pt}
\subfloat{
\includegraphics[width=.235\linewidth]{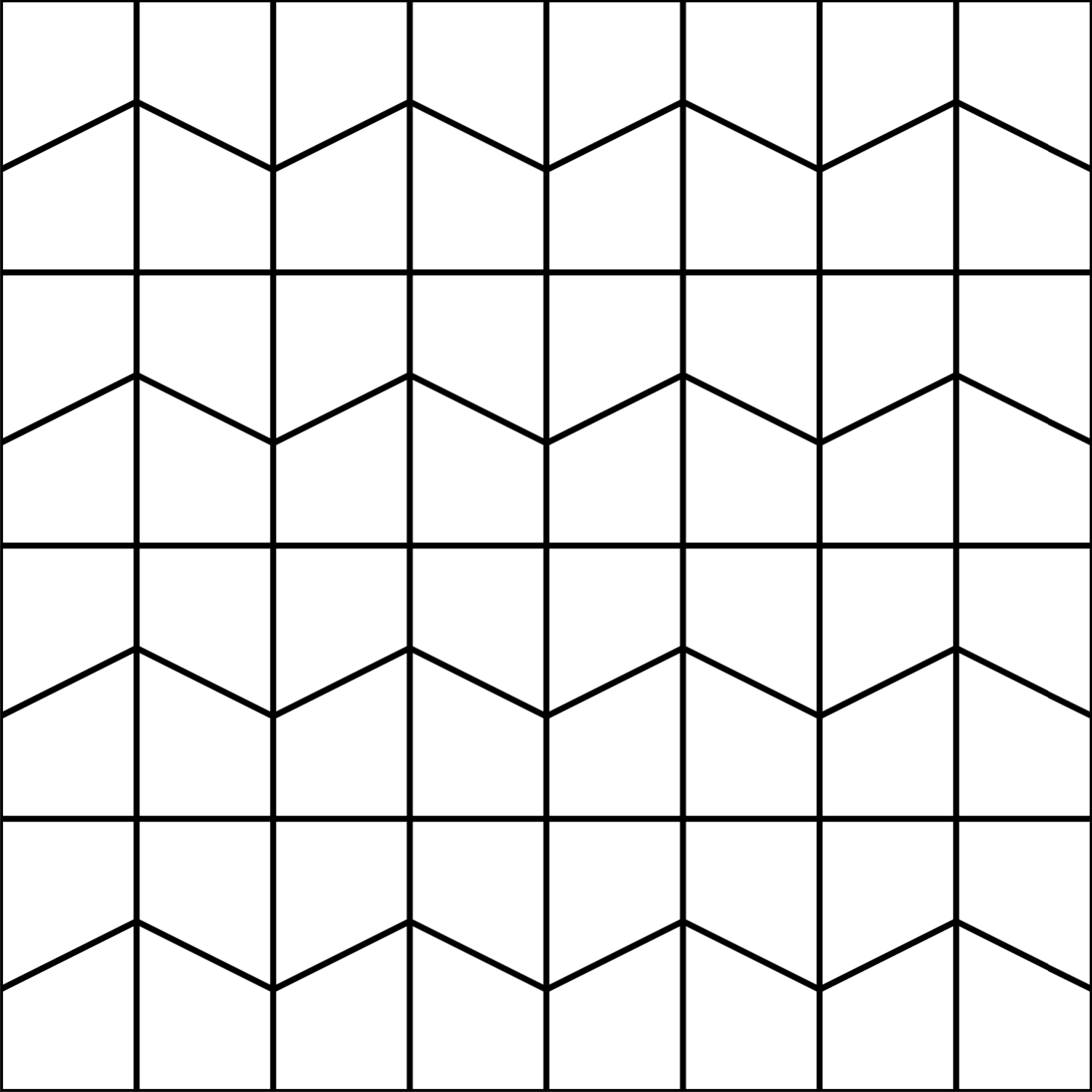}
}
\hspace{2.5pt}
\subfloat{
\includegraphics[width=.235\linewidth]{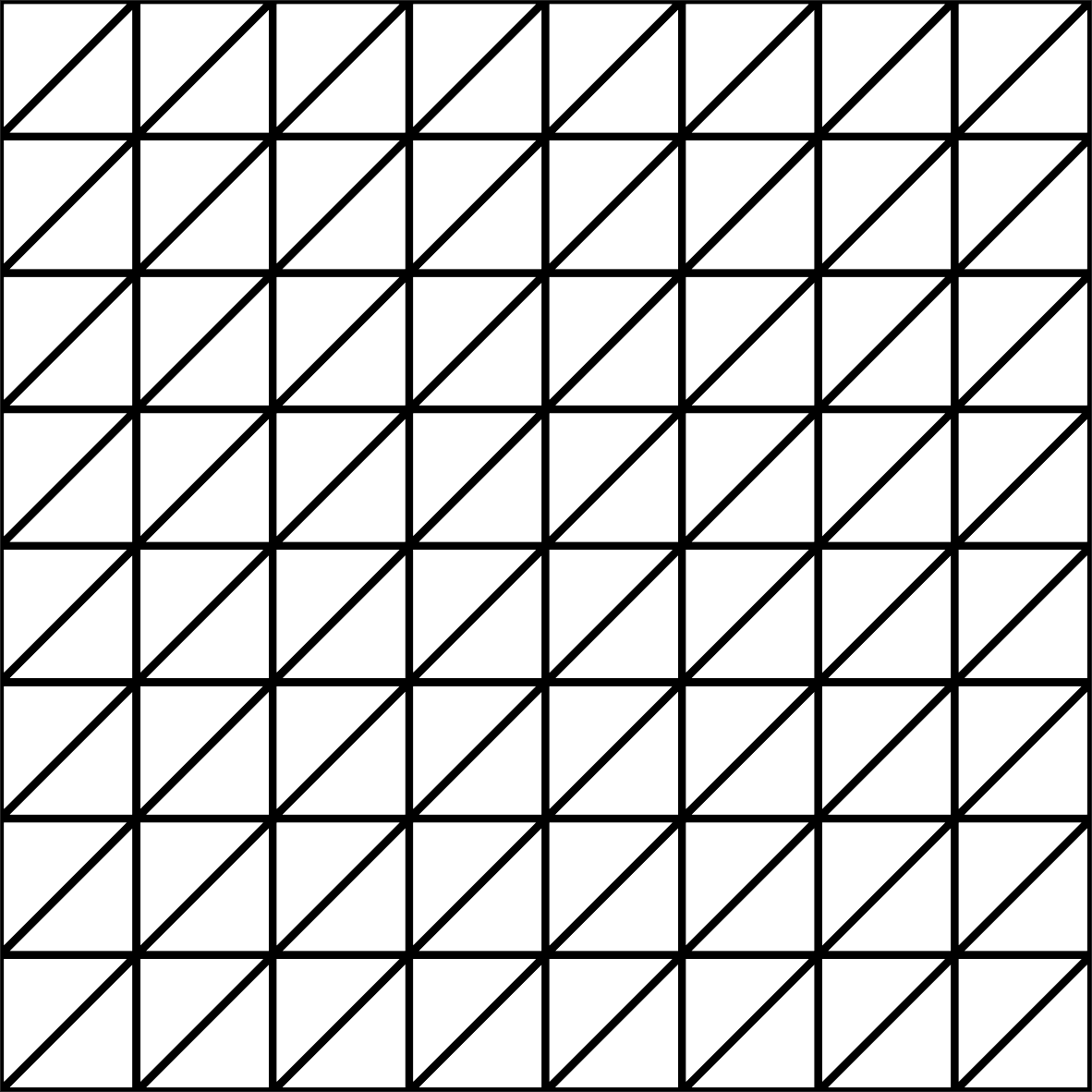}
}
\caption{Square, trapezoidal, and triangular meshes for $n=4$ and $n=8$.}
\label{fig:meshes}
\end{figure}

We consider three mesh types for domain discretization: square, triangular, and trapezoidal.
For the square meshes, the domain is partitioned into $n \times n$ squares, with $n$ being powers of two 
(first column of Figure \ref{fig:meshes}). 
The triangular meshes are obtained by dividing each element of the square meshes through the segment connecting the southwest to the northeast vertices, resulting in $2n^2$ triangular elements (third column of Figure \ref{fig:meshes}).
For the trapezoidal meshes, the $x_2$ coordinate of the vertices of the square meshes is perturbed, leading to $n^2$ congruent trapezoids (second column of Figure \ref{fig:meshes}).
Notice that the square and triangular meshes are generated by affine transformations of the reference element, while the trapezoidal mesh is obtained through non-affine bilinear transformations.

To obtain the approximated solutions, we use the HDP method with the spaces $\mathcal{X}_h^{r+1} \times \mathcal{M}_h^{r} \times \mathcal{P}_h^{r}$, {\color{ok} for $r = 1,2$ and $3$}, as defined in Section \ref{sec:spaces_HDP}.
{\color{ok}
For the implementation of the HDP method, we use the static condensation described in Appendix \ref{ap:condensation}, which allows the global system to be expressed primarily in terms of $\vector{m}_h$.
%
}
Approximations for the stress field are obtained using the post-processing strategy described in Section \ref{sec:post}.
We use the $L^2$ norm to measure the errors $\vector{u} - \vector{u}_h$ and $p - p_h$, and the $H(\div)$ and $\vertiii{ \cdot }_{\mathcal{M}_h}$ norms for the errors $\tensor{\sigma} - \tensor{\sigma}_h$ and $\vector{m} - \vector{m}_h$, respectively.
The results are presented in Tables \ref{tab:conv_trig} to \ref{tab:conv_trap}.
For the three types of mesh used, the numerically measured convergence orders agree with the theoretical \textit{a priori} estimates of Corollary \ref{cor:convergence_rates} and Theorem \ref{theo:stress} (See also remarks \ref{rk:norm_M} and \ref{rk:l2_rate}): $\mathcal{O}(h^{r+2})$ $L^2$ convergence for $\vector{u}_h$ and $\mathcal{O}(h^{r+1})$ convergence for the remaining variables.
\begin{table}[htb]
\centering
\caption{\color{ok}
First test; HDP convergence on triangular meshes.
}
\resizebox{0.8\textwidth}{!}{
\begin{tabular}{c |l  c l  c l  c l c}
\hline
& \multicolumn{2}{c}{$ \|\vector{u}-\vector{u}_h \|_{0} $} & \multicolumn{2}{c}{$  \vertiii{\vector{m} - \vector{m}_h}_{\mathcal{M}_h} $} & \multicolumn{2}{c}{$  \| p - p_h \|_{0} $} & \multicolumn{2}{c}{$  \| \tensor{\sigma} - \tensor{\sigma}_h \|_{H(\div)} $} \\
$n$ & err. & order& err. & order& err. & order& err. & rate\\
\hline
\multicolumn{9}{c}{Space $\mathcal{X}_h^{2} \times \mathcal{M}_h^{1} \times \mathcal{P}_h^{1}$} \Tstrut \Bstrut \\
\hline
8   & 9.84e-03 &  3.0 & 5.56e-01 &  1.9 & 7.00e-02 &  1.9 & 2.42e+00 &  1.9 \\
16  & 1.37e-03 &  2.8 & 1.41e-01 &  2.0 & 1.78e-02 &  2.0 & 6.16e-01 &  2.0 \\
32  & 1.88e-04 &  2.9 & 3.41e-02 &  2.0 & 4.45e-03 &  2.0 & 1.55e-01 &  2.0 \\
64  & 2.46e-05 &  2.9 & 8.32e-03 &  2.0 & 1.11e-03 &  2.0 & 3.87e-02 &  2.0 \\
128 & 3.11e-06 &  3.0 & 2.06e-03 &  2.0 & 2.79e-04 &  2.0 & 9.68e-03 &  2.0 \\
\hline
\multicolumn{9}{c}{Space $\mathcal{X}_h^{3} \times \mathcal{M}_h^{2} \times \mathcal{P}_h^{2}$} \Tstrut \Bstrut \\
\hline
8   & 1.06e-03 &  3.9 & 9.42e-02 &  2.7 & 8.59e-03 &  2.9 & 2.77e-01 &  2.9 \\
16  & 6.41e-05 &  4.0 & 1.27e-02 &  2.9 & 1.10e-03 &  3.0 & 3.51e-02 &  3.0 \\
32  & 3.93e-06 &  4.0 & 1.59e-03 &  3.0 & 1.39e-04 &  3.0 & 4.41e-03 &  3.0 \\
64  & 2.43e-07 &  4.0 & 1.98e-04 &  3.0 & 1.73e-05 &  3.0 & 5.51e-04 &  3.0 \\
128 & 1.52e-08 &  4.0 & 2.46e-05 &  3.0 & 2.17e-06 &  3.0 & 6.89e-05 &  3.0 \\
\hline
\multicolumn{9}{c}{Space $\mathcal{X}_h^{4} \times \mathcal{M}_h^{3} \times \mathcal{P}_h^{3}$} \Tstrut \Bstrut \\
\hline
8   & 5.58e-05 &  4.9 & 6.43e-03 &  4.0 & 7.26e-04 &  3.9 & 2.44e-02 &  3.9 \\
16  & 1.79e-06 &  5.0 & 3.66e-04 &  4.1 & 4.62e-05 &  4.0 & 1.54e-03 &  4.0 \\
32  & 5.64e-08 &  5.0 & 2.11e-05 &  4.1 & 2.90e-06 &  4.0 & 9.68e-05 &  4.0 \\
64  & 1.77e-09 &  5.0 & 1.27e-06 &  4.0 & 1.82e-07 &  4.0 & 6.06e-06 &  4.0 \\
128 & 5.53e-11 &  5.0 & 7.93e-08 &  4.0 & 1.13e-08 &  4.0 & 3.79e-07 &  4.0 \\
\hline
\end{tabular}
}
\label{tab:conv_trig}
\end{table}
\begin{table}[htb]
\centering
\caption{\color{ok}
First test; HDP convergence on square meshes. 
}
\resizebox{0.8\textwidth}{!}{
\begin{tabular}{c|l  c l  c l  c l c}
\hline
& \multicolumn{2}{c}{$ \|\vector{u}-\vector{u}_h \|_{0} $} & \multicolumn{2}{c}{$  \vertiii{\vector{m} - \vector{m}_h}_{\mathcal{M}_h} $} & \multicolumn{2}{c}{$  \| p - p_h \|_{0} $} & \multicolumn{2}{c}{$  \| \tensor{\sigma} - \tensor{\sigma}_h \|_{H(\div)} $} \\
$n$ & err. & order& err. & order& err. & order& err. & rate\\
\hline
\multicolumn{9}{c}{Space $\mathcal{X}_h^{2} \times \mathcal{M}_h^{1} \times \mathcal{P}_h^{1}$} \Tstrut \Bstrut \\
\hline
8   & 7.03e-03 &  3.0 & 4.82e-01 &  1.9 & 7.92e-02 &  1.9 & 3.05e-01 &  2.3 \\
16  & 9.19e-04 &  2.9 & 1.19e-01 &  2.0 & 2.01e-02 &  2.0 & 7.15e-02 &  2.1 \\
32  & 1.18e-04 &  3.0 & 2.93e-02 &  2.0 & 5.06e-03 &  2.0 & 1.75e-02 &  2.0 \\
64  & 1.50e-05 &  3.0 & 7.30e-03 &  2.0 & 1.27e-03 &  2.0 & 4.35e-03 &  2.0 \\
128 & 1.88e-06 &  3.0 & 1.82e-03 &  2.0 & 3.17e-04 &  2.0 & 1.09e-03 &  2.0 \\
\hline
\multicolumn{9}{c}{Space $\mathcal{X}_h^{3} \times \mathcal{M}_h^{2} \times \mathcal{P}_h^{2}$} \Tstrut \Bstrut \\
\hline
8   & 3.15e-04 &  4.0 & 4.43e-02 &  2.8 & 1.00e-02 &  2.9 & 2.26e-02 &  3.1 \\
16  & 1.96e-05 &  4.0 & 5.80e-03 &  2.9 & 1.27e-03 &  3.0 & 2.77e-03 &  3.0 \\
32  & 1.22e-06 &  4.0 & 7.36e-04 &  3.0 & 1.60e-04 &  3.0 & 3.44e-04 &  3.0 \\
64  & 7.59e-08 &  4.0 & 9.26e-05 &  3.0 & 2.00e-05 &  3.0 & 4.29e-05 &  3.0 \\
128 & 4.74e-09 &  4.0 & 1.16e-05 &  3.0 & 2.50e-06 &  3.0 & 5.37e-06 &  3.0 \\
\hline
\multicolumn{9}{c}{Space $\mathcal{X}_h^{4} \times \mathcal{M}_h^{3} \times \mathcal{P}_h^{3}$} \Tstrut \Bstrut \\
\hline
8   & 1.39e-05 &  4.9 & 2.79e-03 &  4.0 & 9.51e-04 &  3.9 & 1.32e-03 &  4.0 \\
16  & 4.48e-07 &  5.0 & 1.54e-04 &  4.2 & 6.04e-05 &  4.0 & 8.24e-05 &  4.0 \\
32  & 1.42e-08 &  5.0 & 8.94e-06 &  4.1 & 3.79e-06 &  4.0 & 5.14e-06 &  4.0 \\
64  & 4.44e-10 &  5.0 & 5.41e-07 &  4.0 & 2.37e-07 &  4.0 & 3.21e-07 &  4.0 \\
128 & 1.39e-11 &  5.0 & 3.34e-08 &  4.0 & 1.48e-08 &  4.0 & 2.00e-08 &  4.0 \\
\hline
\end{tabular}
}
\label{tab:conv_quad}
\end{table}
\begin{table}[htb]
\centering
\caption{
\color{ok}
First test; HDP convergence on trapezoidal meshes. 
}
\resizebox{0.8\textwidth}{!}{
\begin{tabular}{c |l  c l  c l  c l c}
\hline
& \multicolumn{2}{c}{$ \|\vector{u}-\vector{u}_h \|_{0} $} & \multicolumn{2}{c}{$  \vertiii{\vector{m} - \vector{m}_h}_{\mathcal{M}_h} $} & \multicolumn{2}{c}{$  \| p - p_h \|_{0} $} & \multicolumn{2}{c}{$  \| \tensor{\sigma} - \tensor{\sigma}_h \|_{H(\div)} $} \\
$n$ & err. & order& err. & order& err. & order& err. & rate \\
\hline
\multicolumn{9}{c}{Space $\mathcal{X}_h^{2} \times \mathcal{M}_h^{1} \times \mathcal{P}_h^{1}$} \Tstrut \Bstrut \\
\hline
8   & 9.53e-03 &  3.0 & 6.73e-01 &  2.0 & 8.23e-02 &  1.9 & 4.16e-01 &  2.3 \\
16  & 1.22e-03 &  3.0 & 1.73e-01 &  2.0 & 2.10e-02 &  2.0 & 9.77e-02 &  2.1 \\
32  & 1.56e-04 &  3.0 & 4.35e-02 &  2.0 & 5.27e-03 &  2.0 & 2.40e-02 &  2.0 \\
64  & 1.97e-05 &  3.0 & 1.09e-02 &  2.0 & 1.32e-03 &  2.0 & 5.98e-03 &  2.0 \\
128 & 2.47e-06 &  3.0 & 2.73e-03 &  2.0 & 3.30e-04 &  2.0 & 1.49e-03 &  2.0 \\
\hline
\multicolumn{9}{c}{Space $\mathcal{X}_h^{3} \times \mathcal{M}_h^{2} \times \mathcal{P}_h^{2}$} \Tstrut \Bstrut \\
\hline
8   & 4.94e-04 &  3.9 & 7.19e-02 &  2.8 & 1.09e-02 &  2.9 & 3.35e-02 &  3.1 \\
16  & 3.16e-05 &  4.0 & 9.12e-03 &  3.0 & 1.39e-03 &  3.0 & 4.08e-03 &  3.0 \\
32  & 1.99e-06 &  4.0 & 1.12e-03 &  3.0 & 1.74e-04 &  3.0 & 5.03e-04 &  3.0 \\
64  & 1.25e-07 &  4.0 & 1.38e-04 &  3.0 & 2.18e-05 &  3.0 & 6.25e-05 &  3.0 \\
128 & 7.80e-09 &  4.0 & 1.71e-05 &  3.0 & 2.73e-06 &  3.0 & 7.80e-06 &  3.0 \\
\hline
\multicolumn{9}{c}{Space $\mathcal{X}_h^{4} \times \mathcal{M}_h^{3} \times \mathcal{P}_h^{3}$} \Tstrut \Bstrut \\
\hline
8   & 2.22e-05 &  4.9 & 4.98e-03 &  4.0 & 1.11e-03 &  3.9 & 1.89e-03 &  4.0 \\
16  & 7.10e-07 &  5.0 & 3.03e-04 &  4.0 & 7.09e-05 &  4.0 & 1.18e-04 &  4.0 \\
32  & 2.24e-08 &  5.0 & 1.87e-05 &  4.0 & 4.45e-06 &  4.0 & 7.39e-06 &  4.0 \\
64  & 7.01e-10 &  5.0 & 1.16e-06 &  4.0 & 2.79e-07 &  4.0 & 4.62e-07 &  4.0 \\
128 & 2.19e-11 &  5.0 & 7.22e-08 &  4.0 & 1.74e-08 &  4.0 & 2.88e-08 &  4.0 \\
\hline
\end{tabular}
}
\label{tab:conv_trap}
\end{table}

\subsubsection{Comments on approximations in  quadrilaterals} \label{sec:exp_convergence_quad}

In view of Remark \ref{rk:RT_quadrilateral}, we now perform the post-processing \eqref{eq:stress_post-processing} using the quadrilateral Raviart-Thomas based tensor space of index $1$ (see Section 4.1.2 from \cite{taraschi2024global} for its definition).
For that, the HDP solution was previously computed using the $\mathcal{X}_h^{2} \times \mathcal{M}_h^{1} \times \mathcal{P}_h^{1}$ space.
As shown in Table \ref{tab:conv_stress_rt}, RT-based tensor spaces guarantee optimal convergence on square meshes but furnish sub-optimal orders for trapezoidal partitions.
This phenomenon is discussed in depth in \cite{arnold2005quadrilateral, correa2022optimal} and justifies using ABF-based spaces on general non-affine convex quadrilateral meshes.
\begin{table}[htb]
\centering
\caption{
\color{ok}
First test; stress convergence for the RT-based post-processing on square and trapezoidal meshes.
}
\begin{tabular}{c|l  c| l  c}
\hline
\multicolumn{3}{c|}{Square meshes} & \multicolumn{2}{c}{Trapezoidal meshes} \\
\hline
& \multicolumn{2}{c|}{$ \| \tensor{\sigma} - \tensor{\sigma}_h \|_{H(\div)} $} & \multicolumn{2}{c}{$ \| \tensor{\sigma} - \tensor{\sigma}_h \|_{H(\div)} $}  \\
$n$ & err. & order& err. & order\\
\hline
8   & 1.62e+00 &  2.0 & 2.41e+00 &  1.7 \\
16  & 4.08e-01 &  2.0 & 9.00e-01 &  1.4 \\
32  & 1.02e-01 &  2.0 & 4.01e-01 &  1.2 \\
64  & 2.55e-02 &  2.0 & 1.94e-01 &  1.1 \\
128 & 6.38e-03 &  2.0 & 9.62e-02 &  1.0 \\
\hline
\end{tabular}
\label{tab:conv_stress_rt}
\end{table}

Next, to showcase the comments of Remark \ref{rk:mapped_vs_geo}, we solve our test problem using the space $\mathcal{X}_h^{2} \times \mathcal{M}_h^{1} \times \bar{\mathcal{P}}_h^{1}$, where $\bar{\mathcal{P}}_h^{1}$ is constructed by mapping the functions in $P_1(\hat{K}, \mathbb{R}^2)$.
The results of Table \ref{tab:conv_mapped} show that the use of mapped pressure spaces leads to sub-optimal convergence on 
trapezoidal meshes.
We underline that the loss of optimality is verified for all variables, not only the pressure field.
Such degradation is expected whenever mapped pressure approximation spaces are applied on non-affine quadrilateral meshes \citep{boffi2002quadrilateral}.
\begin{table}[htb]
\centering
\caption{
\color{ok}
First test; HDP convergence on square and trapezoidal meshes using mapped pressure spaces.
}
\resizebox{\textwidth}{!}{
\begin{tabular}{c|l c l  c l  c| l  c l  c l  c}
\hline
\multicolumn{7}{c|}{Square meshes} & \multicolumn{6}{c}{Trapezoidal meshes} \\
\hline
& \multicolumn{2}{c}{$ \|\vector{u}-\vector{u}_h \|_{0} $} & \multicolumn{2}{c}{$  \vertiii{\vector{m} - \vector{m}_h}_{\mathcal{M}_h} $} & \multicolumn{2}{c|}{$  \| p - p_h \|_{0} $} & \multicolumn{2}{c}{$ \|\vector{u}-\vector{u}_h \|_{0} $} & \multicolumn{2}{c}{$  \vertiii{\vector{m} - \vector{m}_h}_{\mathcal{M}_h} $} & \multicolumn{2}{c}{$  \| p - p_h \|_{0} $} \\
$n$ & err. & order& err. & order& err. & order& err. & order& err. & order& err. & order\\
\hline
8   & 7.03e-03 &  3.0 & 4.82e-01 &  1.9 & 7.92e-02 &  1.9 & 9.71e-03 &  3.0 & 6.75e-01 &  2.0 & 9.31e-02 &  1.8 \\
16  & 9.19e-04 &  2.9 & 1.19e-01 &  2.0 & 2.01e-02 &  2.0 & 1.27e-03 &  2.9 & 1.77e-01 &  1.9 & 3.02e-02 &  1.6 \\
32  & 1.18e-04 &  3.0 & 2.93e-02 &  2.0 & 5.06e-03 &  2.0 & 1.73e-04 &  2.9 & 4.89e-02 &  1.9 & 1.20e-02 &  1.3 \\
64  & 1.50e-05 &  3.0 & 7.30e-03 &  2.0 & 1.27e-03 &  2.0 & 2.68e-05 &  2.7 & 1.59e-02 &  1.6 & 5.56e-03 &  1.1 \\
128 & 1.88e-06 &  3.0 & 1.82e-03 &  2.0 & 3.17e-04 &  2.0 & 5.16e-06 &  2.4 & 6.41e-03 &  1.3 & 2.72e-03 &  1.0 \\
\hline
\end{tabular}
}
\label{tab:conv_mapped}
\end{table}

\clearpage 
\subsection{Robustness towards locking} \label{sec:exp_locking}

In this second block of experiments, we aim to verify the behavior of the HDP method on nearly incompressible problems.
We adopt the same test problem of \cite{brenner1993nonconforming}, obtained from \eqref{eq:model}
by setting $\Omega = (-1,1) \times (-1,1)$, $\Gamma_D = \partial \Omega$, $\vector{u}_D = 0$, and  the source term 
%
\[
\vector{f}(\vector{x}) = - \begin{bmatrix}
    4 \pi^2 \sin(2 \pi x_2)(2\cos(2\pi x_1) - 1) - \cos(\pi(x_1+x_2)) + \left(\frac{2}{1 + \lambda}\right) \sin(\pi x_1) \sin(\pi x_2) \\
    4 \pi^2 \sin(2 \pi x_1)(1 - 2\cos(2\pi x_2)) - \cos(\pi(x_1+x_2)) + \left(\frac{2}{1 + \lambda}\right) \sin(\pi x_1) \sin(\pi x_2) \\
\end{bmatrix}.
\]
Furthermore, the Lamé constants are given by %
\[
\mu=1 \quad \mbox{ and } \quad \lambda = \frac{\nu}{1 - 2\nu},
\]
%
%
and the exact solution for the resulting problem is
\[
\vector{u}(\vector{x}) = \begin{bmatrix}
    \sin(2\pi x_2)(\cos(2\pi x_1) -1) + \left(\frac{1}{1+\lambda}\right) \sin(\pi x_1) \sin(\pi x_2) \\
    \sin(2\pi x_1)(1 -\cos(2\pi x_2)) + \left(\frac{1}{1+\lambda}\right) \sin(\pi x_1) \sin(\pi x_2)
\end{bmatrix},
\]
which remains bounded as $\lambda$ goes to infinity.
The expression for $\lambda$ derives from \eqref{eq:lame} when $\mu = 1$ is kept fixed.

The domain is then partitioned into triangular and trapezoidal meshes with $n=64$, as described in Section \ref{sec:exp_convergence}.
Numerical solutions are obtained using the HDP method with spaces $\mathcal{X}_h^{r+1} \times \mathcal{M}_h^{r} \times \mathcal{P}_h^{r}$,  for $r=1,2$, and compared with the PH and AP approximations obtained with the spaces $\mathcal{X}_h^{\ph,r+1} \times \mathcal{M}_h^{r}$.
Here, $\mathcal{X}_h^{\ph,r+1}$ denotes the space constructed through \eqref{eq:constr_Xh} by setting $\hat{\mathcal{X}} = P_{r+1}^{\ph}(\hat{K}, \mathbb{R}^2)$ on triangular meshes and $\hat{\mathcal{X}} = Q_{r+1}^{\ph}(\hat{K}, \mathbb{R}^2)$ on quadrilateral partitions, as defined in Section \ref{sec:inf-sup_PH}.
Again, the HDP method was implemented using the  static condensation procedure described in {\color{ok}Appendix} \ref{ap:condensation}.

We apply the local post-processing described in Section \ref{sec:post} to recover stress approximations for the HDP solution.
{\color{ok}Considering Section \ref{sec:post-extensions}}, such a strategy is also available to approximate the stress field from the PH solution.
On the other hand, this stress recovery is not available (at least in the current form) for the AP method since the characterization for its Lagrange multipliers differs (see eq. \eqref{eq:charac_m_AP}).

\begin{figure}[htb]
    \centering
    \begin{tabular}{c c c}
        HDP, $\mathcal{X}_h^2 \times \mathcal{M}_h^1 \times \mathcal{P}_h^1$ & PH, $\mathcal{X}_h^{\ph,2} \times \mathcal{M}_h^1$ & AP, $\mathcal{X}_h^{\ph,2} \times \mathcal{M}_h^1$  \\
        \includegraphics[width=0.3\linewidth]{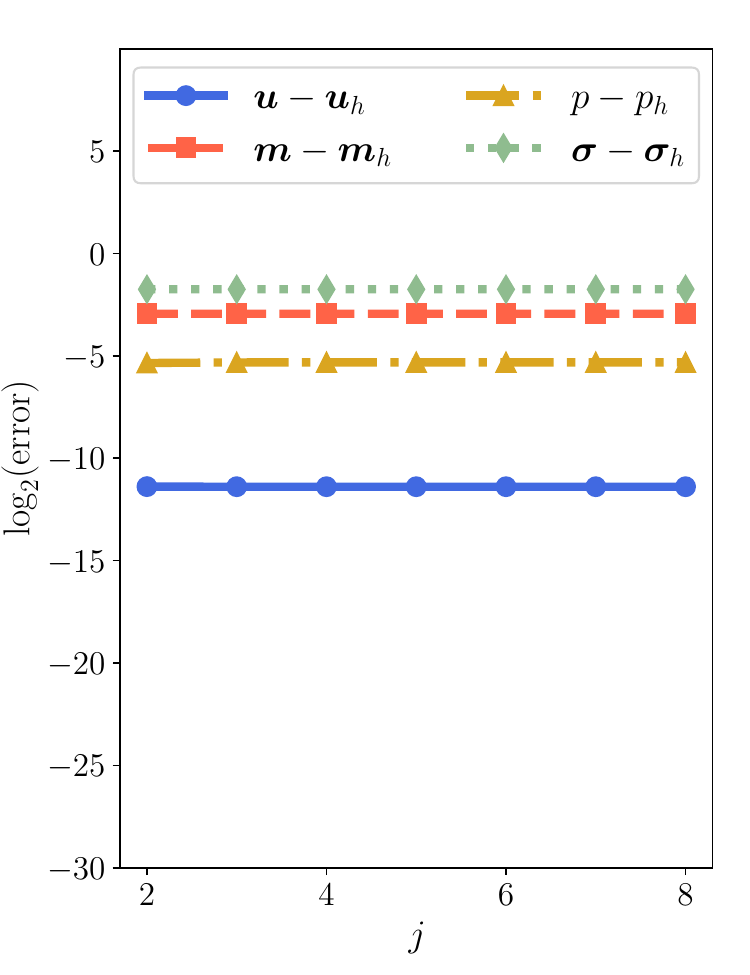} & \includegraphics[width=0.3\linewidth]{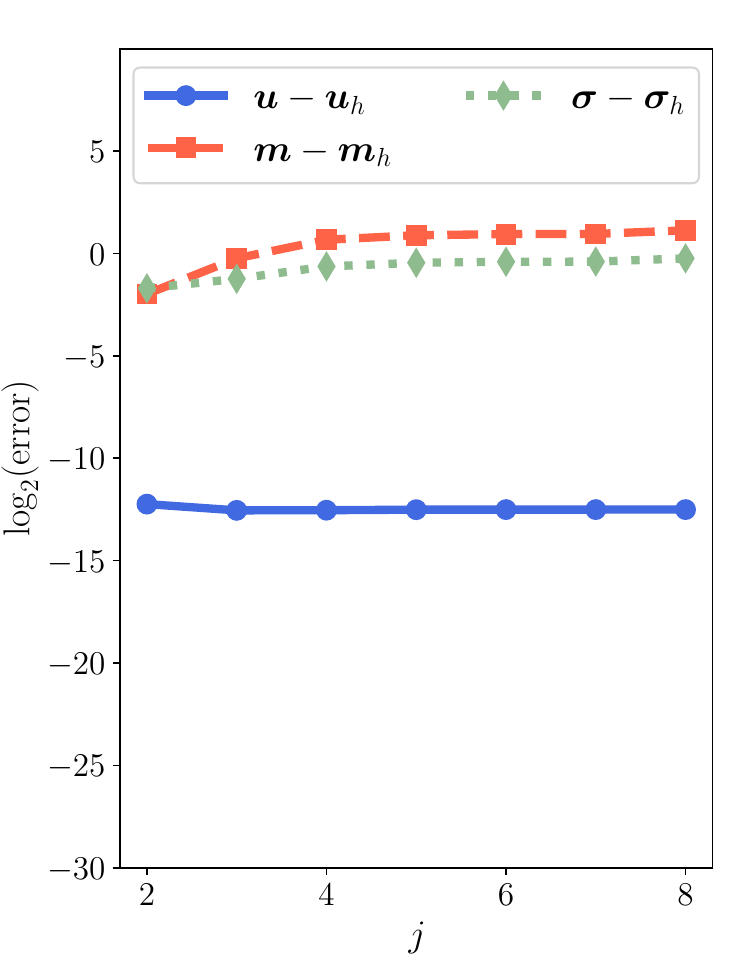} & \includegraphics[width=0.3\linewidth]{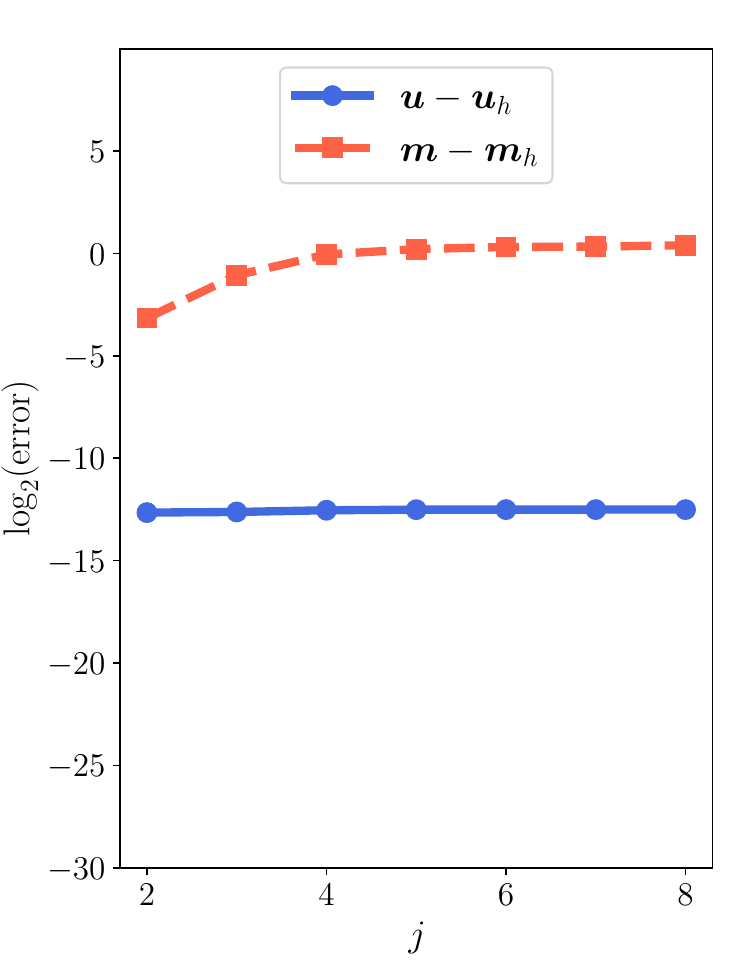}
        \\
        HDP, $\mathcal{X}_h^3 \times \mathcal{M}_h^2 \times \mathcal{P}_h^2$ & PH, $\mathcal{X}_h^{\ph,3} \times \mathcal{M}_h^2$ & AP, $\mathcal{X}_h^{\ph,3} \times \mathcal{M}_h^2$  \\
        \includegraphics[width=0.3\linewidth]{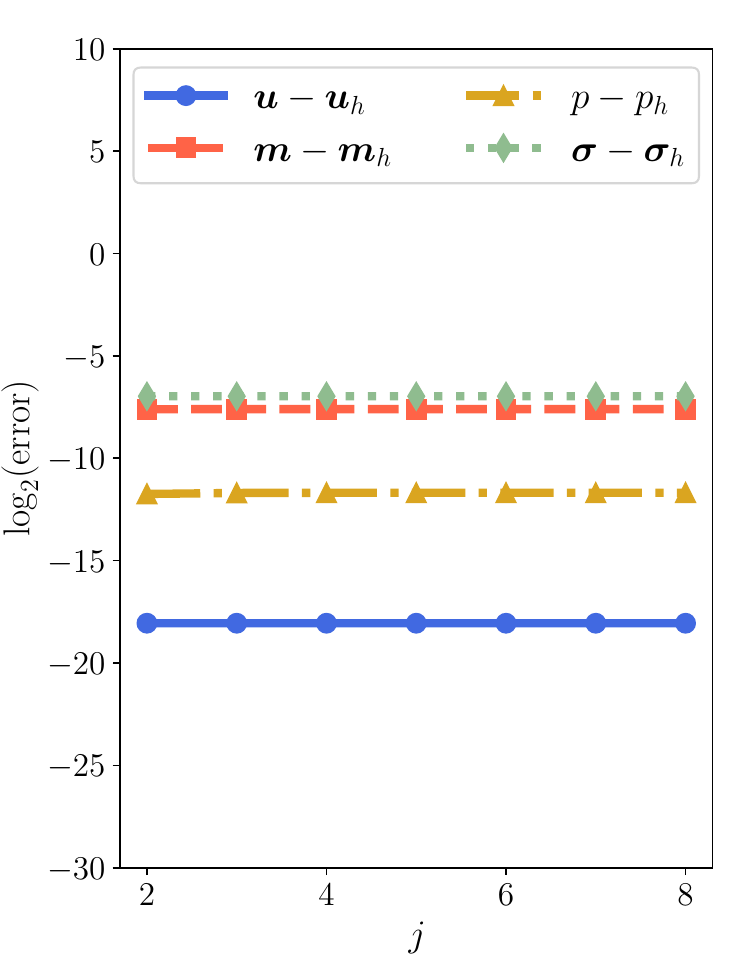} & \includegraphics[width=0.3\linewidth]{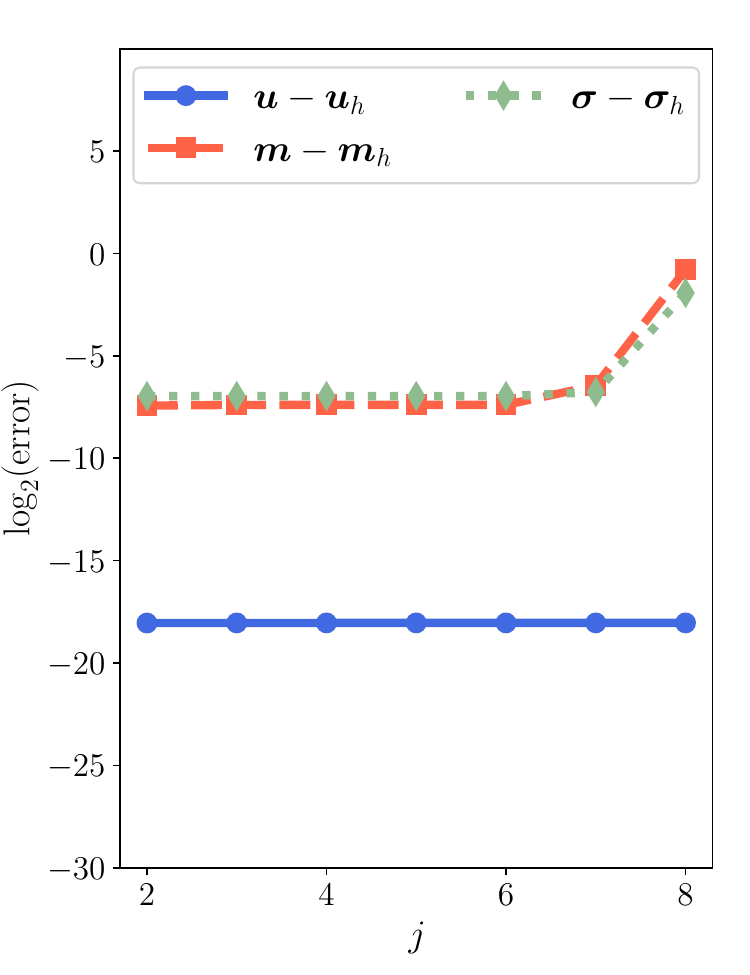} & \includegraphics[width=0.3\linewidth]{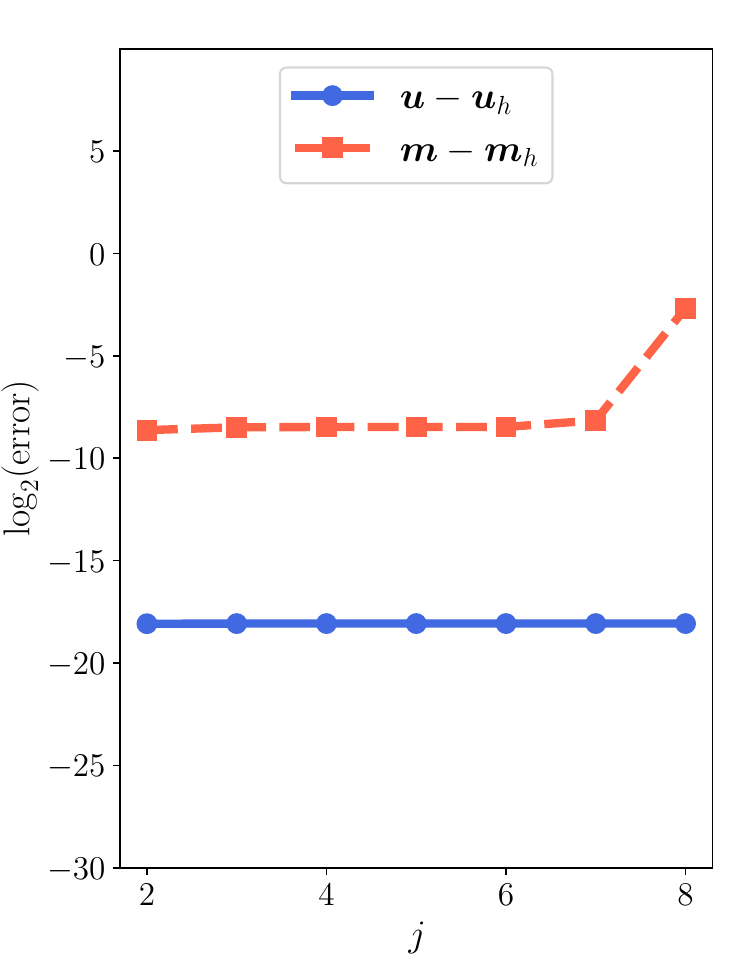}
        \\
        HDP, $\mathcal{X}_h^4 \times \mathcal{M}_h^3 \times \mathcal{P}_h^3$ & PH, $\mathcal{X}_h^{\ph,4} \times \mathcal{M}_h^3$ & AP, $\mathcal{X}_h^{\ph,4} \times \mathcal{M}_h^3$  \\
        \includegraphics[width=0.3\linewidth]{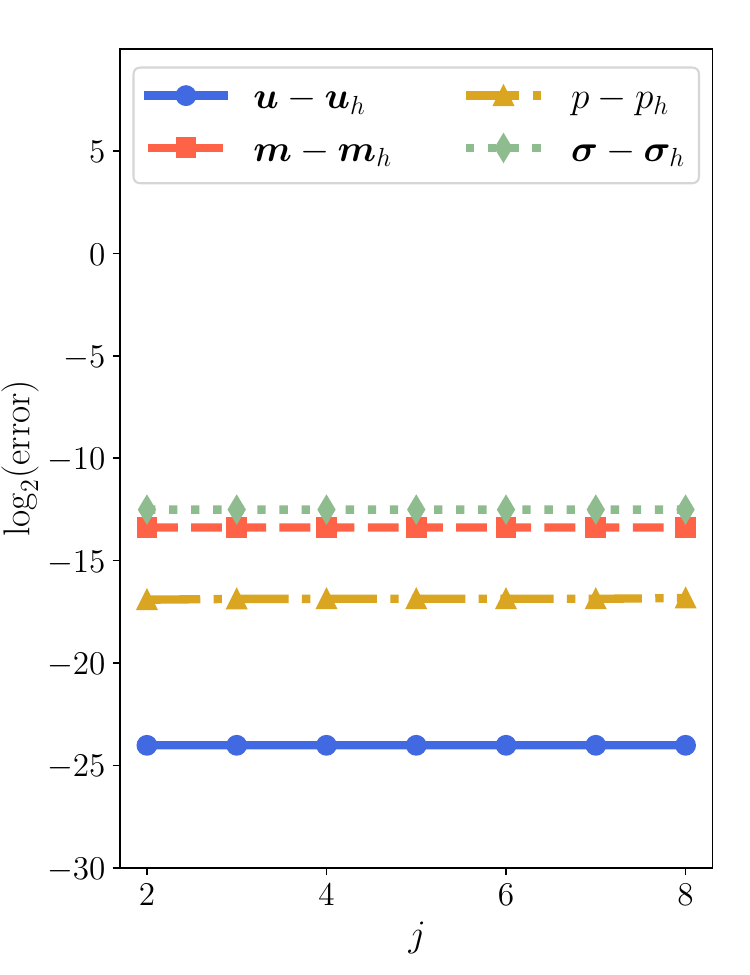} & \includegraphics[width=0.3\linewidth]{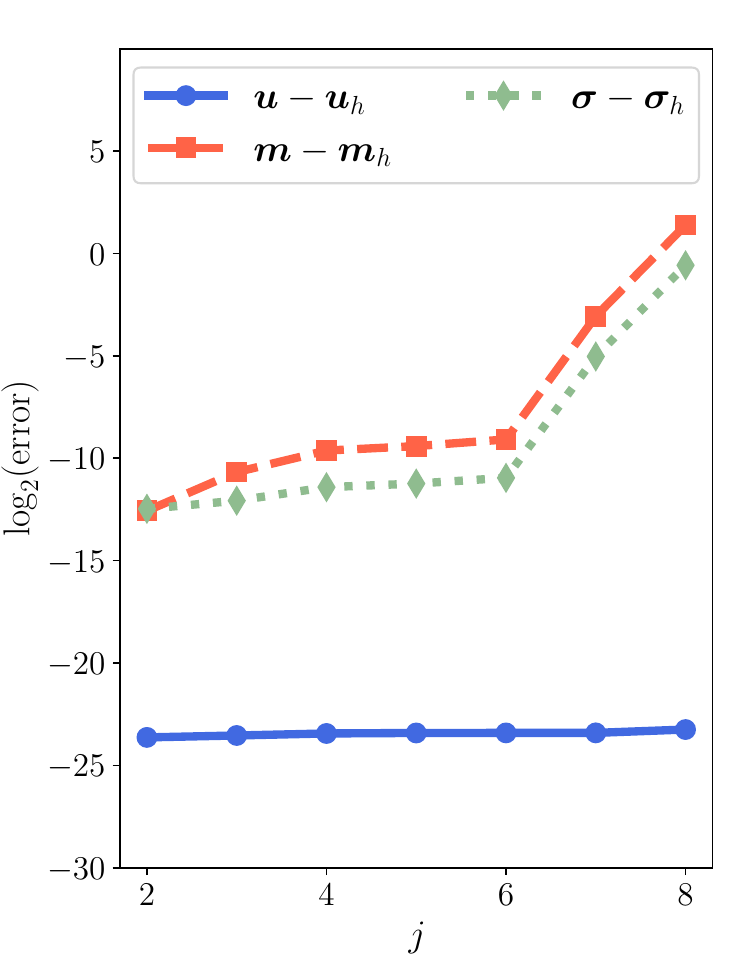} & \includegraphics[width=0.3\linewidth]{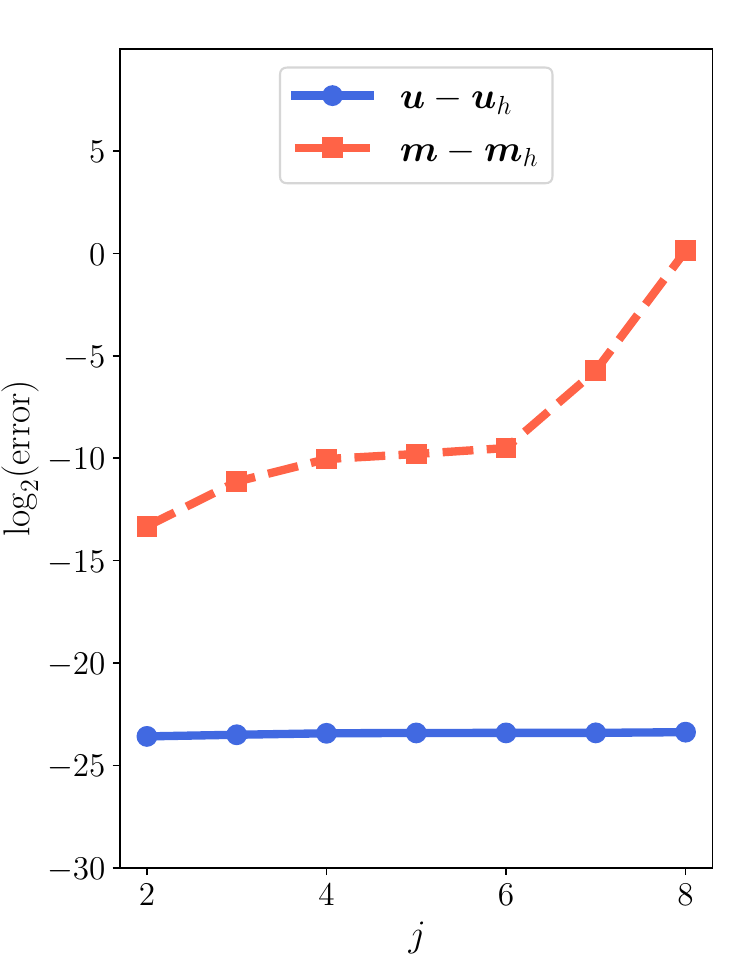}
    \end{tabular}
    \caption{{\color{ok}
    Second test;} results for the locking test on a triangular mesh. The $x$-axis contains the values for $j$ and the $y$-axis the base two log of the approximation errors.
    }
    \label{fig:locking_triang}
\end{figure}

\begin{figure}[htb]
    \centering
    \begin{tabular}{c c c}
        HDP, $\mathcal{X}_h^2 \times \mathcal{M}_h^1 \times \mathcal{P}_h^1$ & PH, $\mathcal{X}_h^{\ph,2} \times \mathcal{M}_h^1$ & AP, $\mathcal{X}_h^{\ph,2} \times \mathcal{M}_h^1$  \\
        \includegraphics[width=0.3\linewidth]{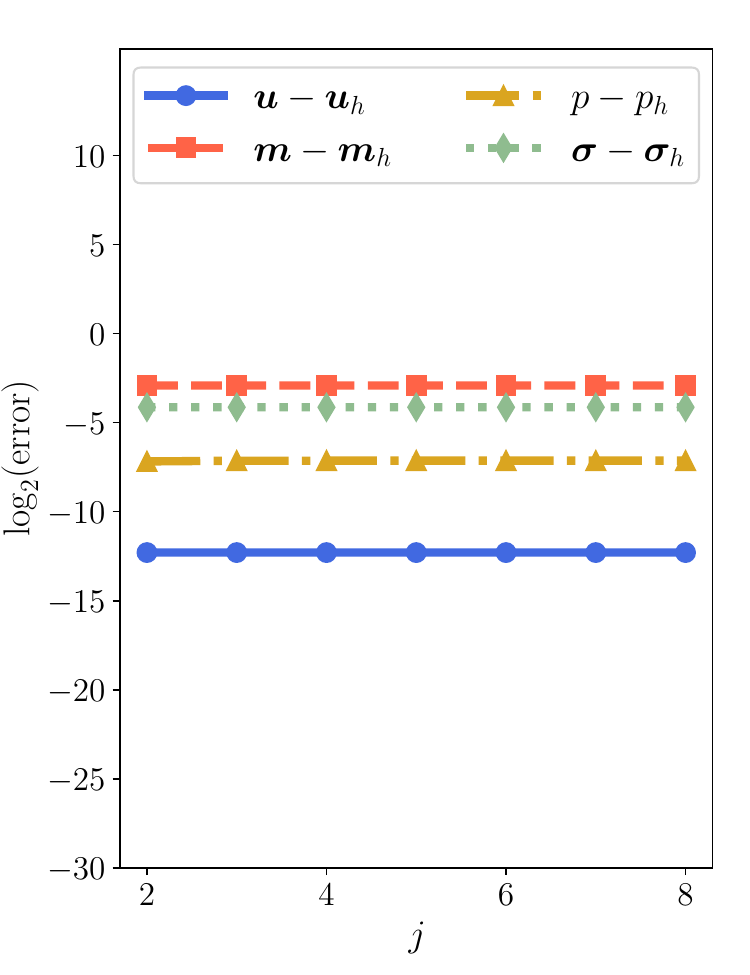} & \includegraphics[width=0.3\linewidth]{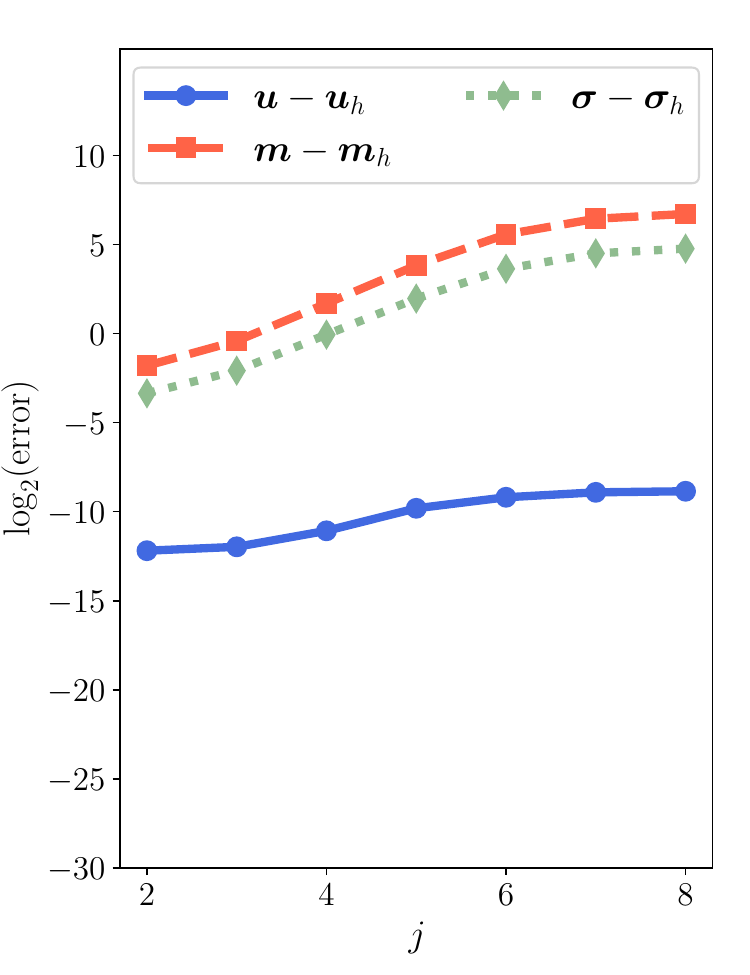} & \includegraphics[width=0.3\linewidth]{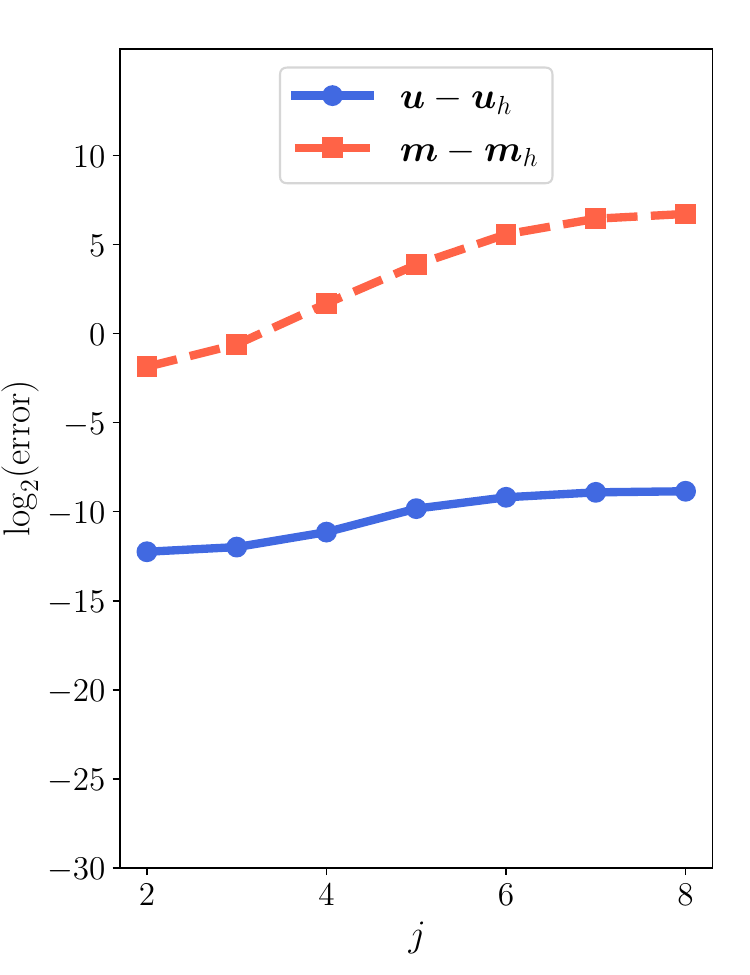}
        \\
        HDP, $\mathcal{X}_h^3 \times \mathcal{M}_h^2 \times \mathcal{P}_h^2$ & PH, $\mathcal{X}_h^{\ph,3} \times \mathcal{M}_h^2$ & AP, $\mathcal{X}_h^{\ph,3} \times \mathcal{M}_h^2$  \\
        \includegraphics[width=0.3\linewidth]{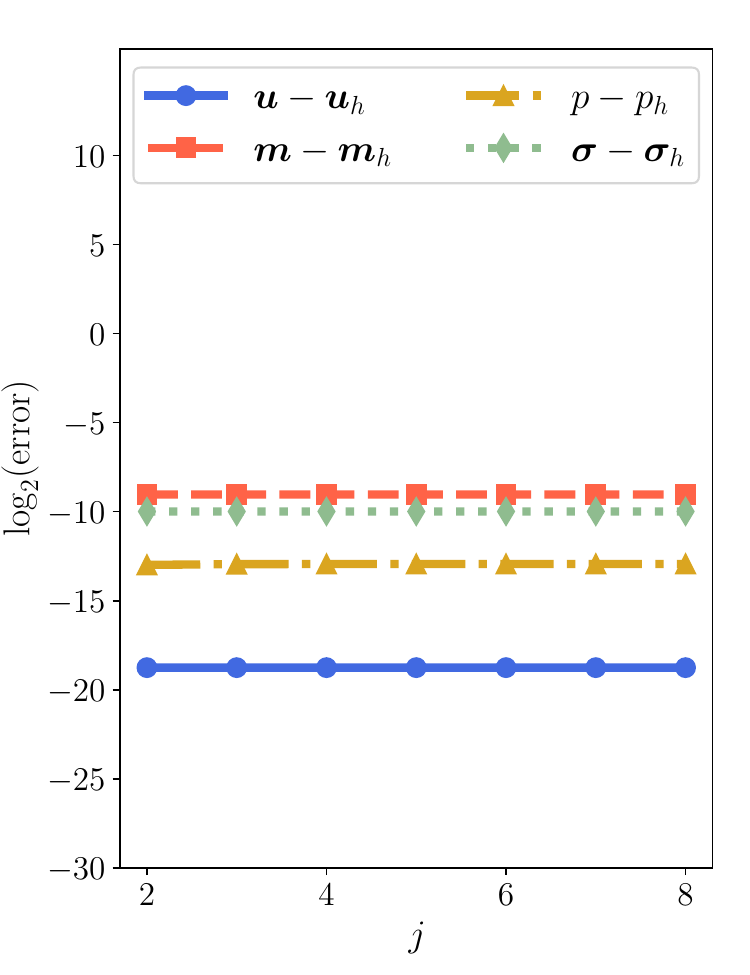} & \includegraphics[width=0.3\linewidth]{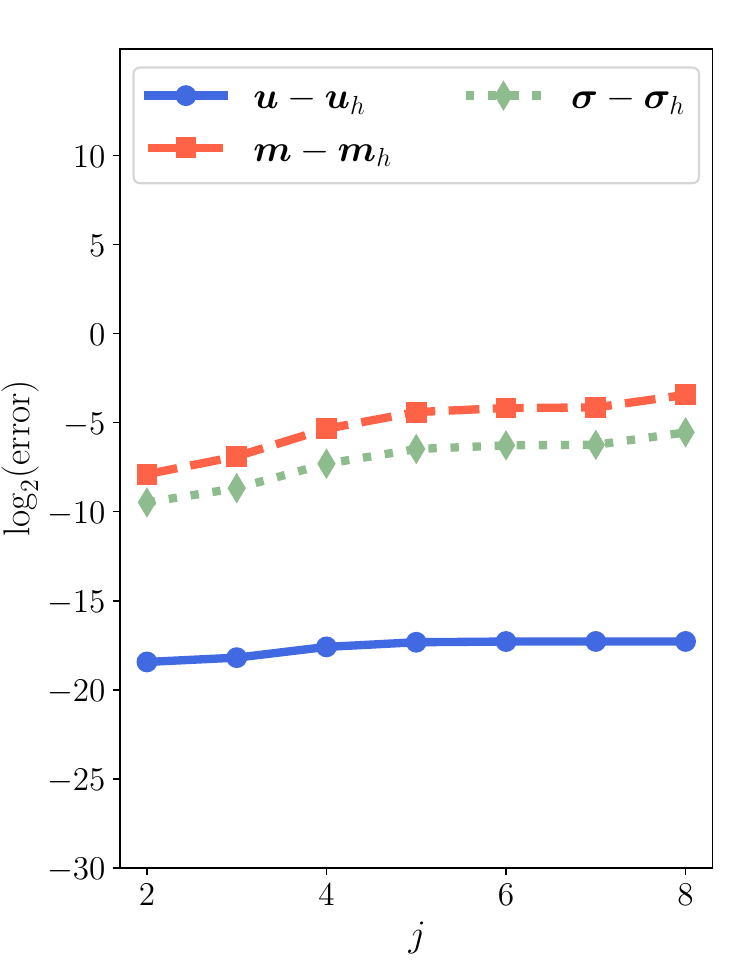} & \includegraphics[width=0.3\linewidth]{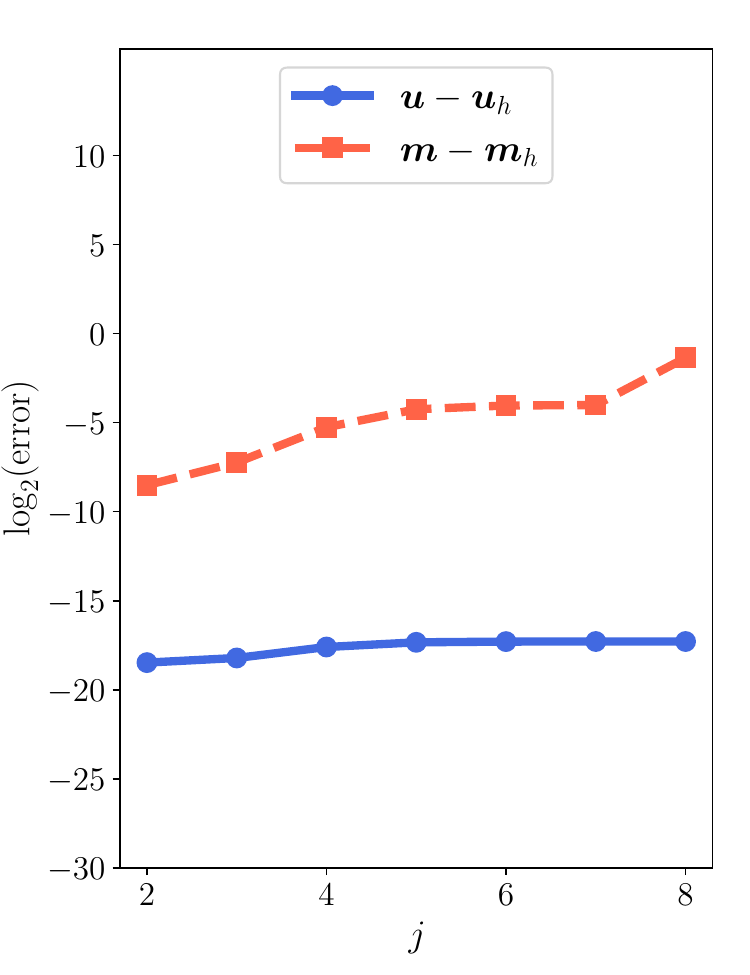}
        \\
        HDP, $\mathcal{X}_h^4 \times \mathcal{M}_h^3 \times \mathcal{P}_h^3$ & PH, $\mathcal{X}_h^{\ph,4} \times \mathcal{M}_h^3$ & AP, $\mathcal{X}_h^{\ph,4} \times \mathcal{M}_h^3$  \\
        \includegraphics[width=0.3\linewidth]{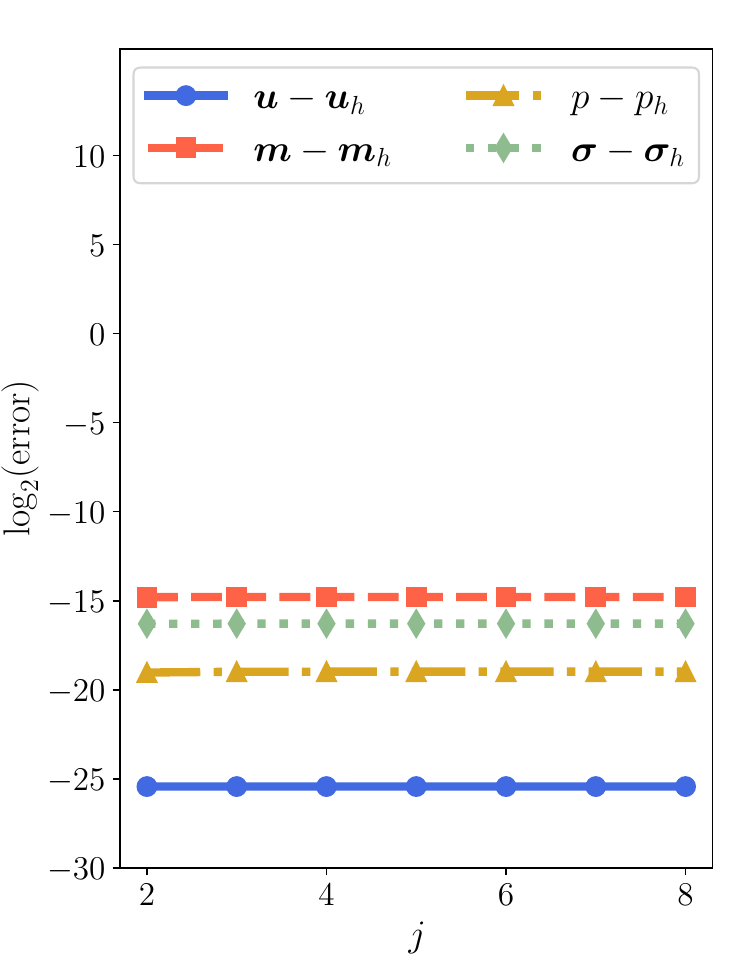} & \includegraphics[width=0.3\linewidth]{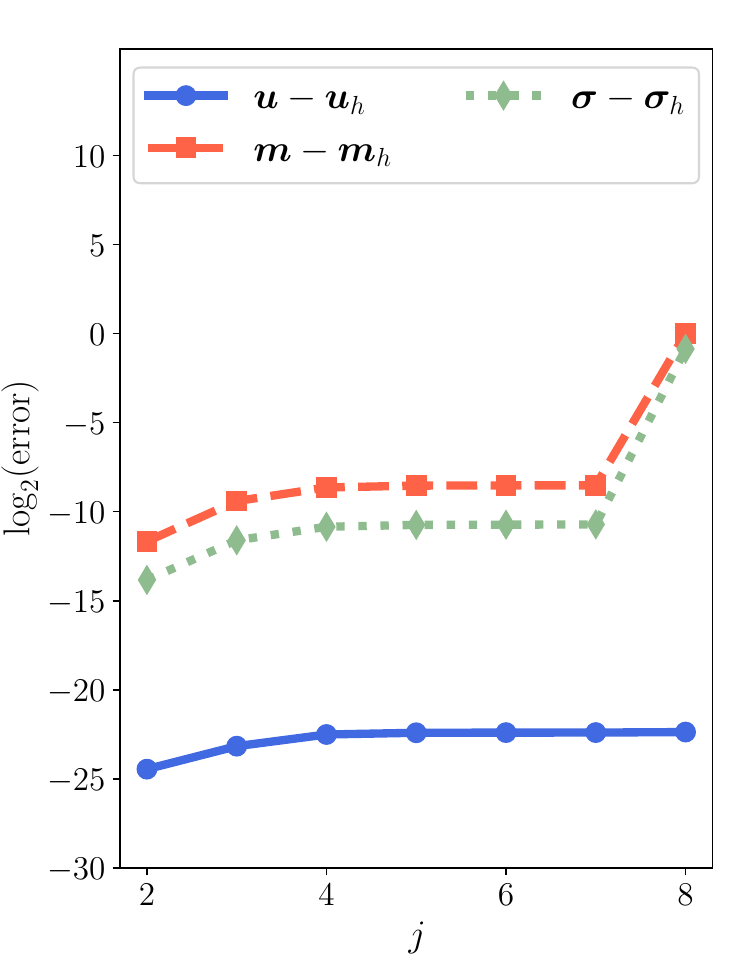} & \includegraphics[width=0.3\linewidth]{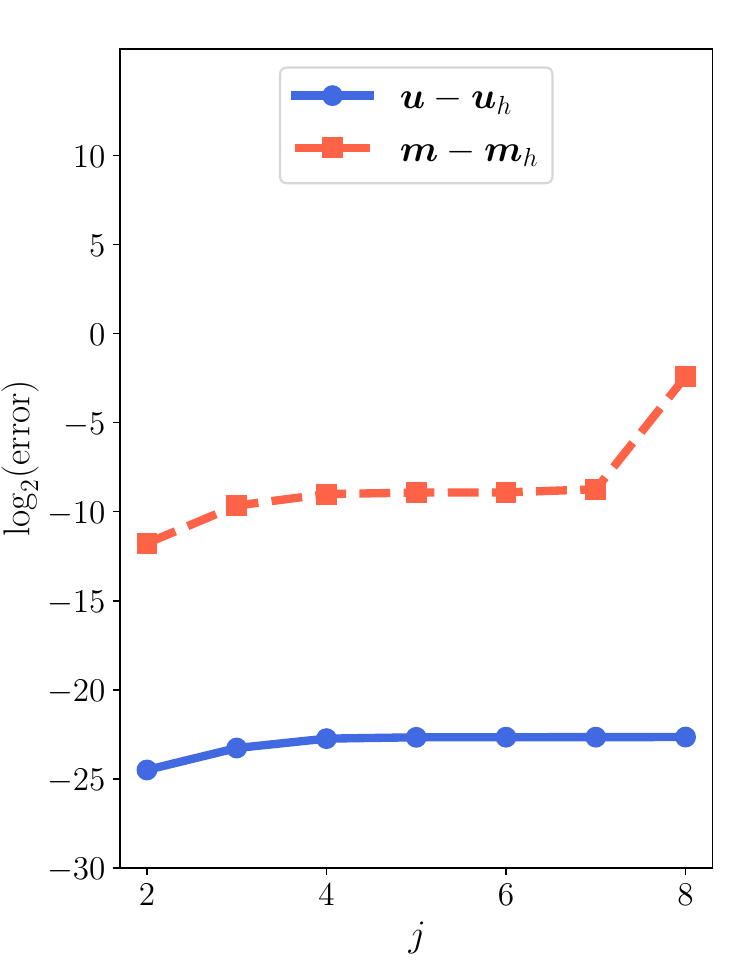}
    \end{tabular}
    \caption{{\color{ok} Second test;} results for the locking test on trapezoidal meshes. The $x$-axis contains the values for $j$ and the $y$-axis the base two log of the approximation errors.}
    \label{fig:locking_trap}
\end{figure}

In this experiment, following the ideas of \cite{taraschi2024global}, we solve the test problem with $\nu = 0.5 - 10^{-j}$, where $j$ goes from $2$ to $8$.
Notice that as $j$ increases, $\nu$ goes to $0.5$, and the problem approaches the incompressibility limit.
In Figures \ref{fig:locking_triang} and \ref{fig:locking_trap}, we plot how the approximation errors vary as a function of $j = -\log_{10}(0.5 - \nu)$ for the triangular and trapezoidal meshes, respectively.

If a given method is robust to locking, we expect its approximation errors for our test problem to stay roughly the same as $j$ increases.
In other words, the plot $j \times \log_2(\mathrm{error})$ should be close to a constant line.
This behavior is verified for the HDP method with the meshes and approximation spaces tested, showing that the proposed method and stress post-processing furnish good approximations even in the more challenging scenario of nearly incompressible problems.
Such a result is in line with the theoretical predictions.

For the standard PH method, the results are more subtle.
On triangular meshes, the displacement approximation seems robust towards locking, but the errors for the Lagrange multiplier increase as $\nu$ approaches $0.5$.
For $r=2$, however, we remark that the multiplier approximation remains robust for moderate incompressible problems.
The results are worse for the trapezoidal meshes, with both the displacement and multiplier approximation errors increasing on nearly incompressible problems.
In both cases, we notice that the Lagrange multiplier is more severely affected by locking, directly impacting the quality of the recovered stress.

\begin{figure}[htb]
    \centering
    \begin{tabular}{c c} AP, Triangles & AP, Trapezoids  \\
        \includegraphics[width=0.3\linewidth]{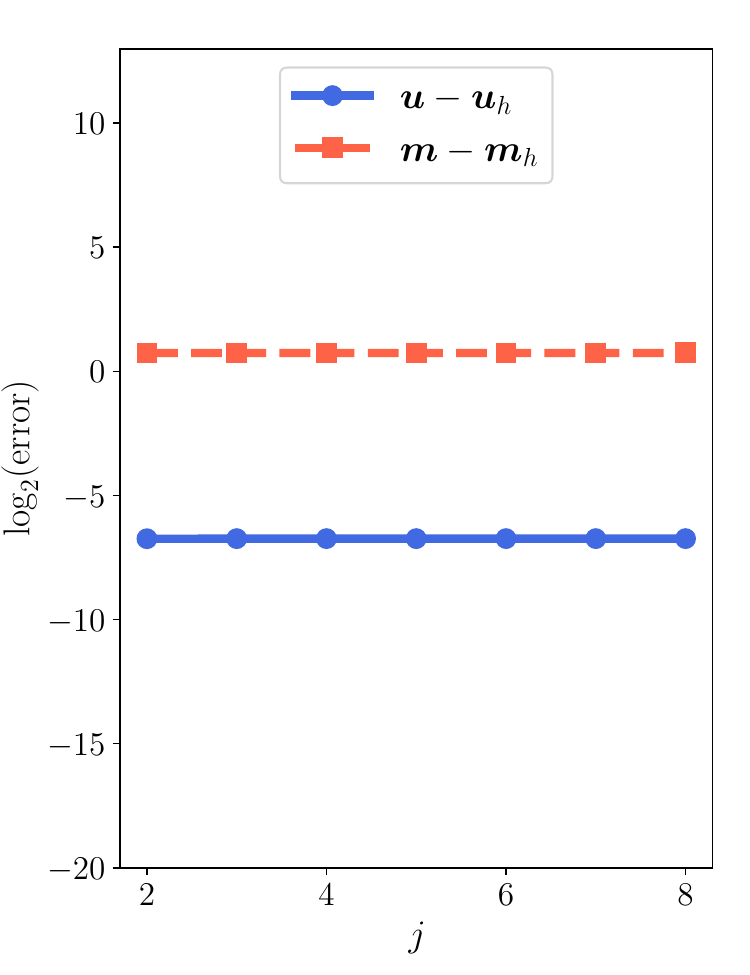} & \includegraphics[width=0.3\linewidth]{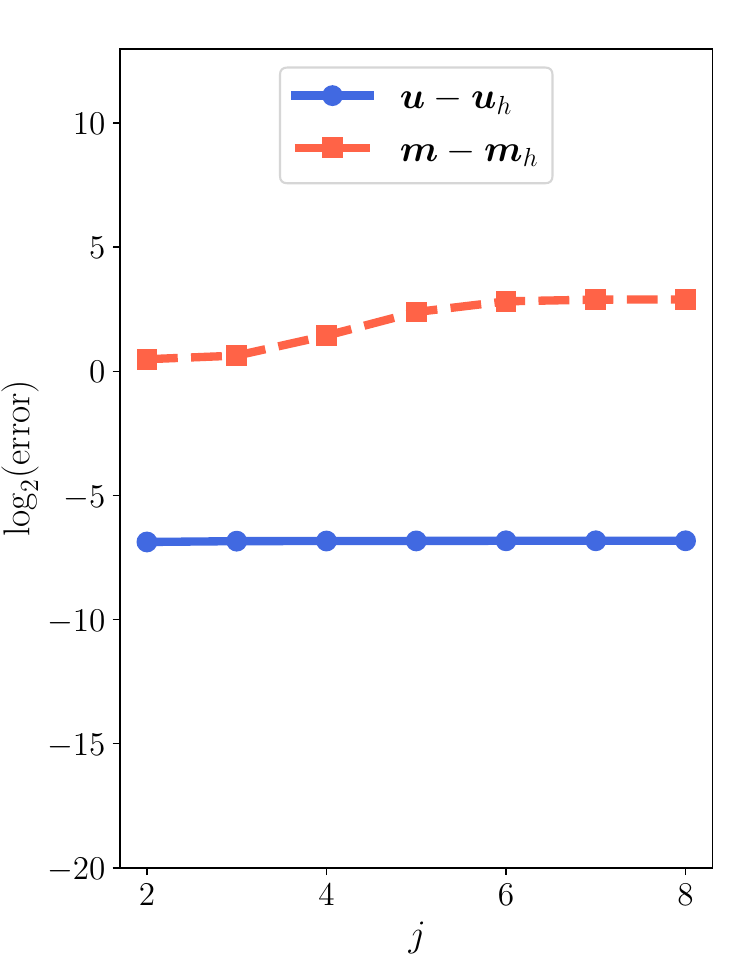}
    \end{tabular}
    \caption{{\color{ok}Second test;} results for the locking test by using the AP method with the lowest-order space $\mathcal{X}_h^{\ph,1} \times \mathcal{M}_h^0$ on triangular and trapezoidal meshes. The $x$-axis contains the values for $j$ and the $y$-axis the base two log of the approximation errors.}
    \label{fig:locking_aph}
\end{figure}

Finally, we verify that the AP method behaves similarly to the PH method when the spaces $\mathcal{X}_h^{\ph,2} \times \mathcal{M}_h^1$ and $\mathcal{X}_h^{\ph,3} \times \mathcal{M}_h^2$ are adopted. This method, however, admits the use of the lowest-order space $\mathcal{X}_h^{\ph,1} \times \mathcal{M}_h^{0}$ \citep{acharya2022primal}. The results for this space, shown in Figure \ref{fig:locking_aph}, indicate robustness near incompressibility when the triangular mesh is adopted.
This result aligns with the analysis and numerical experiments of \cite{acharya2022primal}.
The lowest-order space on the trapezoidal mesh also appears robust for the displacement approximation, but the errors $\vector{m} - \vector{m}_h$ increase.

\begin{figure}[htb]
    \centering
    \includegraphics[width=0.45\linewidth]{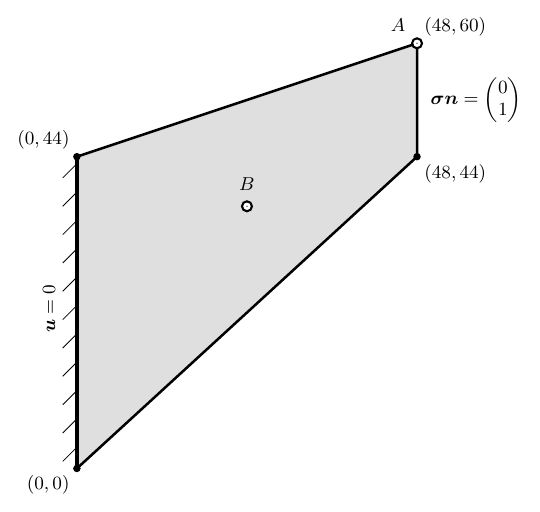}
    \includegraphics[width=0.45\linewidth]{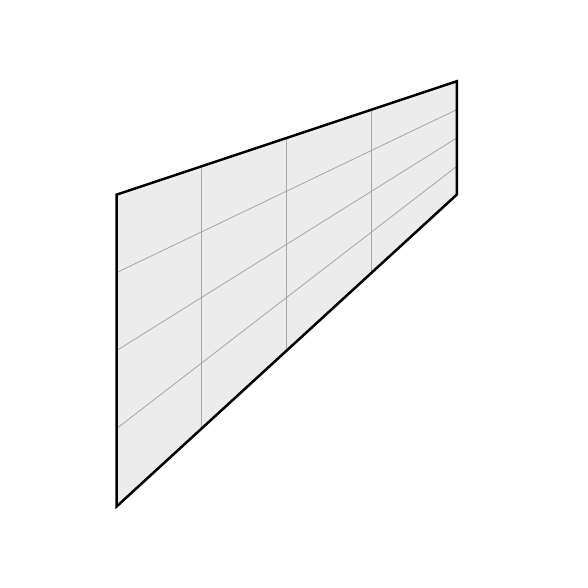}
    \caption{\color{ok}
Domain and boundary conditions representation for the Cook's problem (left) and a $4 \times 4$ quadrilateral mesh (right).}
    \label{fig:cook}
\end{figure}

 \clearpage

{\color{ok}
\subsection{Cook's membrane benchmark}

For the final block of experiments, we continue exploring the behavior of the proposed HDP method in the incompressibility limit, this time using Cook's membrane benchmark problem \cite{cook2007concepts, elguedj2008b}.
The domain is a trapezoid defined by the corners $(0,0)$, $(48,44)$, $(48,60)$, and $(0,44)$, as illustrated in Figure \ref{fig:cook}.
We consider a linear elastic material under plane strain conditions with $E = 70$.
For the Poisson ratio, we consider a compressible case, $\nu = 1/3$, and a nearly-incompressible case, $\nu = 0.5 - 10^{-8}$.
For the boundary conditions, the left side is clamped ($\vector{u} = 0$), while a traction $\tensor{\sigma}\vector{n} = (0,1)^t$ is imposed on the right side.
On the remaining sides, we impose zero traction.


}

{\color{ok}
The domain is partitioned into $n \times n$ squares, as shown in Figure \ref{fig:cook}. The problem is solved using the PH and HDP methods with approximation spaces $\mathcal{X}_h^2 \times \mathcal{M}_h^1$ and $\mathcal{X}_h^2 \times \mathcal{M}_h^1 \times \mathcal{P}_h^1$, respectively.
We measure the vertical displacement in the top-right corner (point A in Figure \ref{fig:cook}) and the  
normal stress (out-of plane) $\sigma_z = \nu({\sigma}_{1,1} + {\sigma}_{2,2})$ 
in the internal point $\vector{x}_\mathrm{B} = (24,37)$ (point B in Figure \ref{fig:cook}) under different mesh refinements.
The solutions are compared with a reference solution obtained on a fine mesh ($n=128$) using a classical mixed displacement-pressure formulation with piecewise cubic displacements and a global stress-recovery strategy, as described in \cite{taraschi2024global}.
Figure \ref{fig:cook_results} shows the obtained results.

\begin{figure}[htb]
    \centering
    \begin{tabular}{c c}
     Displacement, compressible & Stress, compressible \\
     \includegraphics[width=0.45\linewidth]{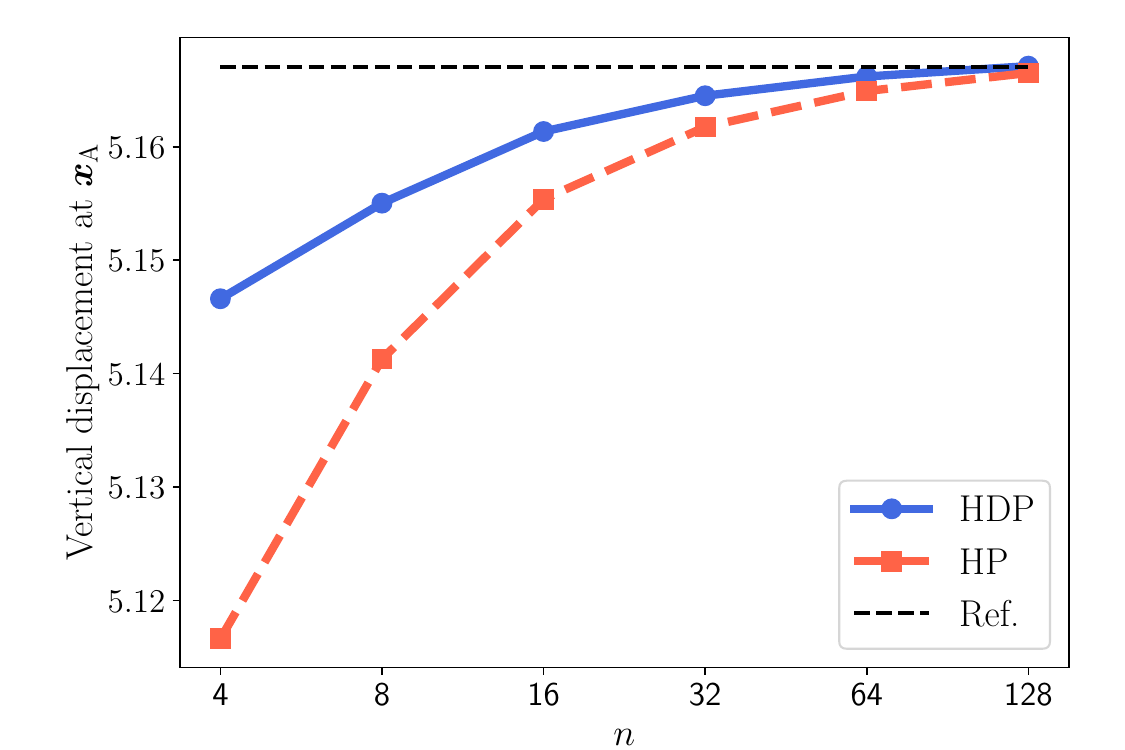} & \includegraphics[width=0.45\linewidth]{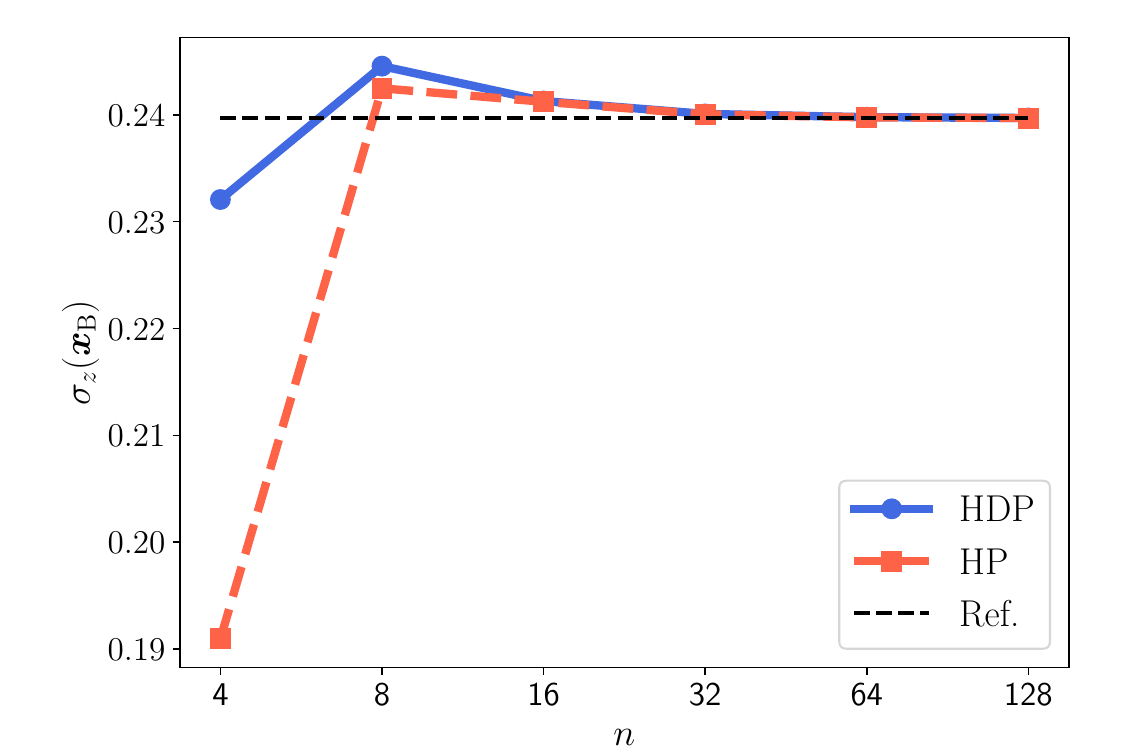} \\
     Displacement, nearly-incompressible & Stress, nearly-incompressible \\
     \includegraphics[width=0.45\linewidth]{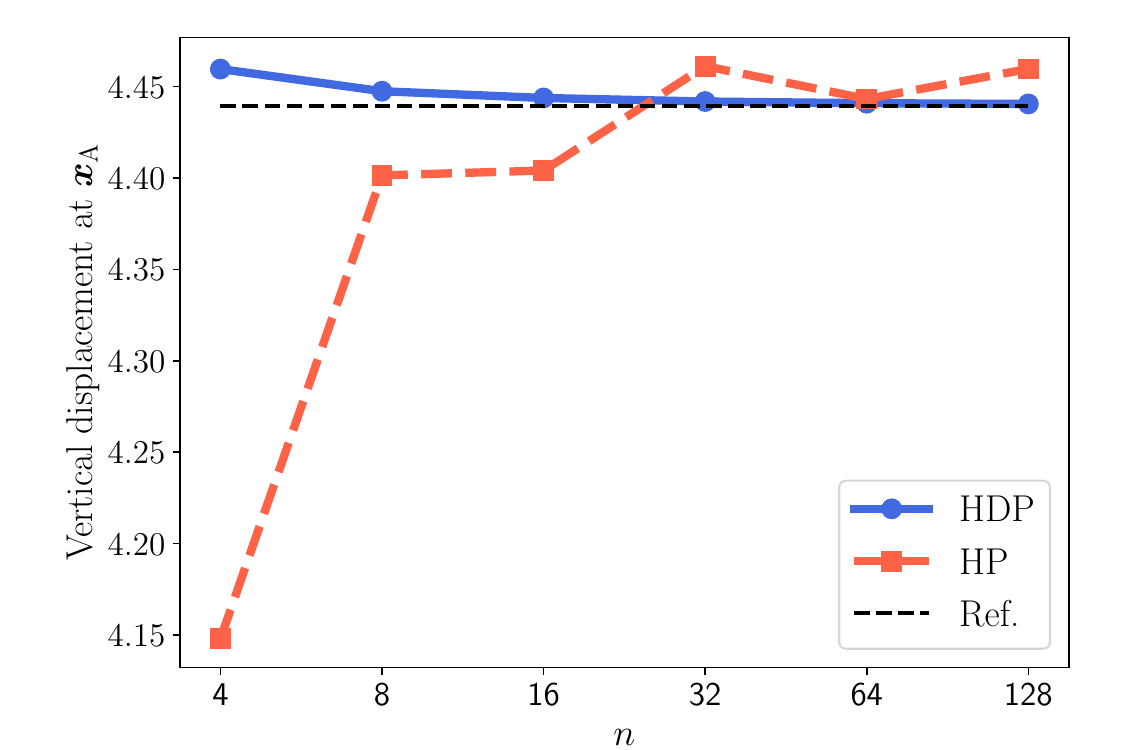} &
     \includegraphics[width=0.45\linewidth]{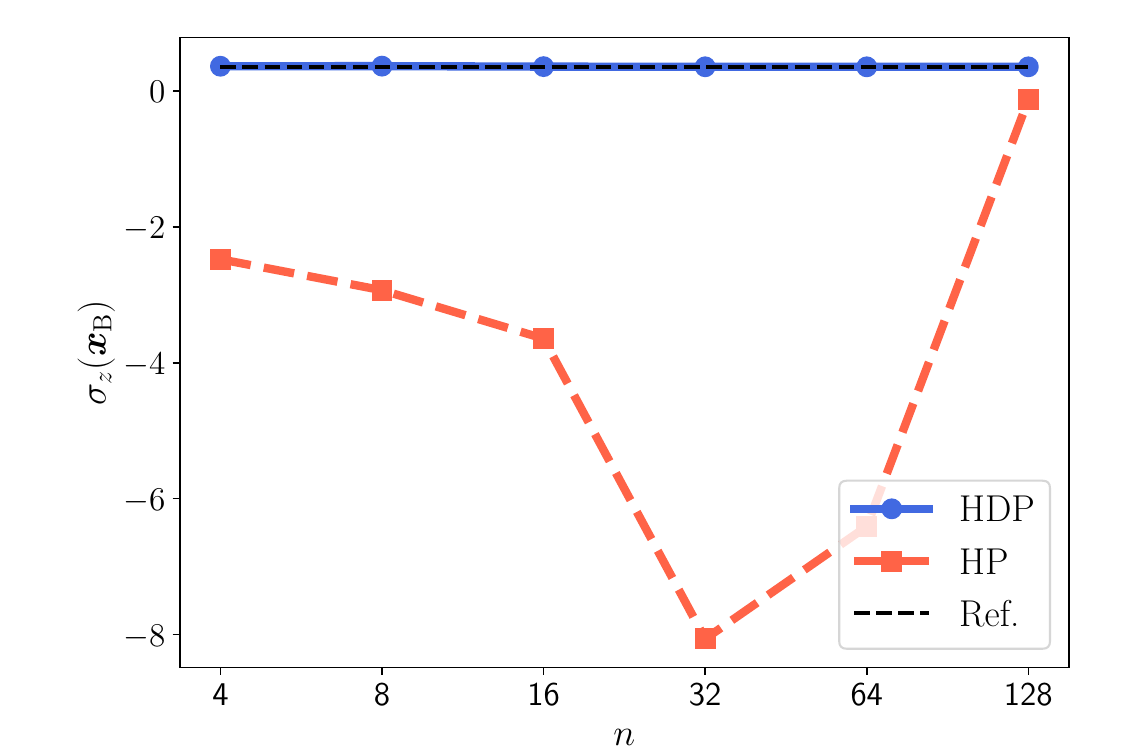} 
    \end{tabular}
    \caption{Displacement and stress results in the compressible in nearly-incompressible cases for the Cook's problem test.}
    \label{fig:cook_results}
\end{figure}

As shown, in the compressible case ($\nu = 1/3$), the PH and HDP results are similar and converge to the reference solution.
For the nearly incompressible case ($\nu = 0.5-10^{-8}$), the HDP approximations remain robust and converge to the reference solution.
The PH displacement now requires finer meshes to approach the reference solution, but the overall result is still acceptable.
The stress, however, remains far from the reference solution.
This result is in line with the results of Section \ref{sec:exp_locking}, in which we verified that locking affects the multiplier (and therefore the post-processed stress) more intensely than the displacement in the PH method.

}

%
%
\subsection*{Acknowledgements}
%
%

The authors thankfully acknowledge the financial support from the research agencies.
This work was funded by the National Council for Scientific and Technological Development - CNPq under grants 403673/2025-9, 307679/2023-3, 304192/2019-8 and 140400/2021-4, the São Paulo Research Foundation (FAPESP) under grant 2013/07375-0, and the Center for Energy and Petroleum Studies (CEPETRO).
%

%
%
\subsection*{Data availability statement}\label{sec:data}
%
%

All data that support the findings of this study are available upon reasonable request.

%
%
\subsection*{Declaration of competing interest}\label{sec:declaration}
%
%

The authors declare that they have no known competing financial interests or personal relationships that could have appeared to influence the work reported in this paper.

\FloatBarrier  
\bibliography{ref}
\bibliographystyle{unsrt}

\clearpage

\appendix

\vspace{1cm}

\begin{center}
{\sc Appendix}
\end{center}

\section{Static condensation} \label{ap:condensation}

Denote by $A$, $B$, and $C$ the local (at the element level) stiffness matrices associated with the bilinear forms $a(\cdot, \cdot)$, $b(\cdot, \cdot)$, and $c(\cdot, \cdot)$, respectively.
Furthermore, we denote by $D$ the local stiffness matrix associated with the term $-\frac{1}{\lambda} (p_h, q)_{\Omega}$ in \eqref{eq:method_HDP_c} and by $F$ and $G$ the local source terms related to the right-hand sides of \eqref{eq:method_HDP_a} and \eqref{eq:method_HDP_b}, respectively.
Notice that $F$ depends on the source term $\vector{f}$ and the Neumann boundary condition $\vector{t}_N$, while $G$ comes from the Dirichlet boundary condition $\vector{u}_D$.
With these definitions, the local stiffness system for the method \eqref{eq:method_HDP} can be written as
%
\begin{equation} \label{eq:system_full}
    \begin{bmatrix}
        A & B^t & C^t \\
        B & 0   & 0   \\
        C & 0   & D   
    \end{bmatrix}
    \begin{bmatrix}
        u_{\mathrm{dof}} \\
        m_{\mathrm{dof}} \\
        p_{\mathrm{dof}}
    \end{bmatrix}
    =
    \begin{bmatrix}
        F \\
        G \\
        0
    \end{bmatrix},
\end{equation} 
where $u_{\mathrm{dof}}$, $m_{\mathrm{dof}}$, and $p_{\mathrm{dof}}$ are the local degrees of freedom associated with the approximated displacement, multiplier, and pressure fields.

Now, assume that $\mathcal{X}_{\mathrm{rm}}^{K} $ is spanned by the first functions of the local basis for $\mathcal{X}_K$.
Setting $k_\mathrm{rm} =\textrm{dim} \mathcal{X}_{\mathrm{rm}}^{K} $ and $\tilde{k} = \textrm{dim} \mathcal{X}_{K} - \textrm{dim} \mathcal{X}_{\mathrm{rm}}^{K} $, consider the following block decomposition
\[
A = \begin{bmatrix}
    0 & 0       \\
    0 & A_{2,2} \\
\end{bmatrix},
\quad
B = \begin{bmatrix}
    B_1 & B_2
\end{bmatrix},
\quad
C = \begin{bmatrix}
    C_1 & C_2
\end{bmatrix}
\quad
F = \begin{bmatrix}
    F_1 \\
    F_2
\end{bmatrix},
\quad
u_{\mathrm{dof}} = \begin{bmatrix}
    u_{\mathrm{dof}}^\mathrm{rm} \\
    \tilde{u}_{\mathrm{dof}}
\end{bmatrix},
\]
where $A_{2,2} \in \mathbb{R}^{\tilde{k} \times \tilde{k}}$; $B_{1} \in \mathbb{R}^{\textrm{dim} \mathcal{M}_K \times k_\mathrm{rm}}$; $B_{2} \in \mathbb{R}^{\textrm{dim} \mathcal{M}_K \times \tilde{k}}$; $C_{1} \in \mathbb{R}^{\textrm{dim} \mathcal{P}_K \times k_\mathrm{rm}}$; $C_{2} \in \mathbb{R}^{\textrm{dim} \mathcal{P}_K \times \tilde{k}}$; $F_{1}, u^\mathrm{rm}_{\mathrm{dof}} \in \mathbb{R}^{k_\mathrm{rm}}$ and $F_{2}, \tilde{u}_{\mathrm{dof}} \in \mathbb{R}^{\tilde{k}}$.
It is possible to show that matrix $A_{2,2}$ is inversible, which is equivalent to show that $a(\cdot, \cdot)$ is coercive on $\mathcal{X}_K \backslash \mathcal{X}_{\mathrm{rm}}^{K} $.
In fact, this can be seen as the discrete counterpart of the solvability of the local sub-problems defined in \cite{gomes2024low}.

Next, by inverting $A_{2,2}$, we can rewrite the system \eqref{eq:system_full} in terms of $u_{\mathrm{dof}}^\mathrm{rm}$, $m_{\mathrm{dof}}$, and $p_{\mathrm{dof}}$ only
\begin{equation} \label{eq:system_semi}
    \begin{bmatrix}
        0 & B_1^t & C_1^t \\
        B_1 & E_{1,1} & E_{2,1}^t \\
        C_1 & E_{2,1} & E_{2,2}
    \end{bmatrix}
    \begin{bmatrix}
        u_{\mathrm{dof}}^\mathrm{rm} \\
        m_{\mathrm{dof}} \\
        p_{\mathrm{dof}}
    \end{bmatrix}
    =
    \begin{bmatrix}
        F_1 \\
        H_1 \\
        H_2
    \end{bmatrix},
\end{equation}
where
$E_{1,1} = -B_2 A_{2,2}^{-1} B_2^t$,  $E_{2,1} = -C_2 A_{2,2}^{-1} B_2^t$, $E_{2,2} = D - C_2 A_{2,2}^{-1} C_2^t$, $H_1 = G - B_2 A_{2,2}^{-1} F_2$, and  $H_2 = -C_2 A_{2,2}^{-1}F_2. $
The remaining degrees of freedom $\tilde{u}_{\mathrm{dof}}$ can be recovered by
\begin{equation} \label{eq:recover1}
    \tilde{u}_{\mathrm{dof}} = A_{2,2}^{-1} (F_2 - B_2^t m_{\mathrm{dof}} - C_2^t p_{\mathrm{dof}}).
\end{equation}
Since the functions in $\mathcal{X}_h$ are discontinuous across inter-element boundaries, such a procedure allows for the static condensation of the majority of the degrees of freedom associated with $\vector{u}_h$, greatly reducing the computational cost of the HDP method.
{\color{ok}
A further condensation, associated with the pressure approximation, is possible by inverting the matrix $E_{2,2}$, leading to the final, fully-condensed system
\begin{equation} \label{eq:HDP_system_cond}
\begin{bmatrix}
R_{1,1} & R_{2,1}^t \\
R_{2,1} & R_{2,2} \\
\end{bmatrix}
\begin{bmatrix}
u_{\mathrm{dof}}^\mathrm{rm} \\
m_{\mathrm{dof}}
\end{bmatrix}
=
\begin{bmatrix}
J_1 \\
J_2
\end{bmatrix},
\end{equation}
with $R_{1,1} = -C_1^t E_{2,2}^{-1} C_1$, $R_{2,1} = B_1 - E_{2,1}^t E_{2,2}^{-1} C_1$, $R_{2,2} = E_{1,1} - E_{2,1}^t E_{2,2}^{-1} E_{2,1}$, $J_1 = F_1 - C_1^t E_{2,2}^{-1} H_2$, and $J_2 = H_1 - E_{2,1}^t E_{2,2}^{-1} H_2.$
Since $E_{2,2}$ depends on $\lambda$, its inversion could be problematic in the nearly-incompressible scenario due to floating-point arithmetic.
In our testing, however, using double  variables and standard inversion routines from LAPACK, this did not pose any problems

Finally, to highlight the benefits of static condensation, Table \ref{tab:condensation} compares the total number of degrees of freedom in the uncondensed and condensed global systems for the meshes used in Section \ref{sec:exp_convergence}.
Also, using the same experimental setup as Section \ref{sec:exp_convergence}, Table \ref{tab:cpu_time} reports the execution CPU times across various stages of the HDP implementation. 
For this, we fully implemented the HDP method in in-house Fortran code using direct solvers from LAPACK.
Computation times were recorded on a personal computer using the native subroutine \textit{cpu\_time} in single-threaded mode.
The evaluated stages include assembly with static condensation, global system solution, solution de-condensation, and the stress post-processing described in Section \ref{sec:post}.
Notably, all of these steps can be parallelized, so the elapsed times could differ significantly in a multithreaded execution.
In particular, de-condensation and post-processing are trivially parallelized and would therefore gain substantial speedup from parallel computing.

\begin{table}[htb]
\centering
\caption{\color{ok} First test; number of degrees of freedom of the HDP method before and after static condensation.}
\begin{tabular}{c|c c c | c c c}
\hline
& \multicolumn{3}{c|}{Triangular meshes} & \multicolumn{3}{c}{Quadrilateral meshes} \\
$n$ & $\mathrm{dof}$ & $\mathrm{dof}^C$ & red. (\%) & $\mathrm{dof}$ & $\mathrm{dof}^C$ & red. (\%) \\
\hline
\multicolumn{7}{c}{Space $\mathcal{X}_h^{2} \times \mathcal{M}_h^{1} \times \mathcal{P}_h^{1}$} \Tstrut \Bstrut \\
\hline
8   & 3264   & 1216   & 37.25 & 2048   & 768    & 37.50 \\
16  & 12928  & 4736   & 36.63 & 8064   & 2944   & 36.51 \\
32  & 51456  & 18688  & 36.32 & 32000  & 11520  & 36.00 \\
64  & 205312 & 74240  & 36.16 & 127488 & 45568  & 35.74 \\
128 & 820224 & 295936 & 36.08 & 508928 & 181248 & 35.61 \\
\hline
\multicolumn{7}{c}{Space $\mathcal{X}_h^{3} \times \mathcal{M}_h^{2} \times \mathcal{P}_h^{2}$} \Tstrut \Bstrut \\
\hline
8   & 5088    & 1632   & 32.08 & 3424   & 1056   & 30.84 \\
16  & 20160   & 6336   & 31.43 & 13504  & 4032   & 29.86 \\
32  & 80256   & 24960  & 31.10 & 53632  & 15744  & 29.35 \\
64  & 320256  & 99072  & 30.94 & 213760 & 62208  & 29.10 \\
128 & 1279488 & 394752 & 30.85 & 853504 & 247296 & 28.97 \\
\hline
\multicolumn{7}{c}{Space $\mathcal{X}_h^{4} \times \mathcal{M}_h^{3} \times \mathcal{P}_h^{3}$} \Tstrut \Bstrut \\
\hline
8   & 7808    & 2048   & 26.23 & 5120    & 1344   & 26.25 \\
16  & 30976   & 7936   & 25.62 & 20224   & 5120   & 25.32 \\
32  & 123392  & 31232  & 25.31 & 80384   & 19968  & 24.84 \\
64  & 492544  & 123904 & 25.16 & 320512  & 78848  & 24.60 \\
128 & 1968128 & 493568 & 25.08 & 1280000 & 313344 & 24.48 \\
\hline
\end{tabular}

\label{tab:condensation}
\end{table}

\begin{table}[htb]
    \caption{\color{ok} First test; CPU time for the different steps in the HDP method execution.}
    \centering
    \resizebox{\textwidth}{!}{
    \begin{tabular}{c|c c| c c| c c| c c}
       \hline
       & \multicolumn{2}{c|}{Assembly} & \multicolumn{2}{c|}{Solve} & \multicolumn{2}{c|}{De-condensation} & \multicolumn{2}{c}{Post-processing} \\
       $n$ & CPU time (s) & \% & CPU time (s) & \% & CPU time (s) & \% & CPU time (s) & \% \\
       \hline
       \multicolumn{9}{c}{$\mathcal{X}^2 \times \mathcal{M}^1 \times \mathcal{P}^1$} \Tstrut \Bstrut \\
       \hline
       8   & 0.90E-02 & 35.07 & 0.22E-02 &  8.65 & 0.89E-02 & 35.03 & 0.54E-02 & 21.25 \\
       16  & 0.36E-01 & 29.35 & 0.30E-01 & 24.04 & 0.36E-01 & 29.27 & 0.21E-01 & 17.35 \\
       32  & 0.15E+00 & 17.92 & 0.43E+00 & 53.58 & 0.15E+00 & 17.94 & 0.86E-01 & 10.56 \\
       64  & 0.58E+00 &  7.31 & 0.64E+01 & 80.75 & 0.61E+00 &  7.64 & 0.34E+00 &  4.30 \\
       128 & 0.23E+01 &  2.03 & 0.11E+03 & 94.46 & 0.26E+01 &  2.30 & 0.14E+01 &  1.21 \\
       \hline
       \multicolumn{9}{c}{$\mathcal{X}^3 \times \mathcal{M}^2 \times \mathcal{P}^2$} \Tstrut \Bstrut\\
       \hline
       8   & 0.29E-01 & 32.70 & 0.59E-02 &  6.67 & 0.28E-01 & 31.74 & 0.25E-01 & 28.90 \\
       16  & 0.11E+00 & 27.95 & 0.80E-01 & 19.90 & 0.11E+00 & 28.33 & 0.95E-01 & 23.82 \\
       32  & 0.45E+00 & 18.16 & 0.12E+01 & 47.89 & 0.46E+00 & 18.58 & 0.38E+00 & 15.36 \\
       64  & 0.18E+01 &  7.47 & 0.19E+02 & 78.50 & 0.19E+01 &  7.74 & 0.15E+01 &  6.29 \\
       128 & 0.71E+01 &  2.14 & 0.31E+03 & 93.72 & 0.78E+01 &  2.33 & 0.60E+01 &  1.81 \\
       \hline
       \multicolumn{9}{c}{$\mathcal{X}^4 \times \mathcal{M}^3 \times \mathcal{P}^3$}\Tstrut \Bstrut\\
       \hline
       8   & 0.80E-01 & 30.99 & 0.13E-01 &  5.01 & 0.80E-01 & 30.88 & 0.86E-01 & 33.12 \\
       16  & 0.32E+00 & 27.47 & 0.18E+00 & 15.23 & 0.32E+00 & 27.94 & 0.34E+00 & 29.36 \\
       32  & 0.13E+01 & 19.55 & 0.26E+01 & 39.31 & 0.13E+01 & 20.10 & 0.14E+01 & 21.04 \\
       64  & 0.51E+01 &  8.95 & 0.41E+02 & 72.38 & 0.52E+01 &  9.16 & 0.54E+01 &  9.51 \\
       128 & 0.20E+02 &  2.68 & 0.69E+03 & 91.66 & 0.21E+02 &  2.78 & 0.22E+02 &  2.88 \\
       \hline
    \end{tabular}
    }
    \label{tab:cpu_time}
\end{table}

}



\end{document}